\documentclass[reqno]{amsart}

\usepackage{tikz-cd}

\usepackage{fullpage,dsfont}

\usepackage[pdfencoding=auto, psdextra]{hyperref}

\usepackage{parskip}
\makeatletter 
\def\thm@space@setup{%
 \thm@preskip=\parskip \thm@postskip=0pt
}
\def\th@remark{%
  \thm@headfont{\itshape}%
  \normalfont 
  \thm@preskip\parskip \thm@postskip=0pt
}
\makeatother

\usepackage[nobysame,alphabetic,initials,msc-links]{amsrefs}

\usepackage[final]{pdfpages}

\usepackage{multirow}

\usepackage{array}

\DefineSimpleKey{bib}{how}
\DefineSimpleKey{bib}{mrclass}
\DefineSimpleKey{bib}{mrnumber}
\DefineSimpleKey{bib}{fjournal}
\DefineSimpleKey{bib}{mrreviewer}

\renewcommand{\PrintDOI}[1]{%
  \href{http://dx.doi.org/#1}{{\tt DOI:#1}}%
}
\renewcommand{\eprint}[1]{#1}
\BibSpec{book}{%
    +{}  {\PrintPrimary}                {transition}
    +{.} { \PrintDate}                  {date}
    +{.} { \textit}                     {title}
    +{.} { }                            {part}
    +{:} { \textit}                     {subtitle}
    +{,} { \PrintEdition}               {edition}
    +{}  { \PrintEditorsB}              {editor}
    +{,} { \PrintTranslatorsC}          {translator}
    +{,} { \PrintContributions}         {contribution}
    +{,} { }                            {series}
    +{,} { \voltext}                    {volume}
    +{,} { }                            {publisher}
    +{,} { }                            {organization}
    +{,} { }                            {address}
    +{,} { }                            {status}
    +{,} { \PrintDOI}                   {doi}
    +{,} { \PrintISBNs}                 {isbn}
    +{}  { \parenthesize}               {language}
    +{}  { \PrintTranslation}           {translation}
    +{;} { \PrintReprint}               {reprint}
    +{.} { }                            {note}
    +{.} {}                             {transition}
    +{}  {\SentenceSpace \PrintReviews} {review}
}
\BibSpec{article}{%
    +{}  {\PrintAuthors}                {author}
    +{,} { \textit}                     {title}
    +{.} { }                            {part}
    +{:} { \textit}                     {subtitle}
    +{,} { \PrintContributions}         {contribution}
    +{.} { \PrintPartials}              {partial}
    +{,} { }                            {journal}
    +{}  { \textbf}                     {volume}
    +{}  { \PrintDatePV}                {date}
    +{,} { \issuetext}                  {number}
    +{,} { \eprintpages}                {pages}
    +{,} { }                            {status}
    +{,} { \PrintDOI}                   {doi}
    +{,} { \eprint}        {eprint}
    +{}  { \parenthesize}               {language}
    +{}  { \PrintTranslation}           {translation}
    +{;} { \PrintReprint}               {reprint}
    +{.} { }                            {note}
    +{.} {}                             {transition}
    +{}  {\SentenceSpace \PrintReviews} {review}
}
\BibSpec{collection.article}{%
    +{}  {\PrintAuthors}                {author}
    +{,} { \textit}                     {title}
    +{.} { }                            {part}
    +{:} { \textit}                     {subtitle}
    +{,} { \PrintContributions}         {contribution}
    +{,} { \PrintConference}            {conference}
    +{}  {\PrintBook}                   {book}
    +{,} { }                            {booktitle}
    +{,} { \PrintDateB}                 {date}
    +{,} { pp.~}                        {pages}
    +{,} { }                            {publisher}
    +{,} { }                            {organization}
    +{,} { }                            {address}
    +{,} { }                            {status}
    +{,} { \PrintDOI}                   {doi}
    +{,} { \eprint}        {eprint}
    +{}  { \parenthesize}               {language}
    +{}  { \PrintTranslation}           {translation}
    +{;} { \PrintReprint}               {reprint}
    +{.} { }                            {note}
    +{.} {}                             {transition}
    +{}  {\SentenceSpace \PrintReviews} {review}
}
\BibSpec{misc}{%
  +{}{\PrintAuthors}  {author}
  +{,}{ \textit}      {title}
  +{.}{ }             {how}
  +{}{ \parenthesize} {date}
  +{,} { available at \eprint}        {eprint}
  +{,}{ available at \url}{url}
  +{,}{ }             {note}
  +{.}{}              {transition}
}
\usepackage{amssymb, amsfonts, amsxtra, amsmath}
\usepackage{mathrsfs}
\usepackage{mathdots}
\usepackage{wasysym}
\usepackage[all]{xy}
\usepackage{bbm}
\usepackage{calc}
\usepackage{accents}
\usepackage{nicematrix}

\usepackage[bbgreekl]{mathbbol}

\usepackage{stackengine}

\usepackage{enumitem}

\numberwithin{equation}{section}

\DeclareSymbolFontAlphabet{\mathbb}{AMSb}
\DeclareSymbolFontAlphabet{\mathbbl}{bbold}

\newtheorem{Theorem}{Theorem}[section]
\newtheorem*{Theorem*}{Theorem}
\newtheorem{Def}[Theorem]{Definition}
\newtheorem*{Def*}{Def}
\newtheorem{Lem}[Theorem]{Lemma}
\newtheorem{Prop}[Theorem]{Proposition}
\newtheorem{Cor}[Theorem]{Corollary}

\newtheorem{Rem}[Theorem]{Remark}

\newtheorem{Exa}[Theorem]{Example}

\newcommand\bp{\begin{proof}}
\newcommand\ep{\end{proof}}

\mathchardef\mhyph="2D

\DeclareMathOperator{\ad}{\mathrm{ad}}

\DeclareMathOperator{\Ad}{\mathrm{Ad}}

\DeclareMathOperator{\End}{\mathrm{End}}

\DeclareMathOperator{\Hom}{\mathrm{Hom}}
\DeclareMathOperator{\id}{\mathrm{id}}

\DeclareMathOperator{\Rep}{\mathrm{Rep}}
\DeclareMathOperator{\rd}{\mathrm{d}\!}

\newcommand{\ev}{\mathrm{ev}}
\newcommand{\cop}{\mathrm{cop}}
\newcommand{\op}{\mathrm{op}}

\newcommand{\red}{\mathrm{red}}

\newcommand{\msF}{\mathscr{F}}

\newcommand{\msM}{\mathscr{M}}
\newcommand{\msN}{\mathscr{N}}

\newcommand{\msR}{\mathscr{R}}

\newcommand{\mcC}{\mathcal{C}}
\newcommand{\mcD}{\mathcal{D}}

\newcommand{\Hsp}{\mathcal{H}}
\newcommand{\Gsp}{\mathcal{G}}
\newcommand{\Ksp}{\mathcal{K}}

\newcommand{\mcK}{\mathcal{K}}

\newcommand{\C}{\mathbb{C}}
\newcommand{\G}{\mathbb{G}}

\newcommand{\T}{\mathbb{T}}

\newcommand{\Z}{\mathbb{Z}}

\newcommand{\opp}{\mathrm{op}}

\newcommand{\Res}{\mathrm{Res}}

\newcommand{\Span}{\mathrm{Span}}

\newcommand{\Corr}{\mathrm{Corr}}

\newcommand{\Compact}{\mathcal{K}}

\newcommand{\slice}{\mathrm{s}}
\newcommand{\wW}{\widetilde{W}}

\begin{document}

\title{Approximation properties for dynamical W$^*$-correspondences}
\author{K.\ De Commer and J.\ De Ro}
\address{Vrije Universiteit Brussel}
\email{kenny.de.commer@vub.be}
\email{joeri.ludo.de.ro@vub.be}

\begin{abstract}
Let $\G$ be a locally compact quantum group, and $A,B$ von Neumann algebras on which $\G$ acts. We refer to these as \emph{$\G$-dynamical W$^*$-algebras}. We make a study of $\G$-equivariant $A$-$B$-correspondences, that is, Hilbert spaces $\Hsp$ with an $A$-$B$-bimodule structure by $*$-preserving normal maps, and equipped with a unitary representation of $\G$ which is equivariant with respect to the above bimodule structure. Such structures are a Hilbert space version of the theory of $\G$-equivariant Hilbert C$^*$-bimodules. We show that there is a well-defined Fell topology on equivariant correspondences, and use this to formulate approximation properties for them. Within this formalism, we then characterize amenability of the action of a locally compact group on a von Neumann algebra, using recent results due to Bearden and Crann. We further consider natural operations on equivariant correspondences such as taking opposites, composites and crossed products, and examine the continuity of these operations with respect to the Fell topology.  
\end{abstract}

\maketitle

\section*{Introduction}

In the setting of representation theory and abstract harmonic analysis, amenability properties can often be characterized by the possibility to approximate a pre-determined restricted class of `trivial representations' within a pre-determined restricted class of `regular representations'. We then say that the former are `weakly contained' in the latter. These classes consist in many cases of just one object. Two examples stand out: 
\begin{itemize}
\item If $G$ is a locally compact group, then $G$ is amenable if and only if the trivial representation is weakly contained in the regular representation on $L^2(G)$.
\item If $A$ is a von Neumann algebra, then $A$ is amenable if and only if the trivial $A$-bimodule $L^2(A)$ is weakly contained in the regular (or \emph{coarse}) $A$-bimodule $L^2(A)\otimes L^2(A)$.  
\end{itemize} 
In this last example, we mean bimodules in the sense of \emph{correspondences} \cite{Con80+} (see also \cite{Pop86}, \cite{Con94}*{Appendix B}): we are given two von Neumann algebras $A,B$ and a Hilbert space $\Hsp$ equipped with a normal unital $*$-representation $\pi$ of $A$ and a normal unital anti-$*$-representation $\rho$ of $B$ for which $\Hsp$ forms an $A$-$B$-bimodule. 

Characterizing amenability in this setting has led to some early breakthrough results in the modern theory of von Neumann algebras by its relation to \emph{injectivity} \cites{Con76,CE77,EL77}, and continues to play a fundamental role in many aspects of operator theory. For example, the notion of \emph{relative amenability}, introduced in \cite{Pop86} and further explored in the general context of correspondences in \cites{A-D90,A-D95}, was only recently shown to be equivalent with the corresponding notion of relative injectivity \cite{BMO20}.

By means of the crossed product construction, close ties can be made between amenability of \emph{actions} of locally compact groups on von Neumann algebras \cite{A-D79} and approximation properties related to correspondences constructed from the von Neumann algebra and its crossed product. In the setting of dynamical systems, such approximation properties are connected to amenability properties for certain associated groupoids \cite{A-DR00}. We also mention the recent key results relating amenability of C$^*$-dynamical systems to amenability of W$^*$-dynamical systems \cites{BEW21,OS21}.

In \cites{KV00,KV03}, a generalization was made of the theory of locally compact groups to that of \emph{locally compact quantum groups} (we refer to these articles for more on the history of the subject). It is based on the notion of a Hopf-von Neumann algebra \cite{Ern67}, which we now recall. 

\begin{Def}
A Hopf-von Neumann algebra consists of a von Neumann algebra $M$ with a unital normal faithful $*$-homomorphism $\Delta: M \rightarrow M \bar{\otimes} M$ satisfying the \emph{coassociativity} condition
\[
(\Delta\otimes \id)\circ \Delta = (\id\otimes \Delta)\circ \Delta.
\]
It is called \emph{involutive} if it is equipped with an involutive anti-$*$-automorphism $R: M\rightarrow M$ satisfying 
\[
\Delta\circ R = (R\otimes R)\circ \Delta^{\opp},
\]
where $\Delta^{\opp}$ is the composition of $\Delta$ with the flip map $\varsigma: x\otimes y \mapsto y\otimes x$.
\end{Def}
If $(M,\Delta)$ is a Hopf-von Neumann algebra, we sometimes write $\G= (M,\Delta)$ and $M = L^{\infty}(\G)$. Then we say that $\G$ is a \emph{locally compact quantum group} if moreover $(M,\Delta)$ possesses invariant weights (Definition \ref{Deflcqg}). Such Hopf-von Neumann algebras are automatically involutive with respect to a canonical involution $R$.

The von Neumann algebraic framework for locally compact quantum groups is very flexible and high-powered. For example, consider the following notions for a Hopf-von Neumann algebra $\G = (M,\Delta)$. 

\begin{Def}\label{DefUnitRep}
A \emph{unitary representation} of $\G$ on a Hilbert space $ \Hsp$ is a unitary $U \in B(\Hsp)\bar{\otimes}M$ such that 
\begin{equation}\label{EqUnitRepDef}
(\id\otimes \Delta)U = U_{12}U_{13}.
\end{equation}
\end{Def}
Here $U_{12} = U\otimes 1$, while $U_{13}$ is $(\id\otimes \varsigma)(U\otimes 1)$. 

\begin{Def}
If $A$ is a von Neumann algebra, an action of $\G$ on $A$ is a unital normal faithful $*$-homomorphism
\[
\alpha: A\rightarrow A\bar{\otimes}M,
\]
called the associated \emph{coaction},  and satisfying the \emph{coaction property} 
\[
(\alpha \otimes \id)\alpha = (\id\otimes \Delta)\alpha.
\]
We also say that $(A,\alpha)$ is a $\G$-dynamical von Neumann algebra, or simply a $\G$-W$^*$-algebra.
\end{Def}

In the general context of Hopf-von Neumann algebras, these two notions have no apparent relation, but when $\G= (M,\Delta)$ is a locally compact quantum group, there is a beautiful generalisation, due to S. Vaes \cite{Vae01}, of Haagerup's canonical implementation of a group action on a von Neumann algebra \cite{Haa78}: there exists a canonical unitary representation $U_{\alpha}$ of $\G$ on the standard Hilbert space $L^2(A)$ such that 
\begin{equation}\label{EqImplmCan}
\alpha(a) = U_{\alpha}(a\otimes 1)U_{\alpha}^*,\qquad \forall a\in A, 
\end{equation}
where we view $A\subseteq B(L^2(A))$ in standard position.

From an abstract and even an aesthetic point of view, one of the main salient features of locally compact quantum groups is that they allow an extension of the Pontryagin duality for arbitrary \emph{abelian} locally compact groups to a Pontryagin duality for arbitrary locally compact \emph{quantum} groups. When restricting to locally compact groups $G$, this is the correspondence sending the function algebra $L^{\infty}(G)$ (with its natural coproduct encoding the group structure) into the group von Neumann algebra $\mathscr{R}(G)$ generated by right translation operators (with its natural diagonal coproduct), and thus not an operation which gives a duality on a single category in this restricted setting. Various approximation properties such as amenability, Haagerup property, property (T),... can also be formulated naturally within the much more general setting of locally compact quantum groups, see e.g.\  \cites{Voi77,Rua96,DQV02,Tom06,Fim10,KyS12,SV14,CN15,DFSW16}.

In this paper, we combine the above two settings by studying \emph{equivariant correspondences}. This notion is a verbatim generalisation of the corresponding notion for $\mathbb{G}$-C$^*$-algebras with $\G$ a compact quantum group, introduced in \cite{AS21}.
\begin{Def}\label{DefEquivCorrMain}
Let $\G= (M,\Delta,R)$ be an involutive Hopf-von Neumann algebra, and let $(A,\alpha),(B,\beta)$ be $\G$-W$^*$-algebras. A \emph{$\G$-equivariant $A$-$B$-correspondence}, or simply a $\G$-$A$-$B$-correspondence, consists of a Hilbert space $\Hsp$ with 
\begin{itemize}
\item a normal unital $*$-representation $\pi$ of $A$ on $\Hsp$, 
\item a normal unital anti-$*$-representation $\rho$ of $B$ on $\Hsp$, and 
\item a unitary representation $U$ of $\G$ on $\Hsp$,
\end{itemize} 
satisfying for all $a\in A$ and $b\in B$ the compatibility relations $\pi(a)\rho(b) = \rho(b)\pi(a)$ and
\begin{equation}\label{EqDefCondGeqCorr}
U(\pi(a)\otimes 1)U^* = (\pi\otimes \id)\alpha(a),\qquad  U^*(\rho(b)\otimes 1)U = (\rho \otimes R)(\beta(b)).
\end{equation}
\end{Def} 

C$^*$-algebraic analogues of such equivariant correspondences, in the context of equivariant Hilbert C$^*$-modules, have already appeared in the original works \cites{BS89,BS93} in relation to equivariant KK-theory, and have remained important in follow-up works (see e.g.\ \cite{BC17}). In \cite{AS21}, concrete use was made rather of the above Hilbert space picture of equivariant correspondences, however within the context of C$^*$-algebras and \emph{compact} quantum groups. On the other hand, in \cite{Vae05} there already appear specific equivariant correspondences, albeit in the particular case arising from inclusions of quantum groups, and using instead the (equivalent) notion of equivariant von Neumann bimodules.

More particular cases of interest have turned up recently within the setting of \emph{quantum automorphism groups} \cite{RV22} (for related results in a more classical setting, see \cite{AV16}), and \emph{quantizations of semisimple Lie groups} \cite{DCDT21}. It is our hope that the framework of this paper will be a convenient one to study analytic approximation properties for these important classes of examples.

Let us now come to the main results of this paper. 

\begin{Theorem*}[Theorem \ref{TheoMMM}]
Let $\G$ be a locally compact quantum group, and $A,B$ two $\G$-W$^*$-algebras. Then there exists a universal C$^*$-algebra $C^{\G}(A,B)$ such that any $\G$-equivariant $A$-$B$-correspondence $\Hsp$ uniquely lifts to a non-degenerate $*$-representation $\theta_{\Hsp}$ of $C^{\G}(A,B)$ on $\Hsp$. 
\end{Theorem*}

The C$^*$-algebra $C^{\G}(A,B)$ is a C$^*$/W$^*$-intermediate of the double cross product C$^*$-algebra appearing in \cite{AS21} (see also \cite{NV02}*{Section 2.6} for a description in a purely algebraic context). 

Due to the above theorem, we can introduce a \emph{Fell topology} on the class of $\G$-equivariant $A$-$B$-correspondences, as well as a notion of weak containment: we say $\Gsp$ is weakly contained in $\Hsp$ if the $*$-representation $\theta_{\Gsp}$ of $C^{\G}(A,B)$ factors through $\theta_{\Hsp}(C^{\G}(A,B))$. When we want to emphasize that we are working in the equivariant setting, we will sometimes use the term \emph{equivariantly weakly contained}.

Our next main result offers another characterisation of equivariant weak containment. If $(\Hsp,U_{\Hsp}),(\Gsp,U_{\Gsp})$ are $\G$-$A$-$B$-correspondences, we use the notation 
\begin{equation}\label{EqHilbertMod}
X_B(\Gsp,\Hsp) = \{x: \Gsp\rightarrow \Hsp\mid x \rho_{\Gsp}(b) = \rho_{\Hsp}(b)x, \forall b\in B\},
\end{equation}
and 
\[
\alpha_{X_B(\Gsp, \Hsp)}: X_B(\Gsp,\Hsp) \rightarrow  X_B(\Gsp,\Hsp) \bar{\otimes} M,\qquad x \mapsto U_{\Hsp}(x\otimes 1)U_{\Gsp}^*.
\]
\begin{Theorem*}[Theorem \ref{TheoMain1}] 
Let $\G$ be a locally compact quantum group, and let $A,B$ be $\G$-W$^*$-algebras. Let $\Hsp,\Gsp$ be $\G$-$A$-$B$-correspondences, and put $\Hsp^{\infty} = \oplus_1^{\infty}\Hsp$. 

Then $\Gsp$ is equivariantly weakly contained in $\Hsp$ if and only if there exists a $1$-bounded net $y_{i} \in X_B(\Gsp,\Hsp^{\infty})$ such that for all $a\in A$
\begin{equation}
(y_i^*\pi_{\Hsp^{\infty}}(a)\otimes 1)\alpha_{X_B(\Gsp,\Hsp)}(y_i) \rightarrow \pi_{\Gsp}(a)\otimes 1 \quad \sigma\textrm{-weakly in }B(\Gsp)\bar{\otimes} M.
\end{equation}
\end{Theorem*}

The above results naturally lead to a notion of amenability for inclusions of $\G$-W$^*$-algebras: if $A\subseteq B$ is $\G$-equivariant, then we can say that this inclusion is \emph{strongly equivariantly amenable} if ${}_AL^2(A)_A$ is equivariantly weakly contained in ${}_AL^2(B)_A$ (these correspondences can easily be made equivariant using the results of \cite{Vae01}, see Example \ref{Exa1}). The following theorem then shows that equivariant strong amenability implies \emph{equivariant amenability} (=equivariant relative injectivity). 

\begin{Theorem*}[Theorem \ref{TheoInclStrong}] Let $\G$ be a locally compact quantum group, and let $(A,\alpha)\subseteq (B,\beta)$ be an equivariant inclusion of $\G$-W$^*$-algebras. If  ${}_AL^2(A)_A$ is equivariantly weakly contained in ${}_AL^2(B)_A$, then there exists a \emph{$\G$-equivariant} conditional expectation $E: B \rightarrow A$, i.e.\
\begin{equation}\label{EqEquivCEE}
(E \otimes \id)\circ \beta = \alpha \circ E.
\end{equation}
\end{Theorem*} 

As a particular case in point, one can consider for a $\G$-W$^*$-algebra $(A,\alpha)$ the equivariant inclusion 
\[
A\cong \alpha(A) \subseteq  A\bar{\otimes}M,
\] 
where $A\bar{\otimes}M$ is endowed with the $\G$-action $\id_A\otimes \Delta$. Our next result is then a converse of the last result above, in the case of a usual locally compact group $G$. It heavily relies on \cite{BC22}. 

\begin{Theorem*}[Theorem \ref{TheoClassCase}]
Let $G$ be a locally compact group, and let $A$ be a $G$-W$^*$-algebra via
\[
\alpha: A\rightarrow A\bar{\otimes} L^{\infty}(G). 
\]
Then ${}_AL^2(A)_A$ is weakly contained in ${}_{\alpha(A)}(L^2(A)\otimes L^2(G))_{\alpha(A)}$ if and only if there exists a $G$-equivariant conditional expectation $E: A\bar{\otimes}L^{\infty}(G) \rightarrow A$.
\end{Theorem*}

As mentioned, an important aspect of locally compact quantum groups concerns their Pontryagin duality theory. We recall (see Section \ref{SecLCQGAct} for details) that if $A$ is a $\G$-dynamical von Neumann algebra, and $\check{\G}$ the Pontryagin dual, one can construct the \emph{crossed product von Neumann algebra} $A\rtimes \G$, which becomes a $\check{\G}$-dynamical von Neumann algebra. We will show that this duality extends to the level of $\G$-equivariant correspondences: to any $\G$-equivariant $A$-$B$-correspondence $\Hsp$ is associated a $\check{\G}$-equivariant $A\rtimes \G$-$B\rtimes \G$-correspondence $\Hsp\rtimes \G$.

Our final main result is then the following topological refinement of the well-known Takesaki-Takai duality. 

\begin{Theorem*}[Theorem \ref{TheoTakTakTakTak} and Theorem \ref{TheoAsMult}] 
The assignment 
\[
\Hsp \mapsto \Hsp\rtimes \G
\] 
is an equivalence between the W$^*$-category of $\G$-equivariant $A$-$B$-correspondences and the W$^*$-category of $\check{\G}$-equivariant $A\rtimes \G$-$B\rtimes \G$-correspondences. Moreover, this correspondence preserves the Fell topology.
\end{Theorem*} 

The precise contents of this paper are as follows. 

In the \emph{first section}, we provide some preliminaries on locally compact quantum groups and their actions. In the \emph{second section}, we introduce equivariant correspondences, provide some elementary examples, and show that the category of equivariant correspondences admits a well-defined and well-behaved Fell topology. In the \emph{third section}, we introduce the notion of weak containment between equivariant correspondences, and provide an equivalent formulation for weak containment in terms of one-sided intertwiners. In the \emph{fourth section}, we introduce the notion of amenability and strong amenability for equivariant correspondences, and show that these notions can be used to introduce corresponding notions for actions of locally compact quantum groups on von Neumann algebras. We show that in the setting of locally compact groups, an amenable action is automatically strongly amenable, using the results of \cite{BC22}. In the \emph{fifth section}, we consider various operations on equivariant correspondences such as opposites, composites and crossed products. The main result here is a generalized Takesaki-Takai-duality, providing an equivalence of W$^*$-bicategories between $\G$-equivariant correspondences and $\check{\G}$-equivariant correspondences. In the \emph{sixth section}, we consider the interplay between the above constructions and the Fell topology, and consider some direct corollaries for the associated approximation properties. We end in the \emph{seventh section} with a brief outlook.

\textbf{Notations and conventions}

All our vector spaces are complex. If $V$ is a vector space and $S \subseteq V$, we denote $\Span(S)$ for the linear span of $S$. If $V,W$ are subspaces of an algebra $A$, we denote 
\[
VW = \Span\{vw\mid v\in V,w\in W\}\subseteq A.
\]
If $V$ is a normed vector space and $S \subseteq V$, we denote $[S]$ for the norm-closure of the linear span of $S$.
We denote $\odot$ for algebraic tensor products (over $\C$). We use $\otimes$ for the tensor product of Hilbert spaces, or the minimal tensor product of C$^*$-algebras or Hilbert C$^*$-modules. We use $\bar{\otimes}$ for the spatial tensor product of von Neumann algebras. The inner product of Hilbert spaces is assumed anti-linear in the first variable. 

When referring to linear maps between operator spaces, we use the abbrevation cb for \emph{completely bounded}, and, when considering linear maps between operator systems, we use the abbreviation (u/c)cp for \emph{(unital/contractive) completely positive}.

We will use the terms W$^*$-algebra and von Neumann algebra interchangeably. When $M$ is a von Neumann algebra, we denote $M_*$ for the predual and $M_+$ for the positive cone. If $\Phi$ is an nsf (= normal, semifinite, faithful) weight on $M$, we use the standard notations.
\[
\mathscr{M}_\Phi^+= \{x\in M^+: \Phi(x)< \infty\}, \quad\msN_{\Phi} = \{x\in M\mid \Phi(x^*x)<\infty\}.
\]

We let $L^2(M)$ be the standard Hilbert space of $M$, with modular conjugation $J = J_M$ and $\Lambda_{\Phi}: \msN_{\Phi}\rightarrow L^2(M)$ the canonical GNS-map of $\Phi$. We denote $\pi_M$ for the associated left representation of $M$ on $L^2(M)$, and 
\begin{equation}\label{EqRightRep}
\rho_M(x) = J\pi_M(x)^*J
\end{equation}
for the associated anti-$*$-representation of $M$ on $L^2(M)$. When there is no ambiguity, we drop the symbol $\pi_M$, and identify $M \subseteq B(L^2(M))$. For the basics of weight theory on von Neumann algebras, we refer to \cite{Tak03}.

If $\Hsp$ is a Hilbert space, we denote $\Compact(\Hsp)$ for the C$^*$-algebra of compact operators on $\Hsp$. If $u \in B(\Hsp)$ is a unitary, we write $\Ad(u)(x) = \Ad_u(x) = uxu^*$ for $x\in B(\Hsp)$. More generally, when $\Hsp,\Gsp$ are Hilbert spaces, and $x\in B(\Gsp,\Hsp)$, we will denote by $\Ad_x$ the completely positive map
\[
\Ad_x: B(\Gsp)\rightarrow B(\Hsp),\qquad T\mapsto xTx^*.
\]

If $\xi,\eta \in \Hsp$, we write
\[
\omega_{\xi,\eta} \in B(\Hsp)_*,\qquad B(\Hsp)\ni x \mapsto \langle \xi,x\eta\rangle \in \C.
\] 
If $\Hsp,\Gsp$ are Hilbert spaces, we denote the \emph{flip map} by  
\[
\Sigma: \Hsp \otimes \Gsp \rightarrow \Gsp \otimes \Hsp,\qquad \xi\otimes \eta\mapsto \eta\otimes \xi.
\]
We use the \emph{leg numbering notation}: if $x$ is an operator on a tensor product of Hilbert spaces, we write 
\[
x_{12} = x\otimes 1,\qquad x_{23}= 1\otimes x,\qquad x_{13}= (\id\otimes \Ad\Sigma)x_{12} = (\Ad\Sigma\otimes \id)x_{23}. 
\]
In practice, there should be no ambiguity on where the unit acts. 

For $C$ a C$^*$-algebra, we write $C^{**}$ for the universal W$^*$-envelope of $C$. Throughout the theory of operator algebras and quantum groups, the multiplier algebra $M(C)$ of $C$ plays a prominent role. However, in this paper, we will encounter a situation where the use of the multiplier $C^*$-algebra is no longer sufficient, and where it is necessary (cfr. Remark \ref{example quasimultiplier}) to work with an even larger object, less familiar than the multiplier algebra. 

The set of \emph{quasi-multipliers} of $C$, originally introduced in \cite{AP73}, is defined as the norm-closed subspace
\[
QM(C) = \{x\in C^{**}\mid CxC \subseteq C\} \subseteq C^{**}.
\]

In general, the space $QM(C)$ is very complicated, partly due to its lack of multiplicative structure. To give some intuition to what kind of elements this space contains, we provide another description for it as follows: Write $\operatorname{Bil}_Q(C)$ for the vector space of bilinear maps $\phi: C\times C \to C$ satisfying
$$\phi(pq,rs) = p\phi(q,r)s, \quad p,q,r,s \in C.$$
It is clear that there is a canonical injective linear map, given by
$$QM(C)\to \operatorname{Bil}_Q(C): x \mapsto (\phi_x: (p,q)\mapsto pxq).$$ It follows from \cite{AP73}*{Proposition 4.2} that this map is bijective. Therefore, we can identify $QM(C)\cong \operatorname{Bil}_Q(C)$ as linear spaces. In very specific examples, an explicit description of the space $QM(C)$ can be given, see e.g. \cite{Li89}*{Section 8}.

We refer to \cite{Ped79}*{Section 3.12} for further basic facts on quasi-multipliers, but we recall below the main properties that we will use. First of all, we emphasize that $QM(C)$, although not necessarily multiplicatively closed, always naturally has the structure of an operator system.
Next, we recall from \cite{Ped79} that whenever $C \subseteq B(\Hsp)$ is a faithful, non-degenerate $*$-representation, the corresponding normal unital $*$-homomorphism $C^{**} \rightarrow B(\Hsp)$ is faithful on $QM(C)$, and we may isometrically identify through this map
\[
QM(C) = \{x\in B(\Hsp) \mid CxC \subseteq C\}.
\]
For a general non-degenerate $*$-representation $\pi: C \rightarrow B(\Hsp)$, we will always extend it to a map $\pi: QM(C) \rightarrow B(\Hsp)$ by restricting the natural normal extension of $\pi$ to $C^{**}$. In particular, we can use the formula
\[
\pi(c)\pi(q)\pi(d) = \pi(cqd),\qquad c,d\in C,q\in QM(C).
\]

Given $C^*$-algebras $C,D$ and a ccp map $\phi: C \to D$, the \emph{multiplicative domain} of $\phi$ \cite{BO08}*{Proposition 1.5.7} is the $C^*$-subalgebra 
$$C_\phi:= \{x\in C: \phi(x^*x) = \phi(x)^*\phi(x) \textrm{\ and\ } \phi(xx^*)= \phi(x)\phi(x)^*\}.$$
We will frequently use that if $x\in C_\phi$ and $y\in C$, then $\phi(xy) = \phi(x)\phi(y)$ and $\phi(yx) = \phi(y)\phi(x).$ 

\section{Locally compact quantum groups and their actions}
\label{SecLCQGAct}

We recall some basic results from the theory of locally compact quantum groups \cites{KV00,KV03,VV03}, and their representations and actions \cites{Kus01,Vae01}. The following is \cite{KV03}*{Definition 1.1}.

\begin{Def}\label{Deflcqg}
A \emph{locally compact quantum group $\G$} consists of a Hopf-von Neumann algebra $\G = (M,\Delta)$ possessing nsf weights $\Phi,\Psi: M_+ \rightarrow [0,+\infty]$ satisfying  the \emph{left} and the \emph{right} invariance conditions 
\begin{equation}\label{EqLeftRightInv}
\Phi((\omega \otimes \id)\Delta(x)) = \Phi(x)\omega(1), \quad \Psi((\id \otimes \omega)\Delta(y)) = \Psi(y)\omega(1), \quad \omega \in M_*, \quad x \in \mathscr{M}_\Phi^+,\quad  y \in \mathscr{M}_\Psi^+.
\end{equation}
\end{Def}
These weights are then necessarily unique up to multiplication by a positive scalar.

Whenever convenient, we write $M = L^{\infty}(\G), M_* = L^1(\G)$ and $L^2(M) = L^2(\G)$. We identify $M \subseteq B(L^2(M))$ by the standard representation, which we write explicitly as $\pi_M$ if we want to emphasize it. We note that a locally compact quantum group automatically satisfies the density conditions
\begin{equation}\label{EqDensGal}
[(M\otimes 1)\Delta(M)]^{\sigma\textrm{-weak}} = [\Delta(M)(1\otimes M)]^{\sigma-\textrm{weak}} = M\bar{\otimes}M.
\end{equation}

We write the induced \emph{convolution product} on $M_*$ as 
\[
\omega*\omega' := (\omega \otimes \omega')\Delta,\qquad \omega,\omega'\in M_*.
\] 

Fundamental to the theory of locally compact quantum groups are the unitaries 
\[
V \in B(L^2(M))\bar{\otimes}M,\qquad W \in M \bar{\otimes} B(L^2(M)),
\]
called respectively \emph{right} and \emph{left} regular unitary representation. They are uniquely characterized by the identities 
\[
(\id\otimes \omega)(V) \Lambda_{\Psi}(x) = \Lambda_{\Psi}((\id\otimes \omega)\Delta(x)),
\qquad \omega \in M_*,x\in \msN_{\Psi},
\]
\[
 (\omega \otimes \id)(W^*)\Lambda_{\Phi}(x) = \Lambda_{\Phi}((\omega\otimes \id)\Delta(x)),\qquad \omega \in M_*,x\in \msN_{\Phi}.
\]
They are \emph{multiplicative unitaries} \cite{BS93} in the sense that as operators on $L^2(M)\otimes L^2(M)$, they satisfy
\[
V_{12}V_{13}V_{23} = V_{23}V_{12},\qquad W_{12}W_{13}W_{23}= W_{23}W_{12}. 
\]
These multiplicative unitaries implement the coproduct of $M$:
\begin{equation}\label{EqComultImpl}
W^*(1\otimes x)W = \Delta(x) = V(x\otimes 1)V^*,\qquad x\in M.
\end{equation}
Moreover, we have the equality
\[
[(\omega\otimes \id)V \mid \omega \in B(L^2(M))_*] = [(\id\otimes \omega)W \mid \omega \in B(L^2(M))_*],
\]
and this defines a $\sigma$-weakly dense C$^*$-subalgebra of $M$ that we will write as 
\[
C_0^{\red}(\G)\subseteq M.
\] 

We recall some aspects of the duality theory of locally compact quantum groups \cite{KV03}*{Section 2}. We write
\[
\hat{M} = [(\omega\otimes \id)W \mid \omega \in M_*]^{\sigma\textrm{-weak}},
\]
which can be shown to be a von Neumann algebra. It becomes a locally compact quantum group for the coproduct 
\[
\hat{\Delta}(x) = \Sigma W(x\otimes 1)W^*\Sigma, \qquad x\in \hat{M}.
\]
When useful for the sake of intuition or terminology, we write $\hat{M} = L^\infty(\hat{\G})$: it is the left group von Neumann algebra of $\G$, viewed as a function algebra on what is to be considered the Pontryagin dual $\hat{\G}$ of $\G$.

As we will need to compare left and right coactions later on, it will also be convenient to introduce independently the dual $\check{M} = M^{\cop\wedge}$ of $M^{\cop} = (M,\Delta^{\opp})$, where $\Delta^{\opp} = \Ad\Sigma \circ \Delta$. Noting that $(M,\Delta^{\opp})$ has the left and right regular unitary representation equal to respectively $V_{21}^*$ and $W_{21}^*$, we have
\[
\check{M} = [(\id\otimes \omega)V \mid \omega \in M_*]^{\sigma\textrm{-weak}}
\]
with coproduct
\[
\check{\Delta}(x) = V^*(1\otimes x)V, \qquad x\in \check{M}.
\]
It can be shown that $\check{M} = \hat{M}'$. We then also occasionally write $\check{M} = \msR(\G) = L^{\infty}(\check{\G})$, and view it as either the right group von Neumann algebra of $\G$, or a function algebra on an `opposite' Pontryagin dual $\check{\G}$ of $\G$. We further note that 
\begin{equation}\label{EqIdSliceMaps}
C_0^{\red}(\hat{\G}) = [(\omega\otimes \id)W \mid \omega \in M_*],\qquad C_0^{\red}(\check{\G}) = [(\id\otimes \omega)V \mid \omega \in M_*].
\end{equation}

There is a canonical identification of $L^2(M)$ and $L^2(\check{M})$, in such a way that  
\[
\check{W}= V
\]
is the left regular representation of $\check{\G}$. With $\check{J} = J_{\check{M}}$ the resulting associated modular conjugation for $\check{M}$ on $L^2(M)$, we have
\[
\check{J}M\check{J} = M.
\]
One calls the induced anti-multiplicative involution
\[
R: M \rightarrow M,\quad x\mapsto \check{J}x^*\check{J}
\]
the \emph{unitary antipode}. It turns $(M,\Delta)$ into an involutive Hopf-von Neumann algebra. The following identities involving the modular conjugations will be useful:
\begin{equation}\label{EqMultUnJhatJ}
(\check{J}\otimes J)W(\check{J}\otimes J) = W^*,\qquad (J\otimes \check{J})V(J\otimes \check{J}) = V^*.
\end{equation}
We also recall from \cite{KV03}*{Corollary 2.12} that one can pick a canonical unimodular number $c\in \C$ such that 
\[
c\check{J}J = \overline{c} J\check{J}.
\]
The resulting self-adjoint unitary $u_{\G}:= c\check{J}J$ is called the \emph{fundamental symmetry}. This unitary allows us to relate $W$ and $V$ directly, for we have 
\begin{equation}\label{EqIdVW}
(u_{\G}\otimes 1)V(u_{\G}\otimes 1)  = W_{21}.
\end{equation}
It also allows us to express the right regular unitary of $\check{M}$ as 
\[
\check{V} = \widetilde{W} :=  (u_{\G}\otimes u_{\G})W(u_{\G}\otimes u_{\G}).
\]
We similarly write $\widetilde{V} = (u_{\G}\otimes u_{\G})V(u_{\G}\otimes u_{\G})$. The operation $M\mapsto \check{M}$ is not on the nose involutive, but we have that 
\[
M^{\vee\vee} = u_{\G}Mu_{\G},\qquad \Delta^{\vee\vee} = \Ad(u_{\G}\otimes u_{\G})\circ \Delta \circ \Ad(u_{\G}). 
\]

Recall from Definition \ref{DefUnitRep} the notion of a \emph{unitary representation} of $\G$ on a Hilbert space $\Hsp$. If $\omega \in M_*$, we write 
\[
U(\omega) = (\id\otimes \omega)U \in B(\Hsp).
\]
Clearly 
\[
U(\omega*\omega') = U(\omega)U(\omega'),\qquad \omega,\omega'\in M_*.
\]

For later use, we recall the following facts, see \cite{Wor96}, \cite{Kus01}*{Theorem 1.6} and \cite{BDS13}*{Theorem 4.12}. 
\begin{Prop}\label{PropImRepC}
If $U$ is a unitary $\G$-representation on a Hilbert space $\Hsp$, then the norm-closure 
\[
C_U = [U(\omega)\mid \omega \in M_*]
\]
is a C$^*$-algebra acting non-degenerately on $\Hsp$, and
\[
U \in M(C_U\otimes C_0^{\red}(\G)) \subseteq M(\Compact(\Hsp)\otimes C_0^{\red}(\G)) \cap M(C_U \otimes \Compact(L^2(\G))).
\]  
\end{Prop} 

The following results from \cite{Vae01} will be crucial. Assume that $A$ is a  $\G$-W$^*$-algebra. The \emph{crossed product} von Neumann algebra is defined as 
\[
A\rtimes \G = \{\alpha(a)(1\otimes x)\mid a\in A, x\in \check{M}\}''  \subseteq A\bar{\otimes} B(L^2(\G)), 
\]
and in fact already
 \[
A\rtimes \G= [\alpha(a)(1\otimes x)\mid a\in A, x\in \check{M}]^{\sigma\textrm{-weak}}.
\]
In the following, we write 
\[
A^{\rtimes} := A\rtimes \G
\] 
as a shorthand. We note that $A^{\rtimes}$ carries the dual right $\check{\G}$-action 
\[
\alpha^{\rtimes}: A^{\rtimes} \rightarrow A^{\rtimes} \bar{\otimes} \check{M},\quad \alpha^{\rtimes}(z) = \widetilde{W}_{23}z_{12}\widetilde{W}_{23}^*,\qquad z\in A^{\rtimes}.
\] 
In particular, one has 
\[
\alpha^{\rtimes}(\alpha(a)) = \alpha(a)\otimes 1,\quad \alpha^{\rtimes}(1\otimes x) = 1\otimes \check{\Delta}(x),\qquad a\in A,x\in \check{M}. 
\]

\begin{Theorem}\label{TheoDualCoact}\label{TheoUnitImpl}
Let $A$ be a $\G$-W$^*$-algebra. Then there is a canonical identification 
\[
L^2(A^{\rtimes})\cong L^2(A)\otimes L^2(\G), 
\] 
in such a way that 
\begin{equation}\label{EqImplDualCoact}
U_{\alpha} = J_{A^{\rtimes}}(J_A\otimes \check{J}) \in B(L^2(A)\otimes L^2(\G))
\end{equation}
is a unitary $\G$-representation implementing $\alpha$ (in the sense of \eqref{EqImplmCan}). Moreover, the associated GNS-representation $\pi_{A^{\rtimes}}$ of $A^{\rtimes}$ becomes the identity representation of $A^{\rtimes}$ on $L^2(A)\otimes L^2(\G)$, while 
\[
\rho_{A^{\rtimes}}(\alpha(a)) = J_{A^{\rtimes}}\pi_{A^{\rtimes}}(a)^*J_{A^{\rtimes}} = \rho_A(a)\otimes 1, \qquad a\in A,
\]
\[
\rho_{A^{\rtimes}}(1\otimes x) = J_{A^{\rtimes}}(1\otimes \pi_{\check{M}}(x)^*)J_{A^{\rtimes}} = U_{\alpha}(1\otimes \rho_{\check{M}}(x))U_{\alpha}^*,\qquad x\in \check{M}.
\]
The dual $\check{\G}$-action has $\widetilde{W}_{23}$ as its canonical unitary implementation.
\end{Theorem} 
Note that we make here the switch from the left coactions of  \cite{Vae01} to right coactions by replacing $(M,\Delta)$ with $(M,\Delta^{\opp})$. 
\begin{proof}
See \cite{Vae01}*{Lemma 3.3, Definition 3.4, Definition 3.6 and Proposition 3.7}. The last claim follows from \cite{Vae01}*{Proposition 4.3}.
\end{proof}

\begin{Rem} 
From \eqref{EqImplDualCoact}, we obtain the identity 
\begin{equation}\label{EqCanImplUni}
(J_A\otimes \check{J})U_{\alpha}(J_A\otimes \check{J}) = U_{\alpha}^*.
\end{equation}
\end{Rem}
\begin{Rem}
When $A =M$ with $\alpha = \Delta$, we have $U_{\Delta} =V$. 
\end{Rem}

We end this section with a concrete computation of a canonical implementing unitary. 

\begin{Exa}\label{ExaAdjAct}
Let $A$ be a $\G$-W$^*$-algebra. We can endow $A^{\rtimes}  = A\rtimes \G$ with the \emph{adjoint $\G$-action}
\begin{equation}\label{EqAdjAction}
\alpha_{\ad}: A^{\rtimes} \rightarrow A^{\rtimes} \bar{\otimes} M,\quad z \mapsto V_{23}z_{12}V_{23}^*,\qquad z\in A^{\rtimes},
\end{equation}
so 
\[
\alpha_{\ad}(\alpha(a)) = (\id\otimes \Delta)\alpha(a) = (\alpha\otimes \id)\alpha(a),\qquad \alpha_{\ad}(1\otimes x) = 1\otimes V(x\otimes 1)V^*,\qquad a\in A,x\in \check{M}.
\]

We claim that the canonical unitary implementation $U_{\ad}$ of $\alpha_{ad}$ is given by 
\begin{equation}\label{EqCanImplAdj}
U_{\ad} = V_{23}W_{32}U_{\alpha,13}.
\end{equation}

Indeed, note first that the crossed product von Neumann algebra $A^{\rtimes}\rtimes \G$ is given by 
\begin{eqnarray*}
A^{\rtimes}\rtimes \G &=& 
[\alpha_{\ad}(\alpha(a)(1\otimes x))(1\otimes 1\otimes y)\mid a\in A,x,y\in \check{M}]^{\sigma\textrm{-weak}}\\
&=& [V_{23}\alpha(a)_{12}x_2V_{23}^* y_3\mid a\in A,x,y\in \check{M}]^{\sigma\textrm{-weak}}\\
&=& [V_{23}\alpha(a)_{12}x_2 \check{\Delta}(y)_{23}V_{23}^*\mid a\in A,x,y\in \check{M}]^{\sigma\textrm{-weak}}\\
&\underset{\eqref{EqDensGal}}{=}& [V_{23}\alpha(a)_{12}(1\otimes x \otimes y)V_{23}^*\mid a\in A,x,y\in \check{M}]^{\sigma\textrm{-weak}}\\
&=& V_{23}(A^{\rtimes}\bar{\otimes} \check{M})V_{23}^*,
\end{eqnarray*} 
Since the dual coaction of $\check{M}$ on $A^{\rtimes}\rtimes \G$ is given by conjugation with $\widetilde{W}_{34}$, we see that $A^{\rtimes}\rtimes \G$ is just a conjugated copy of $A^{\rtimes}\bar{\otimes} \check{M}$ with $\id_{A^{\rtimes}}\otimes \check{\Delta}$ as the coaction. In particular, under the canonical identification \[
L^2(A^{\rtimes}\rtimes \G) = L^2(A^{\rtimes})\otimes L^2(\G) = L^2(A)\otimes L^2(\G)\otimes L^2(\G),
\]
the associated modular conjugation of $A^{\rtimes}\rtimes \G$ is given by 
\[
J_{A^{\rtimes}\rtimes \G} = V_{23}(J_{A^{\rtimes}}\otimes \check{J})V_{23}^*,
\]
and so the canonical unitary implementation of $\alpha_{\ad}$ is by Theorem \ref{TheoUnitImpl} given by 
\begin{eqnarray*}
U_{\alpha_{\ad}} &=& J_{A^{\rtimes}\rtimes \G}(J_{A^{\rtimes}} \otimes \check{J}) \\
&=& V_{23}(J_{A^{\rtimes}}\otimes \check{J})V_{23}^*(J_{A^{\rtimes}} \otimes \check{J}) \\
&=& V_{23}U_{\alpha,12}(J_A\otimes \check{J}\otimes \check{J})V_{23}^*(J_A\otimes \check{J}\otimes \check{J})U_{\alpha,12}^* \\
&\underset{\eqref{EqMultUnJhatJ},\eqref{EqIdVW}}{=}& V_{23}U_{\alpha,12}W_{32}U_{\alpha,12}^* \\
&\underset{\eqref{EqUnitRepDef}}{=}& V_{23} W_{32} U_{\alpha,13},
\end{eqnarray*}
proving \eqref{EqCanImplAdj}.
\end{Exa}

\section{Dynamical W\texorpdfstring{$^*$}{star}-correspondences}\label{SecIntroAndExa}

\subsection{Definition and Examples}
Let $\G$ be a locally compact quantum group, and let $(A,\alpha),(B,\beta)$ be $\G$-dynamical von Neumann algebras. In Definition \ref{DefEquivCorrMain}, we introduced the notion of a $\G$-equivariant $A$-$B$-correspondence $\Hsp = (\Hsp,U,\pi,\rho)$. Note that if $\G$ is the trivial group, we find back the usual notion of an $A$-$B$-correspondence. On the other hand, when $A = B = \C$, we find back the notion of a unitary $\G$-representation.

\begin{Def}
We denote by $\Corr^{\G}(A,B)$ the category of $\G$-$A$-$B$-correspondences, with morphisms from $\Hsp$ to $\Hsp'$ consisting of bounded operators $T: \Hsp \rightarrow \Hsp'$ intertwining the (anti-)representations of $A,B$ and $M_*$. 
\end{Def}
Clearly, $\Corr^{\G}(A,B)$ is a W$^*$-category \cite{GLR85}, and we will always understand it as such. 

Here are some examples of $\G$-equivariant correspondences. 

\begin{Exa}\label{Exa1}
Let $A$ be a $\G$-dynamical von Neumann algebra. The \emph{trivial} or \emph{identity} $\G$-$A$-$A$-correspondence $E_A^{\G}$ is $L^2(A)$ with its canonical representation $\pi_A$ and anti-representation $\rho_A$, and the canonical unitary $\G$-representation $U_{\alpha}$ implementing $\alpha$.
\end{Exa}

The compatibility relations \eqref{EqDefCondGeqCorr} are satisfied by Theorem \ref{TheoUnitImpl} and the identity \eqref{EqCanImplUni}, together with the fact that $\check{J}$ implements the unitary antipode.

Let $(A,\alpha)$ and $(B,\beta)$ be $\G$-dynamical von Neumann algebras. We say that $A \subseteq B$ is an \emph{equivariant inclusion} if $\beta_{\mid A} = \alpha$.
\begin{Exa}\label{ExaGenSub}\label{ExaSub}
 
Let $A\subseteq C$ and $B\subseteq D$ be $\G$-equivariant inclusions. Then for any $\G$-$C$-$D$ correspondence $\Hsp$ we obtain the restriction 
\[
{}_A\Res_B(\Hsp) = {}_A\Hsp_B. 
\]
In particular, if $A \subseteq B$ is an equivariant inclusion, we have equivariant correspondences ${}_AL^2(B)_B$, ${}_BL^2(B)_A$ and ${}_AL^2(B)_A$, each with $U_{\beta}$ as unitary representation for $\beta$ the given coaction on $B$. 
\end{Exa}
\begin{Exa}\label{ExaBasCons}\label{ExaBas}
Let $A,B$ be $\G$-dynamical von Neumann algebras, and let $(\Hsp,U,\pi,\rho)$ be a $\G$-equivariant correspondence. Consider $C := \rho(B)'$. Then $U(C\otimes 1)U^*\subseteq C \bar{\otimes} M$ by the bicommutant theorem, so that 
\[
\gamma: C \rightarrow C\bar{\otimes} M,\qquad c\mapsto U(c\otimes 1)U^*
\]
is a coaction, and $({}_C\Hsp_B,U)$ is a $\G$-equivariant $C$-$B$-correspondence. Moreover, $\pi(A) \subseteq C$ is a $\G$-equivariant inclusion. 

As a particular case, let $\pi:A \rightarrow B$ be a normal unital $\G$-equivariant $*$-homomorphism. Consider $C := \rho_B(\pi(A))'$, so that $\pi_B(\pi(A)) \subseteq \pi_B(B) \subseteq C$ is (by definition) the Jones basic construction. Then 
\[
\gamma: C \rightarrow C\bar{\otimes} M,\qquad c\mapsto U_{\beta}(c\otimes 1)U_{\beta}^*
\]
is a well-defined coaction, and $B\cong \pi_B(B) \subseteq C$ becomes a $\G$-equivariant inclusion of von Neumann algebras. 
\end{Exa}

\begin{Exa}\label{Exa2}
The \emph{coarse} (or regular) $\G$-$A$-$B$-correspondence $C_{A,B}^{\G}$ is $L^2(A)\otimes L^2(\G)\otimes L^2(B)$ endowed with 
\[
\pi(a) = (\pi_A\otimes \pi_M)(\alpha(a))_{12},\qquad \rho(b) =  (\rho_M \otimes \rho_B)(\beta^{\op}(b))_{23},
\]
\[
U = V_{24} \in B(L^2(A)\otimes L^2(\G)\otimes L^2(B) \otimes L^2(\G)).
\]
\end{Exa}
Here, we wrote $\beta^{\op}(b) = \varsigma \beta(b)\in M\bar{\otimes} B$ where $\varsigma: B\bar{\otimes}M\cong M\bar{\otimes} B$ is the flip.
The necessary compatibilities are satisfied by using the coaction property together with the identities \eqref{EqComultImpl} and \eqref{EqMultUnJhatJ}. The same reasoning works for the next examples.

\begin{Exa}\label{ExaSemiCoarsG}\label{Exa3}
Let $(\Hsp,U,\pi,\rho)$ be a $\G$-$A$-$B$-correspondence. We define the \emph{regular amplification} to be the $\G$-$A$-$B$-correspondence $S^{\G}(\Hsp) = \Hsp \otimes L^2(\G)$ endowed with 
\[
\pi_S(a) = (\pi\otimes \pi_M)(\alpha(a)),\qquad \rho_S(b) =  (\rho \otimes \rho_M)(\beta(b)),
\]
\[
U_S = V_{23} \in B(\mathcal{H}\otimes L^2(\G) \otimes L^2(\G)).
\]
It is a restriction of ${}_{A\bar{\otimes} M}(\Hsp \otimes L^2(\G))_{B\bar{\otimes}M}$ where $A\bar{\otimes}M,B\bar{\otimes}M$ are endowed with the coaction $\id\otimes \Delta$, and where we view $L^2(\G) = E_M$.  

In particular, we call $S_A^{\G} := S^{\G}(E_A)$ the \emph{$\G$-semi-coarse} (or $\G$-semi-regular) $\G$-$A$-$A$-correspondence.
\end{Exa}

Note that $S_A^{\G}$ becomes the trivial $A$-$A$-correspondence if $\G$ is trivial, and the regular $\G$-representation if $A$ is trivial. The following example reverses this situation. 

\begin{Exa}\label{ExaSemiCoarsW}\label{Exa4}
The \emph{W$^*$-semi-coarse} (or W$^*$-semi-regular) $\G$-$A$-$B$-correspondence $S_{A,B}^{\G}$ is $L^2(A)\otimes L^2(B)$ endowed with 
\[
\pi(a) = \pi_A(a)\otimes 1,\qquad \rho(b) =  1 \otimes \rho_B(b),
\]
\[
U = U_{\alpha,13}U_{\beta,23} \in B(L^2(A) \otimes L^2(B)\otimes L^2(\G)).
\]
\end{Exa}

\begin{Rem}\label{RemConjCoarse}
We can make the coarse $\G$-$A$-$B$-correspondence $C_{A,B}^{\G}$ resemble more the W$^*$-semi-coarse $\G$-$A$-$B$-correspondence $S^{\G}_{A,B}$ as follows: write 
\begin{equation}\label{EqAltUniCorepr}
\widetilde{U}_{B} = (u_{\G}\otimes 1)U_{\beta,21}^*(u_{\G}\otimes 1)  \in \rho_M(M) \bar{\otimes} B(L^2(B)),
\end{equation}
then from the identities \eqref{EqMultUnJhatJ} and \eqref{EqIdVW} we find that 
\[
\widetilde{U}_{B,12}^* V_{13} \widetilde{U}_{B,12} = V_{13}U_{\beta,23}.  
\]
So, conjugating the coarse $\G$-$A$-$B$-correspondence with the unitary $\widetilde{U}_{B,23}^*U_{\alpha,12}^*$, we find that $C_{A,B}^{\G}$ is equivalent to the $\G$-$A$-$B$-correspondence on $L^2(A)\otimes L^2(\G)\otimes L^2(B)$ given by 
\[
\pi(a) = \pi_A(a)\otimes 1\otimes 1,\qquad \rho(b) = 1\otimes 1\otimes \rho_B(b),
\]
\[
U = U_{\alpha,14}V_{24}U_{\beta,34}. 
\]
\end{Rem}

As a final example, we consider the following. 

\begin{Exa}\label{Exa5}
Let $(\Hsp,U,\pi,\rho)$ be a $\G$-$A$-$B$-correspondence. The \emph{adjoint amplification} is the $\G$-$A$-$B$-correspondence $D^{\G}(\Hsp) = \Hsp\otimes L^2(\G)$ endowed with 
\[
\pi_D(a) = (\pi\otimes \pi_M)(\alpha(a)),\qquad \rho_D(b) = \rho(b)\otimes 1,
\]
\[
U_D = V_{23}W_{32}U_{13}. 
\]
Note that $U_D$ is indeed a unitary $\G$-representation: it is the tensor product of $U_{\ad}$ with $U$, where $U_{\ad} = W_{21}V$ is the adjoint representation of $\G$ on $L^2(\G)$. The compatibility of $U_D$ and $\pi_D$ follows from the fact that $U$ implements $\alpha$, together with the easily computed commutation relation 
\begin{equation}\label{EqBraidedCommInt}
V_{23}W_{32} U_{13}U_{12} = U_{12}U_{13}V_{23}W_{32}. 
\end{equation}
The compatibility of $U_D$ and $\rho_D$ is immediate. 

In particular, we call $D_A^{\G} := D^{\G}(E_A)$ the 
\emph{adjoint} $\G$-$A$-$A$-correspondence. Using Example \ref{ExaAdjAct}, we see that $D_A^{\G}$ is simply the two-sided restriction to $A$ of the identity correspondence $E_{A\rtimes \G}$ for $A\rtimes \G$ with the adjoint $\G$-action. 
\end{Exa}

There are some obvious constructions one can perform on $\G$-equivariant correspondences, such as \begin{itemize}
\item taking (infinite) direct sums, and
\item taking $\G$-equivariant subcorrespondences.
\end{itemize} 
Let us make the latter more precise: if $(\Hsp,U,\pi,\rho)$ is a $\G$-$A$-$B$-correspondence, we define $\Gsp \subseteq \Hsp$ to be a  $\G$-equivariant subcorrespondence if $\Gsp$ is stable under $\pi(A)$, $\rho(B)$ and all $U(\omega)$ for $\omega \in M_*$. It follows from Proposition \ref{PropImRepC} that $\Gsp$ becomes a $\G$-$A$-$B$-correspondence $(\Gsp,U')$ in a unique way such that $U'(\omega)$ is the restriction of $U(\omega)$ to $\Gsp$. Moreover, this description also makes it clear that $\Gsp^{\perp}$ will be a $\G$-equivariant subcorrespondence, with $\Hsp = \Gsp \oplus \Gsp^{\perp}$ as $\G$-$A$-$B$-correspondences.

\begin{Def}
We call a $\G$-$A$-$B$-correspondence $(\Hsp,U,\pi,\rho)$ \emph{cyclic} if there exists $\xi \in \Hsp$ such that 
\[
\Hsp = [\pi(a)U(\omega)\rho(b)\xi\mid a\in A,\omega \in M_*,b\in B].
\]
\end{Def}

\begin{Prop}\label{PropDirectSumCycl}
Any $\G$-$A$-$B$-correspondence $(\Hsp,U,\pi,\rho)$ is a direct sum of cyclic $\G$-$A$-$B$-correspondences. 
\end{Prop}
\begin{proof}
Take $\xi \in\Hsp$. It is sufficient to show that 
\[
\Gsp:=  [\pi(a)U(\omega)\rho(b)\xi\mid a\in A,\omega \in M_*,b\in B]
\]
is a $\G$-equivariant subcorrespondence: since the C$^*$-algebra $C_U=[\{U(\omega)\mid \omega \in M_*\}]$ acts non-degenerately on $\Hsp$, we have $\xi\in \Gsp$, so $\Gsp$ will indeed be cyclic, and an appeal to Zorn's lemma then lets us conclude.

Clearly $\Gsp$ is invariant under $\pi(A)$. If $e_i$ is an orthonormal basis of $L^2(\G)$, then for $\xi,\eta\in L^2(\G)$ and $\omega = \omega_{\xi,\eta} \in B(L^2(\G))_*$, we have
\[
U(\omega)\pi(a) = (\id\otimes \omega)((\pi\otimes \id)\alpha(a)U) = \sum_i \pi((\id\otimes \omega_{\xi,e_i})(\alpha(a))) U(\omega_{e_i,\eta}), 
\]
the latter sum converging e.g.\ $\sigma$-strongly. Hence $\Gsp$ is invariant under the $U(\omega)$. 

Similarly, 
\[
\rho(b)U(\omega) = (\id\otimes \omega)(U (\rho \otimes R)(\beta(b))) = \sum_i U(\omega_{\xi,e_i}) \rho((\id \otimes \omega_{e_i,\eta}R)(\beta(b))),
\]
making it clear that $\Gsp$ is invariant under $\rho(B)$ as well. 
\end{proof} 

It follows easily from the above that we can make a \emph{universal} $\G$-$A$-$B$-correspondence $(\Hsp_u,U_u,\pi_u,\rho_u)$, in the sense that any cyclic $\G$-$A$-$B$-correspondence is unitarily equivalent to a direct summand of $\Hsp_u$. In the following, we fix such a chosen universal model.

\subsection{Associated C\texorpdfstring{$^*$}{star}-algebras}

\begin{Def}\label{DefSliceCat}
Let $A$ be a $\G$-W$^*$-algebra. We define 
\[
A_{\slice} = [(\id\otimes \omega)\alpha(a)\mid \omega \in M_*,a\in A] \subseteq A.
\]
\end{Def}

It is not clear if $A_{\slice}$ will always be a C$^*$-subalgebra of $A$, although it is easily seen to be an operator subsystem.\footnote{In case $\G$ is semi-regular, the proof of \cite{BSV03}*{Proposition 5.7} shows that $A_s$ will be a C$^*$-algebra, but it is not clear if the converse in that Proposition holds in the von Neumann algebraic setting.} In any case, we will not need any multiplicative structure of $A_{\slice}$. It will however be useful to observe that $A_s$ is $\sigma$-weakly dense in $A$ (see e.g. \cite{KaS15}*{Corollary 2.7} for an explicit proof). 

The following observation will be fundamental.

\begin{Theorem}\label{TheoExistCstarAlg}
  If $(\Hsp,U,\pi,\rho)$ is a $\G$-$A$-$B$-correspondence, then 
\begin{equation}\label{EqAnotherEq}
C^{\G}_{\Hsp}(A) := [C_U\pi(A)C_U],\qquad C^{\G}_{\Hsp}(B):= [C_U\rho(B)C_U],\qquad C^{\G}_{\Hsp}(A,B) := [C^{\G}_{\Hsp}(A)C^{\G}_{\Hsp}(B)]
\end{equation}
are C$^*$-algebras, and 
\begin{equation}\label{EqClosDefed}
C^{\G}_{\Hsp}(A) = [\pi(A_{\slice})C_U] = [C_U\pi(A_{\slice})],\quad C^{\G}_{\Hsp}(B) = [\rho(B_{\slice})C_U] = [C_U\rho(B_{\slice})].  
\end{equation}
\end{Theorem} 
\begin{proof}
Clearly $C^{\G}_{\Hsp}(A)$ and $C^{\G}_{\Hsp}(B)$ are closed under the $*$-operation. 

We now first prove \eqref{EqClosDefed}. Observe that by \eqref{EqIdSliceMaps}, we have 
\[
 [(\id\otimes \omega)(V)\mid \omega \in B(L^2(\G))_*] = 
[(\id\otimes \omega)(V^*)\mid \omega \in B(L^2(\G))_*].
\]
We then compute (dropping the notation $\pi$) 
\begin{eqnarray*}
C^{\G}_{\Hsp}(A) &=& [(\id\otimes \omega)(U_{12}a_1 U_{13})\mid a\in A,\omega \in B(L^2(\G)\otimes L^2(\G))_*]\\
&=&[(\id\otimes \omega)(\alpha(a)_{12}U_{12} U_{13})\mid a\in A,\omega \in B(L^2(\G)\otimes L^2(\G))_*]\\
&=& [(\id\otimes \omega)(\alpha(a)_{12}(\id\otimes \Delta)(U))\mid a\in A,\omega \in B(L^2(\G)\otimes L^2(\G))_*]\\
&=& [(\id\otimes \omega)(\alpha(a)_{12}V_{23}U_{12}V_{23}^*)\mid a\in A,\omega \in B(L^2(\G)\otimes L^2(\G))_*]\\
&=& [(\id\otimes \omega)(\alpha(a)_{12}V_{23}U_{12})\mid a\in A,\omega \in B(L^2(\G)\otimes L^2(\G))_*]\\
&=& [(\id\otimes \omega)(\alpha(a)_{12}V_{23}^*U_{12})\mid a\in A,\omega \in B(L^2(\G)\otimes L^2(\G))_*]\\
&=& [(\id\otimes \omega)(V_{23}^*(\id\otimes \Delta)\alpha(a)U_{12})\mid a\in A,\omega \in B(L^2(\G)\otimes L^2(\G))_*]\\
&=& [(\id\otimes \omega)((\alpha\otimes \id)\alpha(a)U_{12})\mid a\in A,\omega \in B(L^2(\G)\otimes L^2(\G))_*]\\
&=& [(\id\otimes \omega)(U_{12}\alpha(a)_{13})\mid a\in A,\omega \in B(L^2(\G)\otimes L^2(\G))_*]\\
&=& [C_UA_{\slice}]. 
\end{eqnarray*}
Since $C_{\Hsp}^{\G}(A)$ is closed under the $*$-operation, we obtain \eqref{EqClosDefed} for $A$. It is now clear that $C_{\Hsp}^{\G}(A)$ is a C$^*$-algebra: 
\[
[C_{\Hsp}^{\G}(A)C_{\Hsp}^{\G}(A)] = [C_UA_{\slice}^2C_U] \subseteq C_{\Hsp}^{\G}(A). 
\]

For $B$, we can simply observe that $B^{\opp}$ has a right coaction by $(M,\Delta^{\opp})$ via 
\[
b^{\opp} \mapsto \beta(b)^{\opp\otimes R},
\]
and that $U^*$ is a unitary representation for $(M,\Delta^{\opp})$ with
\[
U^*(\rho(b)\otimes 1)U = (\rho\otimes R)\beta(b),\qquad b\in B. 
\]
Hence the result for $B$ follows from that for $A$ (since $C_U = C_{U^*}$).

Finally, $C^{\G}_{\Hsp}(A,B)$ is closed under the $*$-structure since 
\[
C^{\G}_{\Hsp}(A,B)^* = [C_U\rho(B_{\slice})\pi(A_{\slice})C_U] = [C_U\pi(A_{\slice})\rho(B_{\slice})C_U] = [C_U\pi(A)C_UC_U\rho(B)C_U] = C^{\G}_{\Hsp}(A,B).
\]
A similar computation shows that $C^{\G}_{\Hsp}(A,B)$ is closed under multiplication. 
\end{proof}

\begin{Rem}\label{RemIdentIntoQM}
It also follows immediately from \eqref{EqAnotherEq} and \eqref{EqClosDefed} that for example
\[
C^{\G}_{\Hsp}(A) = [C_U\pi(A_{\slice})C_U], \qquad C^{\G}_{\Hsp}(A,B) = [\pi(A_{\slice})C_U\rho(B_{\slice})],
\]
and that 
\[
\pi(A) \subseteq QM(C_{\Hsp}^{\G}(A,B)),\qquad \rho(B) \subseteq QM(C_{\Hsp}^{\G}(A,B)). 
\]
This appearance of the space of quasi-multipliers is somewhat surprising and mysterious, but seems unavoidable as the following remark shows. 
\end{Rem}
\begin{Rem}\label{example quasimultiplier}
In general, we will not have that $\pi(A) \subseteq M(C_{\Hsp}^{\G}(A,B))$. A concrete example is obtained as follows: let $A = B(\ell^2(\Z\times \Z))$, let $B = \C$ and let $\G = \T = \{z\in \C\mid |z|=1\}$. Let $e_{m,n}$ be the standard basis of $\Hsp = \ell^2(\Z\times\Z)$, and let $\T$ act on $\Hsp$ via 
\[
U(z)e_{n,m} = z^ne_{n,m},\qquad z\in \T,m,n\in \Z.
\]
We also view $U$ as a unitary in $M(\mathcal{K}(\mathcal{H})\otimes C(\mathbb{T}))\subseteq B(\mathcal{H})\bar{\otimes} L^\infty(\mathbb{T})$ through 
$$(\id \otimes \ev_z)(U) = U(z), \quad z \in \mathbb{T}.$$
Endow $B$ with the trivial action, and endow $A$ with the inner $\T$-action 
\[
\alpha: A \to A \bar{\otimes} L^\infty(\mathbb{T}): x \mapsto U(x\otimes 1)U^*.
\]
This equips $\Hsp$ with the structure of a $\T$-equivariant $A$-$B$-correspondence, and $C^{\T}_{\Hsp}(A,B) = C^{\T}_{\Hsp}(A)$. 

We write the irreducible unitary representations of $\mathbb{T}$ as
$$\chi_n: \mathbb{T}\to \mathbb{C}: z \mapsto z^n, \quad n \in \mathbb{Z}.$$

Consider now the operator $x\in A$ determined by 
\[
xe_{n,m} = \delta_{n,0} e_{m,0}. 
\]
We claim that $x\notin M(C^{\T}_{\Hsp}(A))$. Indeed, let $p_n$ be the projection onto $\ell^2(\{n\}\times \Z)$ for $n\in \mathbb{Z}$. 
Viewing $\chi_n \in L^1(\mathbb{T})= L^\infty(\mathbb{T})_*$, we see that 
$(\id \otimes \chi_n)(U) = p_{-n}$ so that $C_U = [p_n: n \in \mathbb{Z}]$ by the Stone-Weierstrass theorem.
It is then sufficient to prove that $x=xp_0 \notin C^{\T}_{\Hsp}(A) = [A_sC_U]$, or equivalently that $x \notin [A_sp_0]$. However, one easily verifies that the spectral subspace of $B(\mathcal{H})$ associated to the irreducible representation $\chi_m$ is given by
$$A_m= \{y \in B(\Hsp) \mid \forall n \in \Z: yp_n = p_{n+m}yp_n\}.$$
Consequently\[
A_s = [A_m\mid m\in \Z]
\]
It is elementary to see from this description that $\lim_{n\rightarrow \infty} \|p_ny\|=0$ for all $y \in [A_sp_0]$. Since $\|p_nx\|=1$ for all $n \in \Z$, we find that indeed $x \notin [A_sp_0]$. 
\end{Rem} 

\begin{Theorem}\label{TheoMMM}
There is a well-defined semi-norm on $A_{\slice}\odot M_* \odot B_{\slice}$ through 
\begin{equation}\label{DefUniNorm} 
\|\sum_{i=1}^n a_i \otimes \omega_i\otimes b_i\|_u = \sup\{\|\sum_{i=1}^n \pi(a_i)U(\omega_i)\rho(b_i)\|\mid (\Hsp,U)\textrm{ a }\G\textrm{-}A\textrm{-}B\textrm{-correspondence}\}.
\end{equation}
Moreover, the associated separation-completion $C^{\G}(A,B)$ has the unique structure of a C$^*$-algebra such that, for each $\G$-$A$-$B$-correspondence $(\Hsp,U,\pi,\rho)$, the map
\begin{equation}\label{EqThetaFromCorr}
\theta = \theta_{\Hsp}: C^{\G}(A,B) \rightarrow B(\Hsp),\quad a\otimes \omega \otimes b \mapsto \pi(a)U(\omega)\rho(b)
\end{equation}
is a well-defined non-degenerate $*$-representation. 
\end{Theorem} 
\begin{proof}
Clearly $\|\textrm{---}\|_u$ is well-defined as a semi-norm, with 
\[
\|\sum_{i=1}^n a_i \otimes \omega_i\otimes b_i\|_u \leq \sum_{i=1}^n \|a_i\|\|\omega_i\|\|b_i\|. 
\]
Now by Proposition \ref{PropDirectSumCycl}, the supremum in \eqref{DefUniNorm} does not change if we consider only cyclic $\G$-equivariant $A$-$B$-correspondences. Hence we see that $\|-\|_u$ is in fact realized by the universal $\G$-equivariant $A$-$B$-correspondence $(\Hsp_u,U_u,\pi_u,\rho_u)$, i.e.\  $C^{\G}(A,B) := C_{\Hsp_u}^{\G}(A,B)$ satisfies the requirements. 
\end{proof}

In the following, we identify $C^{\G}(A,B) = C_{\Hsp_u}^{\G}(A,B)$ whenever convenient. We do have to be a bit careful when using this identification, as we want to make sure that constructions do not actually depend on the precise model chosen. Here is a case in point: By means of Remark \ref{RemIdentIntoQM}, we have a natural map $\pi_{\Hsp_u}: A \rightarrow QM(C_{\Hsp_u}^{\G}(A,B))$, which through the isomorphism 
\[
\theta_{\Hsp_u}: C^{\G}(A,B) \cong C^{\G}_{\Hsp_u}(A,B)
\]
we can pull back to a map

\begin{equation}\label{EqUnivIncl}
\pi_u: A \rightarrow QM(C^{\G}(A,B)).
\end{equation}
We want to make sure for example that this map does not depend on the chosen model.

\begin{Lem}\label{LemIndepMod}
The map \eqref{EqUnivIncl} is independent of the chosen model. 
\end{Lem}
\begin{proof}
It is enough to show the following: if $(\Hsp,U,\pi,\rho)$ is any $\G$-equivariant $A$-$B$-correspondence, then (using the same symbol to denote extensions of non-degenerate $*$-homomorphisms to their quasi-multipliers)
\[
\theta_{\Hsp}(\theta_{\Hsp_u}^{-1}(\pi_u(a))) = \pi_{\Hsp}(a),\qquad \forall a\in A. 
\]
Since any $\Hsp$ is a direct sum of cyclic ones, it is enough to assume that $\Hsp$ is cyclic as a $\G$-$A$-$B$-correspondence. But then since $\Hsp_u$ contains any cyclic $\G$-$A$-$B$-correspondence up to unitary equivalence, we may as well assume that $\Hsp \subseteq \Hsp_u$. It is now easily seen that if $x \in QM(C^{\G}(A,B))$, then 
\[
\theta_{\Hsp}(x) = \theta_{\Hsp_u}(x)_{\mid \Hsp}.
\]
In particular, taking $x = \pi_u(a)$ for $a\in A$, we see that
\[
\theta_{\Hsp}(\theta_{\mathcal{H}_u}^{-1}(\pi_u(a))) = \pi_{\Hsp_u}(a)_{\mid \Hsp} = \pi_{\Hsp}(a).  
\]\end{proof}
Similarly, one shows that there is a canonical map 
\[
\rho_u: B \to QM(C^\G(A,B)).
\]

\begin{Rem}
\begin{itemize}
\item We do not know if the map $\pi_u$ in \eqref{EqUnivIncl} is multiplicative from $A$ to $C^{\G}(A,B)^{**}$. Indeed, it is not clear if elements of $C^{\G}(A,B)^{**}$ can be separated by (normal extensions of) $*$-representations $\theta$ coming from $\G$-$A$-$B$-correspondences, while it is only in such representations that multiplicativity statements for $A$ can be deduced.
\item On the other hand, it follows from Theorem \ref{TheoExistCstarAlg} that $\pi_u(A_{\slice}) \subseteq M(C^{\G}(A,B))$, and hence the same is true for the C$^*$-algebra generated by $A_{\slice}$. In particular, it \emph{is} hence true that $\pi_u(ab) = \pi_u(a)\pi_u(b)$ for $a,b\in A_{\slice}$.
\end{itemize}
\end{Rem}

We can now deduce a universal property characterizing $*$-representations of $C^{\G}(A,B)$ coming from equivariant correspondences. Let us first consider two extreme cases. 

\begin{Lem}
Assume that $\G$ is the trivial group. Then 
\[
C^{\G}(A,B) \cong A \otimes_{\mathrm{bin}}B^{\opp}
\]
where $\otimes_{\operatorname{bin}}$ denotes the binormal tensor product of von Neumann algebras \cite{EL77}.
\end{Lem}
\begin{proof}
This follows immediately from the discussion in \cite{EL77}*{Section 2}, together with the fact that an $A$-$B$-correspondence is a direct sum of cyclic ones. 
\end{proof}

\begin{Lem}
Assume that $A = B = \C$. Then 
\[
C^{\G}(\C,\C) \cong C^*(\G), 
\]
the usual universal group C$^*$-algebra of a locally compact quantum group \cite{Kus01}.
\end{Lem}
Note: since our conventions on corepresentations are the opposite of the one in \cite{Kus01}, we tacitly modify the definition of $C^*(\G)$ below correspondingly. 
\begin{proof}
Recall from Proposition \ref{PropImRepC} that any unitary $\G$-representation $(\Hsp,U)$ is continuous, i.e.\ 
\[
U \in M(\mcK(\Hsp) \otimes C_0^{\red}(\G)). 
\]

We now recall briefly the construction of $C^*(\G)$. Let $S$ be the (unbounded) antipode for $C_0^{\red}(\G)$, with domain $\mcD(S)$. For $\omega \in L^1(\G)$, denote 
\[
\omega^*: \mcD(S) \rightarrow \C,\quad x \mapsto \overline{\omega(S(x)^*)}.
\]
Then following \cites{Kus01, KV03}, we denote 
\[
L^1_*(\G) := \{\omega \in L^1(\G) \mid \exists \omega' \in L^1(\G) \textrm{ with } \omega'_{\mid \mcD(S)} = \omega^*\}.
\]
It then follows from \cite{Kus01}*{Proposition 5.2} that $C^*(\G)$ may be identified with the universal completion of $L^1_*(\G)$ for the norm 
\[
\|\omega\|_u = 
\sup\{\|(\id\otimes \omega)U\|\mid U \textrm{ unitary }\G\textrm{-representation}\}.
\]
It now follows by the universal properties of the objects involved, that we obtain an isometric inclusion of C$^*$-algebras
\begin{equation}\label{EqIncl}
C^*(\G)\rightarrow C^{\G}(\C,\C),\qquad L^1_*(\G)\ni \omega \mapsto \omega \in L^1(\G).
\end{equation}
It only remains to observe that, through a smoothing argument with respect to the scaling group, the space $L^1_*(\G)$ is norm-dense in $L^1(\G)$ (for the usual pre-dual norm). Hence \eqref{EqIncl} must be an isomorphism of C$^*$-algebras.
\end{proof}

In the following, we will simply identify $C^*(\G) = C^{\G}(\C,\C)$. Then taking a universal model for $C^*(\G)$, we obtain the associated universal multiplicative unitary 
\[
V_u \in M(C^*(\G) \otimes C_0^{\red}(\G))
\]
as in \cite{Kus01}*{Proposition 5.2}. Then for $A,B$ any pair of $\G$-dynamical von Neumann algebras, we have a unique non-degenerate C$^*$-homomorphism 
\[
\theta_{\G}: C^*(\G)\rightarrow C^{\G}(A,B),\qquad (\id\otimes \omega)V_u \mapsto (\id\otimes \omega)U_u,\qquad \omega \in M_*. 
\]

\begin{Prop}
Let $\theta$ be a non-degenerate $*$-representation of $C^{\G}(A,B)$. Then $\theta = \theta_{\Hsp}$ for a $\G$-$A$-$B$-correspondence $(\Hsp,U,\pi,\rho)$ if and only if 
\begin{enumerate}
\item\label{EqCond1Char} $\theta \circ \theta_{\G}$ is a non-degenerate $*$-representation of $C^*(\G)$, 
\item\label{EqCond2Char} $\theta \circ \pi_u$ is a normal unital $*$-representation of $A$, and 
\item\label{EqCond3Char} $\theta\circ \rho_u$ is a normal unital $*$-anti-representation of $B$. 
\end{enumerate}
\end{Prop} 
\begin{proof}
If $\theta = \theta_{\Hsp}$, a similar argument as in Lemma \ref{LemIndepMod} shows that $\theta$ satisfies 
\[
\theta(\pi_u(a)) = \pi(a),\quad \theta(\rho_u(b))  = \rho(b),\quad \theta((\id\otimes \omega)U_u) = (\id\otimes \omega)U,\qquad a\in A,b\in B,\omega \in M_*,
\]
so $\theta$ satisfies the conditions \eqref{EqCond1Char}, \eqref{EqCond2Char} and \eqref{EqCond3Char}.

Assume now conversely that $(\Hsp,\theta)$ satisfies the conditions \eqref{EqCond1Char}, \eqref{EqCond2Char} and \eqref{EqCond3Char}. Write 
\[
\pi = \theta \circ \pi_u,\quad \rho = \theta\circ \rho_u,
\]
and let $U$ be the unitary $\G$-representation associated to the non-degenerate $*$-representation $\theta\circ \theta_{\G}$ of $C^*(\G)$. Then we claim that
\begin{enumerate}[label=(\alph*)]
\item\label{EqNeedToProve1} $\pi$ and $\rho$ commute elementwise. 
\item\label{EqNeedToProve2} $U$ implements $(\pi\otimes \id)\alpha$ on $\pi(A)$.
\item\label{EqNeedToProve3} $U^*$ implements $(\rho \otimes R)\beta$ on $\rho(B)$. 
\end{enumerate}
This will make $(\Hsp,U,\pi,\rho)$ into a $\G$-equivariant $A$-$B$-correspondence, and the equality $\theta = \theta_{\Hsp}$  follows from 
\[
\theta(a\otimes \omega\otimes b) = \theta(\pi_u(a)\pi_{\G}((\id\otimes \omega)(V_u))\rho_u(b)) = \pi(a)(\id\otimes \omega)(U)\rho(b),\qquad a\in A_{\slice},b\in B_{\slice},\omega \in M_*.
\]

Let us tackle first \ref{EqNeedToProve1}.

Since $A_{\slice}$ is $\sigma$-weakly dense in $A$, and $B_s$ is $\sigma$-weakly dense in $B$, it is by the assumption of normality of $\pi$ and $\rho$ sufficient to prove that $\pi(A_s) = \theta(\pi_u(A_{\slice}))$ and $\rho(B_s) = \theta(\rho_u(B_{\slice}))$ commute elementwise. But $\pi_u(A_{\slice}),\rho_u(B_{\slice})$ lie in the multiplier algebra of $C^{\G}(A,B)$, so we can check the elementwise commutation of $\pi_u(A_{\slice})$ and $\rho_u(B_{\slice})$ in a universal model, where this is obvious. 

Let us now tackle the second item (the third item is similar).

Again by normality of $\pi$, it is sufficient to show that 
\begin{equation}\label{EqToProveEquivv}
U(\pi(a)\otimes 1)U^* = (\pi\otimes \id)(\alpha(a)) ,\qquad a\in A_{\slice}. 
\end{equation}
Then since $\pi_u(A_s) \subseteq M(C^{\G}(A,B))$ and 
\[
U_u = (\theta_{\G}\otimes \id)V_u \in M(C^{\G}(A,B)\otimes C_0^{\red}(\G)), 
\]
we have that 
\[
(\pi_u \otimes \id)(\alpha(a))  = U_u(\pi_u(a)\otimes 1)U_u^* \in M(C^{\G}(A,B)\otimes C_0^{\red}(\G)),\qquad a\in A_{\slice},
\]
using a universal model $(\Hsp_u,U_u,\pi_u,\rho_u)$ to make sense of the left hand side as an operator on $\Hsp_u \otimes L^2(\G)$ a priori. Applying $\id\otimes \omega$ for $\omega\in M_*$, we then find for $a\in A_{\slice}$ and $\omega \in M_*$ that
\[
\pi((\id\otimes \omega)\alpha(a)) = \theta(\pi_u((\id\otimes \omega)\alpha(a))) = \theta((\id\otimes \omega)(U_u(\pi_u(a)\otimes 1)U_u^*)) = (\id\otimes \omega)(U(\pi(a)\otimes 1)U^*).
\]
The identity \eqref{EqToProveEquivv} follows immediately from this.
\end{proof}

\subsection{Fell topology}

For $C$ a C$^*$-algebra, we denote by $\Rep_*(C)$ the W$^*$-category of its non-degenerate representations.

The following lemma states that the functor $\Corr^{\G}(A,B)\rightarrow \Rep_*(C^{\G}(A,B))$ (determined by \eqref{EqThetaFromCorr}) is full.

\begin{Lem}\label{LemMorphDoNotChange}
If $(\Hsp,U,\pi,\rho),(\widetilde{\Hsp},\widetilde{U},\widetilde{\pi},\widetilde{\rho})$ are $\G$-$A$-$B$-correspondences, then an element $x \in B(\Hsp,\widetilde{\Hsp})$ is an intertwiner of $\G$-$A$-$B$-correspondences if and only if it is an intertwiner of $C^{\G}(A,B)$-correspondences. 
\end{Lem} 
\begin{proof}
Clearly an intertwiner of $\G$-$A$-$B$-correspondences is also an intertwiner of $C^{\G}(A,B)$-representations, by construction of $\theta_{\Hsp}$. Conversely, if $x$ intertwines $\theta_{\Hsp}$ and $\theta_{\widetilde{\Hsp}}$, we have (again by construction) that
\[
x \pi(a)U(\omega)\rho(b) = \widetilde{\pi}(a)\widetilde{U}(\omega)\widetilde{\rho}(b)x,\qquad \forall a\in A_s,\omega\in M_*,b\in B_s. 
\]Taking $a= 1$ and $b=1$, we find $xU(\omega)= \widetilde{U}(\omega)x$ for all $\omega \in M_*$. Also,  (taking $b=1$)
$$x\pi(a)U(\omega) = \widetilde{\pi}(a)\widetilde{U}(\omega)x = \widetilde{\pi}(a)xU(\omega), \qquad \forall a \in A_s, \omega \in M_*$$
so that by non-degeneracy of $C_U$ we find $x\pi(a)= \widetilde{\pi}(a)x$ for $a\in A_s$. Since $A_s$ is $\sigma$-weakly dense in $A$, we find $x\pi(a) = \widetilde{\pi}(a)x$ for $a\in A$. Similarly $x\rho(b) = \widetilde{\rho}(b)x$ for $b\in B$. Thus, $x$ is an intertwiner of $\G$-$A$-$B$-correspondences.
\end{proof}

We now define a topology on the objects of $\Corr^{\G}(A,B)$. Since objects in a category $\mcC$ such as this do not form a set, we need to be careful in specifying what we mean by this. We use the formalism that a topology for objects in a category $\mcC$ is determined by a  \emph{set-indexed local open neighborhood system}: to each object $X\in \mcC$ is assigned a directed \emph{set} $\msF_X$ and, to each $F \in \msF_X$, a collection of objects $\mcC_F(X)$ in $\mcC$, such that 
\begin{enumerate}
\item $X$ is an object of $\mcC_F(X)$ for each $F\in \msF_X$,
\item $\mcC_{F'}(X)$ is contained in $\mcC_{F}(X)$ if $F\leq F'$, and
\item if $Y$ is an object of some $\mcC_F(X)$, there exists $F'\in \msF_Y$ such that $\mcC_{F'}(Y)$ is contained in $\mcC_F(X)$. 
\end{enumerate}
We then say that two such set-indexed local open neighborhood systems determine the same topology $\tau$, if at each object $X$ the respective open neighborhood systems can be nested in each other.

\begin{Def}\label{DefFellTop}
Let $(\Hsp,U,\pi,\rho)$ be a $\G$-$A$-$B$-correspondence. Let $F$ be a collection of a scalar $\varepsilon>0$ and finite subsets $F_A,F_B,F_{\G}$ and $F_{\Hsp} = \{\xi_1,\ldots,\xi_n\}$ in resp.\ $A,B,M_*$ and $\Hsp$. 

We say that a $\G$-$A$-$B$-correspondence $(\Hsp',U',\pi',\rho')$ lies in the \emph{$F$-neighborhood} $U(\Hsp;F)$ of $(\Hsp,U,\pi,\rho)$ if there exist $\xi_1',\ldots,\xi_n'$ in $\Hsp'$ with $\|\xi_i'\|\leq \|\xi_i\|$ for all $i$, such that, for all $a\in F_A, b\in F_B,\omega\in  F_{\G}$ and $1\leq i,j\leq n$, 
\begin{equation}\label{EqFundNeighb}
|\langle \xi_i,\pi(a)U(\omega)\rho(b)\xi_j\rangle - \langle\xi_i', \pi'(a)U'(\omega)\rho'(b)\xi_j'\rangle |<\varepsilon. 
\end{equation}
\end{Def}
Such $F$ as above form a directed set in a natural way, by putting 
\[
F \leq F' \quad \iff\quad \varepsilon \geq \varepsilon', F_A \subseteq F_A', F_B\subseteq F_B',F_{\G}\subseteq F_{\G}',F_{\Hsp}\subseteq F_{\Hsp}'.
\]
It is then easily verified that the neighborhoods defined in Definition \ref{DefFellTop} define a set-indexed local open neighborhood system, and hence a topology $\tau$ on $\Corr^{\G}(A,B)$. Moreover, for any finite collection of $\G$-equivariant correspondences $\Hsp_1,\ldots,\Hsp_n$ and neighborhoods $U(\Hsp_i;F_i)$, there exists $F$ such that  \begin{equation}\label{EqIntersect}
U(\oplus_{i=1}^n \Hsp_i;F) \subseteq \cap_{i=1}^n U(\Hsp_i;F_i).
\end{equation}

\begin{Rem}
The ordering of the factors in \eqref{EqFundNeighb} is important: in general, one seems to get a different (and less natural) notion if for example the positions of $U(\omega)$ and $\pi(a)$ are interchanged. However, one may interchange the positions of $\pi(a)$ and $\rho(b)$ at no cost, as will follow easily from Theorem \ref{TheoEquivTop}. 
\end{Rem}

Recall now from Theorem \ref{TheoMMM} that any $\G$-$A$-$B$-correspondence $\Hsp$ leads to a (unique) non-degenerate $*$-representation $(\Hsp,\theta_{\Hsp})$ of $C^{\G}(A,B)$. We will show that the $\tau$-topology on $\Corr^{\G}(A,B)$ is simply the restriction of the corresponding `Zimmer topology' on the W$^*$-category $\Rep_*(C^{\G}(A,B))$ \cite{Zim84} (see also \cite{BdlHV08}{Appendix F}), passing from the group case to the C$^*$-algebra setting.

\begin{Def}\label{DefFellC}
Let $C$ be a C$^*$-algebra, and let $(\Hsp,\theta)$ be a non-degenerate $*$-representation of $C$. Let $F$ consist of $\varepsilon>0$ and finite sets $F_C \subseteq C$ and $F_{\Hsp} = \{\xi_1,\ldots,\xi_n\}\subseteq \Hsp$. 

If $(\Hsp',\theta')$ is a non-degenerate $*$-representation of $C$, we say that $(\Hsp',\theta')$ lies in the $F$-neighborhood $V(\Hsp;F)$ of $(\Hsp,\theta)$ if there exist $\xi_1',\ldots,\xi_n' \in \Hsp'$ such that $\|\xi_i'\|\leq \|\xi_i\|$ for all $i$, and, for all $c\in F_C$ and $1\leq i,j\leq n$, 
\[
|\langle \xi_i,\theta(c)\xi_j\rangle - \langle \xi_i',\theta'(c)\xi_j'\rangle |<\varepsilon. 
\]
\end{Def}

We denote the resulting topology on $\Rep_*(C)$ by $\sigma$. 

\begin{Rem}
By means of a contractive approximate unit in $C$, it is easily verified that if in Definition \ref{DefFellC} we drop the requirement that $\|\xi_i'\|\leq \|\xi_i\|$ for all $i$, the resulting topology coincides with $\sigma$. We do not know if this still holds for Definition \ref{DefFellTop}.
\end{Rem} 

Note now that the definition of $V(\Hsp;F)$ makes sense if $F_C \subseteq C^{**}$, using the unique normal extension of $\theta$ to $C^{**}$. We will then use the following general fact. 

\begin{Lem}\label{LemFellTop}
Let $S \subseteq QM(C)$ be a subset such that $C \subseteq [S]$. Let $\sigma_S$ be the topology on $\Rep_*(C)$ obtained by using as local neighborhoods the $V(\Hsp;F)$ for finite $F_C\subseteq S$. Then 
\[
\sigma_S= \sigma.
\]
\end{Lem}

\begin{proof}
Let $(\Hsp,\theta)\in \Rep_*(C)$ and $F = (\varepsilon,F_C,F_{\Hsp} = \{\xi_1,\ldots,\xi_n\})$ with $F_C \subseteq S$. Choose $d \in C$ with $\|d\|\leq 1$ and 
\[
\max_{i}\{\|\xi_i -\theta(d)\xi_i\|\}\max_{c\in F_C}\{ \|c\|\}\max_{j}\{\|\xi_j\|\}<\varepsilon/4.
\] 
Then for all $1\leq i,j\leq n$ and all $c\in F_C$, we have 
\[
|\langle \xi_i,\theta(c)\xi_j\rangle - \langle \theta(d)\xi_i,\theta(c)\theta(d)\xi_j\rangle| <\varepsilon/2.
\]

If now $\widetilde{F} = (\varepsilon/2,\widetilde{F}_C,F_{\Hsp})$ with $\widetilde{F}_C = \{d^*cd\mid c\in F_C\} \subseteq C$, then $V(\Hsp;\widetilde{F}) \subseteq V(\Hsp;F)$. Indeed, if $(\Hsp',\theta')\in V(\Hsp;\widetilde{F})$, we can choose $\xi_1',\ldots,\xi_n' \in \Hsp'$ with $\|\xi_i'\|\leq \|\xi_i\|$ such that for all $1\leq i,j\leq n$ and $c \in F_C$,
\[
|\langle \theta(d)\xi_i,\theta(c)\theta(d)\xi_j\rangle - \langle \theta'(d)\xi_i',\theta'(c)\theta'(d)\xi_j'\rangle|<\varepsilon/2, 
\]
and then clearly 
\[
|\langle \xi_i,\theta(c)\xi_j\rangle - \langle \theta'(d)\xi_i',\theta'(c)\theta'(d)\xi_j'\rangle|<\varepsilon. 
\]
In particular, this shows that 
\[
\sigma_S \subseteq  \sigma.
\]

We now claim conversely that $\sigma_S$ is finer than $\sigma$. In fact, we have the following easy observations: the topology $\sigma_S$ does not change if we add to $S$ 
\begin{itemize}
\item scalar multiples of elements of $S$,
\item sums of elements of $S$, or 
\item normlimits of elements of $S$.
\end{itemize}
In particular, we may replace $S$ by $[S]$, and so the assumption $[S]\supseteq C$ immediately gives the claim. 
\end{proof} 
 
\begin{Theorem}\label{TheoEquivTop}
Let $\G$ be a locally compact quantum group. The topology $\tau$ on $\Corr^{\G}(A,B)$ is induced by the corresponding topology $\sigma$ on $\Rep_*(C^{\G}(A,B))$ under the correspondence $(\Hsp,U,\pi,\rho) \mapsto (\Hsp,\theta_{\Hsp})$. 
\end{Theorem}
\begin{proof}
Observing that 
\[
C^{\G}(A,B) = [\pi_u(A_s)C_{U_u}\rho_u(B_s)] \subseteq [\pi_u(A)\{(\id\otimes \omega)U_u\mid \omega \in M_*\}\rho_u(B)] \subseteq QM(C^{\G}(A,B)),
\] 
this follows immediately from Lemma \ref{LemFellTop}.
\end{proof}

The topology $\tau$ itself is not the most useful: rather, we want to consider the \emph{Fell topology}\footnote{It is an instructive exercise to check that this indeed determines a well-defined topology, in the sense described earlier.} $\tau^{\infty}$ which has local neighorhoods
\[
U^{\infty}(\Hsp;F) = \{\Gsp \mid \Gsp^{\infty} \in U(\Hsp;F)\}
\]
where here and in the following we use the shorthand notation 
\[
\Gsp^{\infty} := \oplus_1^{\infty}\Gsp.
\]
The corresponding topology $\sigma^{\infty}$ on $\Rep_*(C^{\G}(A,B))$ indeed coincides with the usual Fell (=inner hull-kernel) topology \cite{Fel62}.

\section{Weak containment}

In this section, we fix again a locally compact quantum group $\G$ and $\G$-W$^*$-algebras $A,B$.

\begin{Def}
Let $\Hsp,\Hsp'$ be $\G$-$A$-$B$-correspondences. We say that 
$\Hsp'$ is weakly contained in $\Hsp$ iff $\Hsp'$ lies in the $\tau^{\infty}$-closure of $\{\Hsp\}$.
\end{Def} 
We then write
\[
\Hsp' \preccurlyeq \Hsp. 
\]
In other words, $\Hsp' \preccurlyeq \Hsp$ if and only if for any $F$ as in Definition \ref{DefFellTop}, we have $\oplus_1^{\infty} \Hsp \in U(\Hsp';F)$. 

Clearly, weak containment is a transitive and reflexive relation. We also note that in terms of the C$^*$-algebra $C^{\G}(A,B)$, we have $\Gsp \preccurlyeq \Hsp$ if and only if the representation $\theta_{\Gsp}$ of $C^{\G}(A,B)$ factors through $\theta_{\Hsp}(C^{\G}(A,B))$ (see e.g.\ \cite{Dix77}*{Section 3.4}).  

\begin{Rem}\label{RemMulti}
Let $I,J$ be arbitrary sets, and let $\Hsp,\Hsp'$ be $\G$-$A$-$B$-correspondences. Then it follows immediately from the above description that 
\[
\Hsp' \preccurlyeq \Hsp\qquad \Leftrightarrow \qquad \oplus_{j\in J} \Hsp' \preccurlyeq \oplus_{i\in I} \Hsp, 
\]
i.e.\ weak containment does not witness multiplicities. More intrinsically, we can write this as
\[
\Hsp' \preccurlyeq \Hsp\qquad \Leftrightarrow \qquad   l^2(J)\otimes\Hsp'  \preccurlyeq  l^2(I) \otimes\Hsp, 
\]
and we may then as well replace $l^2(I),l^2(J)$ with arbitrary multiplicity Hilbert spaces $\Gsp,\Gsp'$.
\end{Rem}

In this section, we will provide a different characterization of weak containment. This will generalize to the equivariant setting some of the techniques and ideas of \cites{A-D90,A-DH90,A-D95}. 

Before we continue, we prove a small preparatory lemma on von Neumann algebras and nsf weights. Although it must be well-known, we do not know a convenient reference, so we briefly provide a proof. 

Recall that if $B$ is a von Neumann algebra and $\varphi$ an nsf weight on $B$, we can endow the opposite von Neumann algebra $B^{\opp}$ with the nsf weight 
\[
\varphi^{\opp}: B^{\opp}_+\rightarrow [0,+\infty],\quad b^{\opp}\mapsto \varphi(b).
\]
Then $\mathscr{N}_{\varphi^{\opp}} = (\mathscr{N}_{\varphi}^*)^{\opp}$, and we can identify $L^2(B)$ with $L^2(B^{\opp})$ via 
\[
\Lambda_{\varphi^{\opp}}(b^{\opp})= J\Lambda_{\varphi}(b^*),\qquad b\in \mathscr{N}_{\varphi}^*. 
\]

\begin{Lem}\label{LemPrelResvN}
Let $B$ be a von Neumann algebra. If $\xi_1,\ldots,\xi_n \in L^2(B)$, there exists an nsf weight $\varphi$ on $B$, a vector $\xi\in L^2(B)$ and $b_1,\ldots, b_n \in B$ such that $b_i^{\opp} \in \mathscr{N}_{\varphi^{\opp}}$ for all $i$ and 
\begin{equation}\label{EqIdWeight}
\xi_i =\rho_B(b_i)\xi = \Lambda_{\varphi^{\opp}}(b_i^{\opp}),\qquad \forall i.
\end{equation}
\end{Lem} 
\begin{proof}
Let $\Hsp =[\rho_B(B)\xi_1 + \ldots + \rho_B(B)\xi_n]$. Then $\Hsp = \pi_B(p)L^2(B)$ for some selfadjoint projection $p\in B$. Put $q= 1-p$, and let $\psi$ be an nsf weight on $qBq$. Let 
\[
\omega: B \rightarrow \C,\qquad b\mapsto \sum_{i=1}^n \langle \xi_i,\pi_B(b)\xi_i\rangle,
\]
and put
\[
\varphi(b) = \omega(pbp) + \psi(qbq),\qquad b\in B_+.
\]
Then clearly $\varphi$ is an nsf weight on $B$, with $p$ integrable. Moreover, $p$ is invariant under the modular automorphism group of $\varphi$, so we can define
\[
\xi := \Lambda_{\varphi}(p) = \Lambda_{\varphi^{\opp}}(p^{\opp}). 
\]
Now an easy calculation gives that 
\[
\langle \xi,\pi_B(b)\xi\rangle = \omega(b),\qquad b\in B. 
\]
Recalling the notation \eqref{EqRightRep}, it follows that there exists $b_i \in B$ such that $\rho_B(b_i)\pi_B(b)\xi = \pi_B(b)\xi_i$. In particular, $\rho_B(pb_i) \xi = \xi_i$, and clearly $\rho_B(pb_i)\xi = \Lambda_{\varphi^{\opp}}((pb_i)^{\opp})$. 
\end{proof} 

Assume now that $\Hsp,\Gsp$ are Hilbert spaces with a normal anti-$*$-representation $\rho_{\Hsp},\rho_{\Gsp}$ by the von Neumann algebra $B$. Recall the notation \eqref{EqHilbertMod}. In particular, we write 
\begin{equation}\label{EqInterStand}
X_B(\Hsp) := X_B(L^2(B),\Hsp).
\end{equation}
We can view $X_B(\Hsp)$ as a self-dual (right) Hilbert $B$-module \cites{Pas73,Rie74}, and more generally $X_B(\Gsp,\Hsp)$ as a W$^*$-TRO \cite{Zet93}, but we will not really need this more abstract viewpoint for now. 

Also the following lemma is elementary. 
\begin{Lem}\label{LemEasyPeasy}
Let $\Gsp$ be a Hilbert space with a normal anti-$*$-representation $\rho_{\Gsp}$ by $B$. Then there exists a net of contractions $x_i \in \cup_{n=1}^{\infty} X_B(L^2(B)^n,\Gsp)$ such that the maps 
\[
B(\Gsp) \rightarrow B(\Gsp),\qquad y \mapsto x_ix_i^*yx_ix_i^*
\]
converge pointwise $\sigma$-weakly to the identity map. 
\end{Lem}
\begin{proof}
It follows immediately from the fact that we can identify 
\[
(\Gsp,\rho_{\Gsp}) \cong ((\pi_{B}\otimes \id)(p)(L^2(B)\otimes l^2(I)),\rho_B\otimes \id)
\]
for some set $I$ and some $p \in B\bar{\otimes} B(l^2(I))$. 
\end{proof}

We will also need a slight extension/variation of \cite{A-DH90}*{Lemma 2.2}. Although the proof almost applies ad verbatim, we prefer to include it for the convenience of the reader.

\begin{Lem}\label{LemRehash}
Let $B$ be a von Neumann algebra, and let $\Hsp,\Gsp$ be Hilbert spaces with a normal anti-$*$-representation $\rho_{\Hsp},\rho_{\Gsp}$ of $B$. Let $F \subseteq X_B(\Gsp,\Hsp)$ be a subset which is closed under right multiplication with $\rho_{\Gsp}(B)'$. Let $S \subseteq B(\Hsp)$ be a subset with $1\in S$, and assume that $\{\Ad_{y^*}: S\rightarrow B(\Gsp)\mid y\in F\}$ is convex. Let $x \in X_B(\Gsp,\Hsp)$ with $\|x\|\leq 1$, and assume that there exists a net of $x_i \in F$ such that for each $s\in S$
\begin{equation}\label{Eqwwww}
x_i^*sx_i \rightarrow x^*sx \quad \textrm{weakly in }B(\Gsp).
\end{equation}
Then there exists a net of \emph{contractive} $x_i \in F$ such that \eqref{Eqwwww} holds for the $\sigma$-weak topology.
\end{Lem} 

\begin{proof}
As in \cite{A-DH90}*{Lemma 2.2}, we begin with proving the lemma under the extra condition that $x^*x = 1$. 

Assume that a net of $x_i \in F$ is given such that \eqref{Eqwwww} holds for each $s\in S$. By the assumption that $\{(\Ad_{y^*})_{\mid S}\mid  y \in F\}$ is convex, we may choose this net so that \eqref{Eqwwww} in fact holds with respect to the \emph{strong} topology. 

Put  
\[
y_i = 2x_i (1+x_i^*x_i)^{-1} \in F. 
\]
Following further the computation as in \cite{A-DH90}*{Lemma 2.2}, we find 
\[
(x_i-y_i)^* (x_i-y_i) \leq (x_i^*x_i-1)^2,
\]
and for each vector state $\omega = \omega_{\xi}$ on $B(\Gsp)$ and $s\in S$ we have
\begin{equation}\label{EqEstADnew}
|\omega(x_i^*sx_i -y_i^*sy_i)|\leq \|s\| \omega((x_i^*x_i-1)^2)^{1/2}(\omega(x_i^*x_i)^{1/2} + \omega(y_i^*y_i)^{1/2}). 
\end{equation}
Since $\|y_i\|\leq 1$ and since $\Ad_{x_i^*}(1)=x_i^*x_i$ converges strongly to $x^*x = 1$ by assumption, we find that \eqref{EqEstADnew} goes to zero in the limit. It follows that $y_i$ is a contractive net in $F$ such that for each $s\in S$ the net $\Ad_{y_i^*}(s)$ converges weakly (and hence $\sigma$-weakly) to $x^*sx$, proving the lemma in this special case.

Let us now consider the case where $x\in X_B(\Gsp,\Hsp)$ is an arbitrary contraction. Let $e \in \rho_{\Gsp}(B)' \subseteq B(\Gsp)$ be the support projection for $x^*x$. Then we may in fact view $\Ad_{x^*}$ as a map $B(\Hsp) \rightarrow B(e\Gsp)$. By multiplying a net $x_i$ as in \eqref{Eqwwww} with $e$ on the right, we see that without loss of generality we may in fact assume that $e = 1$. 

Let then $x = u|x|$ be the polar decomposition of $x$, where $u\in X_B(\Gsp,\Hsp)$ is isometric by the above assumption. Then $u$ lies in the strong closure of the unit ball of $x\rho_{\Gsp}(B)'$. It follows from this last property, together with the fact that $F$ is closed under right multiplication with $\rho_{\Gsp}(B)'$, that $\Ad_{u^*}$ can be approximated pointwise weakly on $S$ with maps $\Ad_{w_i^*}$ for $w_i \in F$. By the first part of the proof, we may assume the $w_i$ to be contractive. But then the $y_i = w_i|x|$ will also be contractive, with $\Ad_{y_i^*}$ converging to $\Ad_{x^
*}$ pointwise weakly and hence $\sigma$-weakly on $S$. 
\end{proof}

Assume now again that $(A,\alpha_A)$ and $(B,\alpha_B)$ are $\G$-W$^*$-algebras.

\begin{Lem}\label{LemApproxSep}
Let $\Hsp,\Gsp$ be two $\G$-equivariant $A$-$B$-correspondences. Assume that $\Gsp$ is weakly contained in $\Hsp$. Then there exists a contractive net $y_{i} \in X_B(\Gsp,\Hsp^{\infty})$ such that for all $c\in C^{\G}(A,B)$
\begin{equation}\label{EqSliceConvLemPrel1}
y_i^*\theta_{\Hsp^{\infty}}(c)y_i \rightarrow \theta_{\Gsp}(c) \quad \sigma\textrm{-weakly in }B(\Gsp).
\end{equation}
\end{Lem} 
\begin{proof}

We first prove the following: if, for some $n$, we are given a contractive map $T = (T_1,\ldots,T_n)\in X_B(L^2(B)^n,\Gsp)$, there exists a $1$-bounded net $y_{i} \in X_B(L^2(B)^n,\Hsp^{\infty})$ such that for all $c\in C^{\G}(A,B)$
\begin{equation}\label{EqSliceConvLemPrelWeak}
y_i^*\theta_{\Hsp^{\infty}}(c)y_i \rightarrow T^*\theta_{\Gsp}(c)T \quad \sigma\textrm{-weakly in }B(L^2(B)^n).
\end{equation}

Indeed, pick a finite set $F_{L^2(B)} \subseteq L^2(B)$ and a finite subset $F_C \subseteq C^{\G}(A,B) \cup \{1\} \subseteq QM(C^{\G}(A,B))$. Take $\varepsilon>0$. We first look for elements $x_1,\ldots,x_n\in X_B(\Hsp^{\infty})$ such that 
\begin{equation}\label{EqIneqNeedPrelLemHere}
|\langle T_i\xi,\theta_{\Gsp}(z)T_j\eta\rangle - \langle x_i\xi,\theta_{\Hsp^{\infty}}(z)x_j\eta\rangle|<\varepsilon,\qquad \forall z \in F_C,\forall \xi,\eta\in F_{L^2(B)},\forall i,j:1\rightarrow n.
\end{equation}
By Lemma \ref{LemPrelResvN}, we can choose $\zeta \in L^2(B)$ and an nsf weight $\varphi$ on $B$ together with a finite set $F_B \subseteq B$ such that each $b^{\opp}\in \mathscr{N}_{\varphi^{\opp}}$ for $b\in F_B$ and for each $\xi \in F_{L^2(B)}$ there exists $b\in F_B$ with 
\[
\xi = \rho_B(b)\zeta = \Lambda_{\varphi^{\opp}}(b^{\opp}). 
\]
Since $(\Hsp,\theta_{\Hsp})$ lies in an arbitrary $\tau^{\infty}$-neighborhood of $(\Gsp,\theta_{\Gsp})$, and since $\rho_u(B) \subseteq QM(C^{\G}(A,B))$ implies that also $\rho_u(B)(C^{\G}(A,B)\cup \{1\})\rho_u(B)\subseteq QM(C^{\G}(A,B))$,  Lemma \ref{LemFellTop} allows us to pick $\zeta_i'\in \Hsp^{\infty}$ with $\|\zeta_i'\|\leq \|T_i\zeta\|\leq \|\zeta\|$ and, for all $b,c\in F_B,z\in F_C$, 
\begin{multline*}
|\langle T_i\rho_{B}(b)\zeta,\theta_{\Gsp}(z)T_j\rho_{B}(c)\zeta \rangle - 
\langle \rho_{\Hsp^{\infty}}(b)\zeta_i',\theta_{\Hsp^{\infty}}(z)\rho_{\Hsp^{\infty}}(c)\zeta_j' \rangle|
\\ =  
|\langle T_i\zeta,\theta_{\Gsp}(\rho_u(b)^*z\rho_u(c))T_j\zeta \rangle - 
\langle \zeta_i',\theta_{\Hsp^{\infty}}(\rho_u(b)^*z\rho_u(c))\zeta_j' \rangle| < \varepsilon.
\end{multline*}
By norm-perturbing the $\zeta_i'$ if necessary, we may assume that they are left bounded for $\varphi$ (in the sense that $\zeta_i' \in \mathfrak{D}(\Hsp^{\infty},\varphi)$, see \cite{Tak03}*{Section IX.3}). But then with $x_i = L_{\zeta_i'}\in X_B(\Hsp^{\infty})$, we find that $\rho_{\Hsp^{\infty}}(b)\zeta_i' = x_i \Lambda_{\varphi^{\opp}}(b^{\opp}) = x_i \rho_B(b)\zeta$, so that \eqref{EqIneqNeedPrelLemHere} is satisfied.

It follows that we can find a net $x_i \in X_B(L^2(B)^n,\Hsp^{\infty})$ such that $\Ad_{x_i^*}\circ \theta_{\Hsp^{\infty}}$ converges weakly to $\Ad_{T^*}\circ \theta_{\Gsp}$ pointwise on $C^{\G}(A,B)\cup \{1\}$. Now the set
\[
F= X_{B}(L^2(B)^n,\Hsp^{\infty}) \subseteq X_{B}(L^2(B)^n,\Gsp \oplus \Hsp^{\infty})
\]
is stable under right multiplication with $\rho_{L^2(B)^n}(B)'$, and with $\mathrm{diag}_{\infty}: B(\Hsp) \rightarrow B(\Hsp^{\infty})$ the diagonal inclusion, we have that 
\[
\{\Ad_{x^*}: B(\Gsp)\oplus \mathrm{diag}_{\infty}( B(\Hsp))\rightarrow B(L^2(B)^n)\mid x\in F\}
\]
is convex. By viewing also $T \in X_{B}(L^2(B)^n,\Gsp) \subseteq X_B(L^2(B)^n,\Gsp\oplus \Hsp^{\infty})$ and 
applying Lemma \ref{LemRehash} to $F \subseteq X_B(L^2(B)^n,\Gsp\oplus \Hsp^{\infty})$ and $S = (\theta_{\Gsp}\oplus \theta_{\Hsp^{\infty}})(C^{\G}(A,B)\cup\{1\})$, we can conclude that we can also find a net $y_i \in X_B(L^2(B)^n,\Hsp^{\infty})$ with $\|y_i\|\leq 1$ such that $\Ad_{y_i^*}\circ\theta_{\Hsp^{\infty}}$ converges $\sigma$-weakly to $\Ad_{T^*}\circ \theta_{\Gsp}$ pointwise on $C^{\G}(A,B)\cup\{1\}$. This proves in particular \eqref{EqSliceConvLemPrelWeak}.

We can now prove the lemma. Pick $T \in X_B(L^2(B)^n,\Gsp)$ contractive for some $n$, and choose a contractive net $y_i \in X_B(L^2(B)^n,\Hsp^{\infty})$ such that $\Ad_{y_i^*}\circ \theta_{\Hsp^{\infty}}$ converges weakly to $\Ad_{T^*}\circ \theta_{\Gsp}$ pointwise on $C^{\G}(A,B)$. Then $y_iT^* \in X_B(\Gsp,\Hsp^{\infty})$ is still contractive, and for all $c\in C^{\G}(A,B)$ we have
\[
(y_iT^*)^*\theta_{\Hsp^{\infty}}(c)(y_iT^*) \rightarrow TT^*\theta_{\Gsp}(c)TT^* \quad \sigma\textrm{-weakly in }B(\Gsp).
\]
By Lemma \ref{LemEasyPeasy}, the pointwise $\sigma$-weak closure of the maps $B(\Gsp) \ni x \mapsto TT^*xTT^* \in B(\Gsp)$ contains the identity map, proving \eqref{EqSliceConvLemPrel1}.
\end{proof}

Given $\G$-equivariant $A$-$B$-correspondences $\Hsp,\Gsp$, we obtain a coaction 
\begin{equation}\label{EqCoactMod}
\alpha_{X_B(\Gsp,\Hsp)}: X_B(\Gsp,\Hsp)\rightarrow X_B(\Gsp,\Hsp)\bar{\otimes} M,\quad x\mapsto U_{\Hsp}(x\otimes 1)U_{\Gsp}^*.
\end{equation}

Here, the tensor product $X_B(\mathcal{G}, \mathcal{H})\bar{\otimes} M$ can be interpreted in multiple ways. We can view it as the Fubini tensor product of operator spaces, or alternatively, we can view $X_B(\mathcal{G}, \mathcal{H})\odot M$ as a corner of 
$$X_B(\mathcal{G}\oplus\mathcal{H}, \mathcal{G}\oplus \mathcal{H})\bar{\otimes} M$$
and then define $X_B(\mathcal{G}, \mathcal{H})\bar{\otimes} M$ to be the $\sigma$-weak closure of $X_B(\mathcal{G}, \mathcal{H})\odot M$ inside $X_B(\mathcal{G}\oplus\mathcal{H}, \mathcal{G}\oplus \mathcal{H})\bar{\otimes} M$. Both interpretations agree and we can make sense of the notion of a $\G$-coaction on $X_B(\mathcal{G}, \mathcal{H})$, see e.g. \cite{PS17}*{Section 1 \& 3} for more details and information.

\begin{Theorem}\label{TheoMain1}
Let $A,B$ be $\G$-W$^*$-algebras, and let $\Hsp,\Gsp$ be $\G$-$A$-$B$-correspondences. Then $\Gsp$ is equivariantly weakly contained in $\Hsp$ if and only if there exists a $1$-bounded net $y_{i} \in X_B(\Gsp,\Hsp^{\infty})$ such that for all $a\in A$
\begin{equation}\label{EqSliceConv}
(y_i^*\pi_{\Hsp^{\infty}}(a)\otimes 1)\alpha_{X_B(\Gsp,\Hsp^{\infty})}(y_i) \rightarrow \pi_{\Gsp}(a)\otimes 1 \quad \sigma\textrm{-weakly in }B(\Gsp\otimes L^2(\G)).
\end{equation}
\end{Theorem} 
\begin{proof}
Assume first that $\Gsp$ is weakly contained in $\Hsp$. By Lemma \ref{LemApproxSep}, we can find a contractive net $y_{i} \in X_B(\Gsp,\Hsp^{\infty})$ such that for all $c\in C^{\G}(A,B)$ we have
\begin{equation}\label{EqSliceConvLemPrel}
y_i^*\theta_{\Hsp^{\infty}}(c)y_i \rightarrow \theta_{\Gsp}(c) \quad \sigma\textrm{-weakly in }B(\Gsp).
\end{equation}

By passing to a subnet, we may assume that
\begin{equation}\label{EqTyi}
\Ad_{y_i^*}: B(\Hsp^{\infty}) \rightarrow B(\Gsp),\qquad x\mapsto y_i^*xy_i.
\end{equation}
converges pointwise $\sigma$-weakly to a ccp map 
\[
\widetilde{E}: B(\Hsp^{\infty}) \rightarrow B(\Gsp).
\]
Put 
\begin{equation}\label{EqLimitMapE}
E = \widetilde{E}\circ \mathrm{diag}_{\infty}: B(\Hsp) \rightarrow B(\Gsp),
\end{equation}
where $\mathrm{diag}_{\infty}: B(\Hsp) \rightarrow B(\Hsp^{\infty})$ is the diagonal inclusion. Then clearly $E$ is a ccp map with 
\[
E(\theta_{\Hsp}(z)) = \theta_{\Gsp}(z),\qquad \forall z\in C^{\G}(A,B).
\]
In particular, $\theta_{\Hsp}(C^{\G}(A,B))$ belongs to the multiplicative domain of $E$. But then also 
\[
\theta_{\Gsp}(x)E(\theta_{\Hsp}(z))\theta_{\Gsp}(y) = E(\theta_{\Hsp}(xzy)) = \theta_{\Gsp}(x)\theta_{\Gsp}(z)\theta_{\Gsp}(y),\qquad x,y\in C^{\G}(A,B),z\in QM(C^{\G}(A,B)). 
\]
It follows that
\begin{equation}\label{EqQuasiMultEprel}
E(\theta_{\Hsp}(z)) = \theta_{\Gsp}(z),\qquad \forall z\in QM(C^{\G}(A,B)). 
\end{equation}
Since $\pi_u(A) C_{U_u} \subseteq QM(C^{\G}(A,B))$, we find 
\begin{equation}\label{EqAuxConv}
y_i^*\pi_{\Hsp^{\infty}}(a)U_{\Hsp^{\infty}}(\omega)y_i \overset{\sigma\textrm{-weakly}}{\underset{i\rightarrow \infty}{\rightarrow}} \pi_{\Gsp}(a)U_{\Gsp}(\omega),\qquad\forall a\in A,\omega \in M_*.
\end{equation} 
Hence for all $\xi,\eta\in \Gsp$ and all $f,g\in L^2(\G)$ and $a\in A$, we have 
\[
\langle y_i\xi\otimes f, (\pi_{\Hsp^{\infty}}(a)\otimes 1)U_{\Hsp^{\infty}}(y_i\eta\otimes g)\rangle \rightarrow \langle \xi\otimes f,(\pi_{\Gsp}(a)\otimes 1)U_{\Gsp}(\eta\otimes g)\rangle.
\]
By boundedness of the $y_i$, this means that for all $a\in A$ we have
\begin{equation}\label{EqWeakConvaU}(y_i^*\pi_{\Hsp^{\infty}}(a) \otimes 1)U_{\Hsp^{\infty}}(y_i\otimes 1)
\underset{i\rightarrow \infty}{\overset{\sigma\textrm{-weakly}}{\longrightarrow}}(\pi_{\Gsp}(a)\otimes 1)U_{\Gsp}.
\end{equation}
Multiplying to the right with $U_{\Gsp}^*$, we find \eqref{EqSliceConv}. 

Let us now go back. Assume that we have chosen contractive $y_k\in X_B(\Gsp,\Hsp^{\infty})$ satisfying \eqref{EqSliceConv}. Choose finitely many $\xi_i \in \Gsp$ as well as a finite set $F_C \subseteq C^{\G}(A,B)$ and $\varepsilon>0$. By Theorem \ref{TheoEquivTop}, it is enough find $\xi_i'\in \Hsp^{\infty}$ such that $\|\xi_i'\|\leq \|\xi_i\|$ and 
\[
|\langle \xi_i,\theta_{\Gsp}(z)\xi_j\rangle - \langle \xi_i',\theta_{\Hsp^{\infty}}(z)\xi_j'\rangle|<\varepsilon,\qquad \forall z\in F_C,i,j\in \{1,\ldots,n\}.  
\] 
Now by assumption, we have for all $\omega \in M_*$ and $a\in A$ that
\[
y_k^*\pi_{\Hsp^{\infty}}(a)U_{\Hsp^{\infty}}(\omega) y_k \underset{k\rightarrow \infty}{\overset{\sigma-\textrm{weakly}}{\longrightarrow}} \pi_{\Gsp}(a)U_{\Gsp}(\omega). 
\]
Then clearly also for all $\omega \in M_*, a\in A$ and $b\in B$ we have 
\[
y_k^*\pi_{\Hsp^{\infty}}(a)U_{\Hsp^{\infty}}(\omega)\rho_{\Hsp^{\infty}}(b) y_k = y_k^*\pi_{\Hsp^{\infty}}(a)U_{\Hsp^{\infty}}(\omega)y_k \rho_{\Gsp}(b) \underset{k\rightarrow \infty}{\overset{\sigma-\textrm{weakly}}{\longrightarrow}} \pi_{\Gsp}(a)U_{\Gsp}(\omega)\rho_{\Gsp}(b). 
\]
Since $[\pi_u(A_{\slice})C_{U_u}\rho_u(B_{\slice})]= C^{\G}(A,B)$, the uniform boundedness of the $y_k$ ensures that for all $z\in C^{\G}(A,B)$ we have
\[
y_k^*\theta_{\Hsp^{\infty}}(z) y_k \underset{k\rightarrow \infty}{\overset{\sigma-\textrm{weakly}}{\longrightarrow}} \theta_{\Gsp}(z).
\]
Hence for our finite collection of $\xi_i$ and $z\in F_C$, we have 
\[
|\langle \xi_i,\theta_{\Gsp}(z)\xi_j\rangle - \langle y_k\xi_i,\theta_{\Hsp^{\infty}}(z)y_k\xi_j\rangle| \underset{k\rightarrow \infty}{\longrightarrow} 0,
\]
so we can take $\xi_i' = y_k\xi_i$ for $k$ large enough.
\end{proof}

Let us recuperate part of the proof of Theorem \ref{TheoMain1} as a separate corollary.

\begin{Cor}\label{CorWeakContCCP}
If $\Gsp$ is equivariantly weakly contained in $\Hsp$,  there exists a ucp map 
$E: B(\Hsp) \rightarrow B(\Gsp)$
such that 
\begin{equation}\label{EqQuasiMultE}
E(\pi_{\Hsp}(a)) = \pi_{\Gsp}(a),\qquad E(\rho_{\Hsp}(b)) = \rho_{\Gsp}(b), \qquad E(U_{\Hsp}(\omega)) = U_{\Gsp}(\omega)\qquad \forall a\in A, b\in B,\omega \in M_*.
\end{equation}
\end{Cor}

\begin{proof} 
Since the map $E$ of \eqref{EqLimitMapE} satisfies \eqref{EqQuasiMultEprel}, and since $\pi_u(A),\rho_u(A)\subseteq QM(C^{\G}(A,B))$, we immediately find that \eqref{EqQuasiMultE} is satisfied.
\end{proof}

\section{Amenability for equivariant correspondences and \texorpdfstring{$\G$}{G}-actions}

Let $A,B,C$ be von Neumann algebras. Given a ucp map $E: A\to B$ (not necessarily normal!), we use the notation
\begin{equation}\label{EqFubiniCPMap}
E\otimes \id: A\bar{\otimes} C \rightarrow B\bar{\otimes} C
\end{equation}
to denote the unique ucp map such that 
\[
E((\id\otimes \omega)z) = (\id\otimes \omega)((E\otimes \id)z),\qquad \forall z\in A\bar{\otimes} C, \omega \in C_*,
\]
using for example the identification 
\[
B\bar{\otimes} C \cong CB(C_*,B),
\]
where the right hand side means the set of all cb maps \cite{ER00}*{Proposition 8.1.2}. This allows us to give meaning to \eqref{EqEquivCEE}, which we then refer to as $\G$-equivariance of the corresponding ucp map.

\begin{Def}\label{DefAmenable}
Let $\pi: A \rightarrow B$ be a $\G$-equivariant normal unital $*$-homomorphism.

We say that $\pi:A \rightarrow B$ is \emph{strongly equivariantly amenable} if $E_A = {}_{A}L^2(A)_A$ is equivariantly weakly contained in ${}_{\pi(A)}L^2(B)_{\pi(A)}$, which we view as $\G$-$A$-$A$-correspondence. More precisely, ${}_{\pi(A)}L^2(B)_{\pi(A)}$ is the $\G$-$A$-$A$-correspondence $(L^2(B), \pi_B\circ \pi: A \to B(L^2(B)), \rho_B\circ \pi: A \to B(L^2(B)), U_\beta)$.

We say that $\pi: A\rightarrow B$ is \emph{equivariantly amenable} if there exists a $\G$-equivariant ucp map
\[
P: B \rightarrow A
\]
which is left inverse to $\pi$. 
\end{Def}
We may view equivariant amenability as a relative version of $\G$-equivariant injectivity. Note also that if $\pi$ is (strongly) equivariantly amenable, then $\pi$ is automatically faithful, so we are actually dealing with a $\G$-equivariant inclusion $A \subseteq B$.

\begin{Rem}\label{RemTermin}
The notion of equivariant amenability was already introduced in \cite{Moa18}*{Definition 4.4} in the setting of discrete quantum group actions. Note that we follow here the terminology of \cite{BMO20}, which is dual to the one followed for example in \cite{A-D90}.
\end{Rem}

\begin{Theorem}\label{TheoInclStrong}
Let $\pi: A \rightarrow B$ be a $\G$-equivariant unital normal $*$-homomorphism between $\G$-W$^*$-algebras. If $\pi$ is strongly $\G$-amenable, then $\pi: A \rightarrow B$ is $\G$-amenable. 
\end{Theorem}
\begin{proof}
By assumption, the trivial equivariant $A$-$A$-correspondence $L^2(A)$ is equivariantly weakly contained in the equivariant $A$-$A$-correspondence ${}_{\pi(A)}L^2(B)_{\pi(A)}$. Consider the ucp map
\[
E: B(L^2(B)) \rightarrow B(L^2(A))
\] 
as in Corollary \ref{CorWeakContCCP}. Then $E(\pi_B(B))$ will lie in $\pi_A(A)$, as these elements commute pointwise with $\rho_A(A)$ by a multiplicative domain argument. We can hence uniquely define a ccp map $P:B\rightarrow A$ by 
\[
\pi_A(P(z)) := E(\pi_B(z)),\qquad z\in B.
\]
It is then clear that $P(\pi(a)) = a$ for all $a\in A$. Moreover, we have by \eqref{EqQuasiMultE} that $U_{\beta}$, the unitary $\G$-representation for the $A$-$A$-correspondence ${}_AL^2(B)_A$, lies in the multiplicative domain of $E\otimes \id_{M}$, hence for $z\in B$ and $\omega \in M_*$ we get 
\begin{eqnarray*} 
\pi_A(P((\id\otimes \omega)\beta(z))) &=& (\id\otimes \omega)(E\otimes \id)(U_{\beta}(\pi_B(z)\otimes 1)U_{\beta}^*) \\
&=&  (\id\otimes \omega)(U_{\alpha} (E(\pi_B(z))\otimes 1)U_{\alpha}^*)\\
&=& \pi_A((\id\otimes \omega)\alpha(P(z))),
\end{eqnarray*}
so that $P$ is $\G$-equivariant.
\end{proof}

We can use the definition of (strong) amenability for equivariant normal unital $*$-homomorphisms to define a notion of left (strong) amenability for equivariant correspondences. Recall the constructions made in Example \ref{ExaBasCons}.

\begin{Def}\label{DefEquivariantlyStrong}
We say that a $\G$-equivariant $A$-$B$-correspondence $(\Hsp,U,\pi,\rho)$ is \emph{(strongly) left $\G$-amenable} if $\pi: A \rightarrow C= \rho(B)'$ is (strongly) $\G$-amenable.
\end{Def}
\begin{Rem}
In \cite{A-D95}*{Definition 2.1}, strong left amenability (called left amenability loc.\ cit.) is rather introduced in terms of the Connes fusion product. We use for now the shortcut provided by the remarks following \cite{A-D95}*{Definition 2.1}, and will return to this point in  Remark \ref{RemSauv}.
\end{Rem}

\begin{Rem}\label{RemAmpliIsMute}
If $\Gsp$ is a multiplicity Hilbert space, then it is easily seen that $(\Hsp,U,\pi,\rho)$ is (strongly) left $\G$-amenable if and only if $(\Gsp\otimes \Hsp,\id_{\Gsp}\otimes U,\id_{\Gsp}\otimes\pi,\id_{\Gsp}\otimes \rho)$ is (strongly) left $\G$-amenable. Indeed, if $\rho(B)' = C$, then $(\id_{\Gsp}\otimes \rho(B))' = B(\Gsp)\bar{\otimes} C$, and by Remark \ref{RemMulti} we know that 
\[
E_A \preccurlyeq {}_{\pi(A)}L^2(C)_{\pi(A)} \quad \iff\quad E_A \preccurlyeq L^2(B(\Gsp))\otimes {}_{\pi(A)}L^2(C)_{\pi(A)} = {}_{(\id_{\G}\otimes \pi)(A)}L^2(B(\Gsp)\bar{\otimes} C)_{(\id_{\G}\otimes \pi)(A)},
\]
showing the claim for strong left $\G$-amenability. For left $\G$-amenability, we notice that if $P: C \rightarrow A$ is a $\G$-equivariant ucp map which splits $\pi$, and $\omega$ an arbitrary normal state on $B(\Gsp)$, then 
\begin{equation}\label{EqSplittAmpl}
\widetilde{P}: B(\Gsp)\bar{\otimes} C \rightarrow A,\quad z \mapsto P((\omega\otimes \id)z)
\end{equation}
is a $\G$-equivariant ucp map which splits $\id_{\Gsp}\otimes \pi$. Conversely, given \eqref{EqSplittAmpl}, the map 
\[
P: C \rightarrow A,\quad c\mapsto \widetilde{P}(\id_{\Gsp}\otimes c)
\]
will be a $\G$-equivariant ucp map which splits $\pi$. This shows that left $\G$-equivariant (non-)amenability is unperturbed by amplification. 
\end{Rem} 

\begin{Theorem}\label{TheoAmenacorr}
If $\Hsp$ is strongly left equivariantly amenable, then it is left equivariantly amenable.
\end{Theorem}
\begin{proof}
This is a direct consequence of Theorem \ref{TheoInclStrong}.
\end{proof}

In \cite{BMO20}*{Corollary A.2}, it was proven that in the non-equivariant setting, the above can be strengthened: strong left amenability of $\Hsp$ is in fact \emph{equivalent} with left amenability. In the equivariant setting, this is currently yet unknown, as the following special situation demonstrates. 

\begin{Def}\label{EqDefAmAct}
Let $(A,\alpha)$ be a $\G$-dynamical von Neumann algebra. We say that the action of $\G$ on $A$ is \emph{(strongly) $\G$-amenable} if the inclusion $\alpha: A \rightarrow A\bar{\otimes} M$ is (strongly) $\G$-amenable. 
\end{Def}
Concretely, recalling the notation introduced in Example \ref{Exa3}, strong $\G$-amenability of $\alpha$ means that $E_A^{\G}$ is $\G$-equivariantly weakly contained in $S_A^{\G}$ and $\G$-amenability of $\alpha$ means that there exists a ucp map 
\begin{equation}\label{EqDefPDisc}
P: A\bar{\otimes} M \rightarrow A
\end{equation}
such that $P(\alpha(a)) = a$ for all $a\in A$ and 
\begin{equation}\label{EqEquiDisQu}
P((\id\otimes (\id\otimes \omega)\Delta)z) = (\id\otimes \omega)\alpha(P(z)),\qquad \forall z\in A\bar{\otimes}M,\omega \in M_*.
\end{equation}
It is not clear to us if the definition of $\G$-amenability has explicitly appeared in print before in the general setting of locally compact quantum groups, but it is the obvious one, see e.g.\ \cite{Moa18}*{Definition 4.1} for the setting with $\G$ discrete, or \cites{Voi77, Rua96, BT03, DQV02} when $A= \C$ . On the other hand, for $\G=G$ an ordinary locally compact group, the above is equivalent with the existence of a ucp map 
\begin{equation}\label{EqEquivProj}
P: A\bar{\otimes} L^{\infty}(G) \rightarrow A
\end{equation}
such that $P(a\otimes 1)=a$ for all $a\in A$ and which is equivariant with respect to the diagonal action of $G$ on $A\bar{\otimes} L^{\infty}(G)$, using the action $\rho$ of $G$ on $L^{\infty}(G)$ by right shifts. This is the notion of amenability of the action $G \curvearrowright A$ as it appears in \cite{A-D79}*{Definition 3.4}.\footnote{In \cite{A-D79} one uses left shifts, but this is of course inconsequential.} 

Specializing Theorem \ref{TheoInclStrong}, we obtain:
\begin{Theorem}\label{TheoConsStrong}
Assume $A$ is a von Neumann algebra with a strongly $\G$-amenable $\G$-action. Then the action of $\G$ on $A$ is $\G$-amenable. 
\end{Theorem}
\begin{proof}
This is just a special case of Theorem \ref{TheoInclStrong}.
\end{proof}

If we take $A = \C$, then the above simply states that if the trivial $\G$-representation is weakly contained in the regular $\G$-representation (= co-amenability of $\hat{\G}$ \cite{BT03}), there exists a $\G$-invariant state $M \rightarrow \C$ (= amenability of $\G$). The converse is a well-known open problem, making it hence difficult to ascertain if the converse of Theorem \ref{TheoAmenacorr} or even Theorem \ref{TheoConsStrong} will hold. Let us now show that the converse of the result in Theorem \ref{TheoConsStrong} holds for locally compact groups.

\begin{Theorem}\label{TheoClassCase}
If $G$ is a locally compact group, and $A$ a $G$-W$^*$-algebra, then the action of $G$ on $A$ is strongly $G$-amenable if and only if it is $G$-amenable.
\end{Theorem}
\begin{proof}
Note first that an action of $G$ on $A$ is the same thing as a coaction
\[
\alpha: A\rightarrow A\bar{\otimes}L^{\infty}(G),\qquad a \mapsto \alpha(a),
\]
where we view $\alpha(a) \in B(L^2(A)\otimes L^2(G)) = B(L^2(G,L^2(A)))$ as multiplication by the function $g \mapsto \pi_A(\alpha_g(a))$. Note that $M= L^{\infty}(G)$ is endowed with the structure of a locally compact quantum group by putting $\Delta(f)(g,h) = f(gh)$ for $f\in L^{\infty}(G)$ (this definition makes sense locally almost everywhere). 

By Theorem \ref{TheoConsStrong}, it is enough to prove that $G$-amenability of the action $G \curvearrowright A$ implies strong amenability, i.e.\ the weak containment of $E_A$ in the  $G$-semi-coarse correspondence $S^G_A$. 

Let $V_g$ for $g\in G$ denote the right shift of $g$ on $L^2(G)$ (so we work with respect to $dg$ the right invariant Haar measure). Note that this is simply the $G$-representation corresponding to the right regular representation $V$ of $L^{\infty}(G)$ as a locally compact quantum group. Similarly, we write $g \mapsto U_g$ for the unitary $G$-representation associated to $U_{\alpha}$. Let $Z(A)$ be the center of $A$, and let $Z(A)_{\slice}$ be as in Definition \ref{DefSliceCat}. In this case, $Z(A)_{\slice}$ is actually the C$^*$-algebra of elements which are norm-continuous for the $G$-action \cite{Ped79}*{Lemma 7.5.1}. Let then $C_c(G,Z(A)_{\slice})$ be the $*$-algebra of norm-continuous maps $G \rightarrow Z(A)_{\slice}$ with compact support. By \cite{BC22}*{Theorem 1.1.(1) $\Rightarrow$ (2)}, we can choose a net
\[
x_i \in C_c(G,Z(A)_{\slice}) \subseteq Z(A)\bar{\otimes} L^2(G) \subseteq B(L^2(A),L^2(A)\otimes L^2(G))
\]
which is contractive for the $B(L^2(A),L^2(A)\otimes L^2(G))$-norm and such that 
\begin{equation}
x_i^*(\alpha_g\otimes V_g)(x_i) \overset{\sigma\textrm{-weakly}}{\underset{i\rightarrow\infty}{\longrightarrow}} 1
\end{equation}
uniformly on compacts. If we define $z_i(g) = x_i^*(\alpha_g\otimes V_g)(x_i)$, it is elementary to verify that $z_i \in C_c(G,Z(A)_{\slice})$, and that hence 
\[
\int_G f(g) z_i(g) \rd g \overset{\sigma\textrm{-weakly}}{\underset{i\rightarrow\infty}{\longrightarrow}} \left(\int_G f(g)\rd g\right)1
\]
for each $f\in C_c(G)$. A density argument now allows to conclude that 
\begin{equation}\label{EqGlobalis}
(x_i^*\otimes 1)V_{23}(\alpha\otimes \id)(x_i)_{132} \overset{\sigma\textrm{-weakly}}{\underset{i\rightarrow\infty}{\longrightarrow}} 1.
\end{equation}

Put $y_i = U_{\alpha}x_i$. Then $y_i \in  X:=X_A(L^2(A)\otimes L^2(G))$ (using the right $A$-representation as in Example \ref{ExaSemiCoarsG}). Moreover, the identity $U_{\alpha,12}^*V_{23}U_{\alpha,12} = U_{\alpha,12} V_{23}$ transforms \eqref{EqGlobalis} into 
\begin{equation}\label{EqCoactConv}
(y_i^*\otimes 1)\alpha_{X}(y_i) \overset{\sigma\textrm{-weakly}}{\underset{i\rightarrow\infty}{\longrightarrow}} 1.
\end{equation}
But since the first leg of $x_i$ lies in $Z(A)$, multiplying \eqref{EqCoactConv} to the left with $\pi_A(a)\otimes 1$ for $a\in A$ entails that 
\begin{equation}\label{EqCoactConv2}
(y_i^*(\pi_A\otimes \pi_M)\alpha(a)\otimes 1) \alpha_{X}(y_i) \overset{\sigma\textrm{-weakly}}{\underset{i\rightarrow\infty}{\longrightarrow}} (\pi_A\otimes \pi_M)\alpha(a)\otimes 1,\qquad \forall a\in A.
\end{equation}
It now follows from Theorem \ref{TheoMain1} that $E_A$ is weakly contained in $S_A^G$. 
\end{proof}

To end this section, we note that Theorem \ref{TheoMain1} seems to allow for a stronger and more symmetric condition of weak containment. Namely, let $A,B$ be $\G$-W$^*$-algebras, and let $\Hsp,\Gsp$ be $\G$-$A$-$B$-correspondences. Denote 
\[
{}_AX_B(\Gsp,\Hsp) = \{x\in X_B(\Gsp,\Hsp)\mid x\pi_{\Gsp}(a) = \pi_{\Hsp}(a)x \textrm{ for all }a\in A\}.
\]
Let us then say that $\Hsp$ \emph{symmetrically weakly} contains $\Gsp$ if there exists a contractive net $y_{i} \in  {}_AX_B(\Gsp,\Hsp^{\infty})$ such that
\begin{equation}\label{EqSliceConvSymm}
(y_i^*\otimes 1)\alpha_X(y_i) \rightarrow 1\quad \sigma\textrm{-weakly in }B(\Gsp\otimes L^2(\G)).
\end{equation}
Clearly, if $\Hsp$ symmetrically weakly contains $\Gsp$, then it weakly contains $\Gsp$,  but in the characterisation of weak containment in Theorem \ref{TheoMain1}, the $x_i$ only satisfy $x_i\pi_{\Gsp}(a) = \pi_{\Hsp}(a)x_i$ \emph{approximately}. 

If we can find $y_i  \in {}_AX_B(\Gsp,\Hsp^{\infty})$ satisfying \eqref{EqSliceConvSymm}, then it is easily seen that ${}_A\Hsp_B^{\infty}$ must contain a (not necessarily equivariant) isometric copy of ${}_A\Gsp_B$ as an $A$-$B$-correspondence. It is hence clear that symmetric weak containment is in general much stronger than ordinary weak containment, and already puts a very strong condition on the underlying correspondences. However, Theorem \ref{TheoClassCase} shows that these notions \emph{do coincide} when considering the weak inclusion of $E_A$ inside the $\G$-semi-coarse correspondence $S_A^{G}$ for $G$ a locally compact group. We do not know if, for a general locally compact quantum group $\G$, the equivariant weak containment $E_A \preccurlyeq S_A^{\G}$ will imply the symmetric equivariant weak containment. Note that the latter, when framed appropriately in the C$^*$-context, coincides with how amenability is defined for actions of discrete quantum groups on unital C$^*$-algebras in \cite{VV07}*{Definition 4.1}.

\section{Operations on equivariant correspondences and the Takesaki-Takai duality}\label{SecOpEquiCorr}

In this section, we consider particular constructions which can be performed on equivariant correspondences.

\subsection{Opposites of \texorpdfstring{$\G$}{G}-equivariant correspondences}

Fix a locally compact quantum group $\G = (M,\Delta)$ and $\G$-W$^*$-algebras $A,B$.

\begin{Def}\label{DefOppositeCorr}
Let $\Hsp$ be a $\G$-equivariant $A$-$B$-correspondence. Let $\overline{\Hsp}$ be the complex conjugate Hilbert space, and let $C_{\Hsp}: \Hsp \rightarrow \overline{\Hsp}$ be the canonical anti-unitary conjugation map $\xi\mapsto \overline{\xi}$. 
The \emph{opposite} $B$-$A$-correspondence is obtained by endowing  $\overline{\Hsp}$ with the bimodule structure 
\[
\overline{\pi}(b) = C_{\Hsp}\rho(b^*)C_{\overline{\Hsp}},\qquad \overline{\rho}(a) = C_{\Hsp}\pi(a^*)C_{\overline{\Hsp}} 
\]
and the unitary $\G$-representation 
\[
\overline{U}= (C_{\Hsp}\otimes \check{J})U^*(C_{\overline{\Hsp}}\otimes \check{J}).
\]
\end{Def}

We can easily see that $\overline{\Hsp}$ is again a $\G$-equivariant correspondence, using the defining relations of the $\G$-equivariant correspondence $\Hsp$ and the fact that $\check{J}$ implements the unitary antipode.

Let us consider this operation for the examples appearing in Section \ref{SecIntroAndExa}. Of course, we may in Definition \ref{DefOppositeCorr} use \emph{any} anti-unitary operator $C_{\Hsp}$ from $\Hsp$ to (another or the same) Hilbert space, as this will clearly give an equivariant copy. 

\begin{Exa}\label{ExaOppExasCan}
We have $\overline{E_A} \cong E_A$ by taking $C_{\Hsp}= J_A$ and using \eqref{EqCanImplUni}.
\end{Exa}
\begin{Exa}\label{ExaResRestr}
Let $A\subseteq C$ and $B\subseteq D$ be $\G$-equivariant inclusions of $\G$-W$^*$-algebras, and ${}_C\Hsp_D$ a $\G$-equivariant $C$-$D$-correspondence. Then
\[
{}_B\Res_A(\overline{\Hsp}) = \overline{{}_A\Res_B(\Hsp)}.
\]
\end{Exa}

\begin{Exa}
With $D^{\G}(-)$ the adjoint amplification, we see that 
$D^{\G}(\overline{\Hsp}) \cong \overline{D^{\G}(\Hsp)}$. Indeed, implementing the latter by means of the anti-unitary map 
\[
C_{D^{\G}(\Hsp)} = C_{\Hsp}\otimes \check{J}: \Hsp \otimes L^2(\G) \rightarrow \overline{\Hsp}\otimes L^2(\G),
\]
we find the unitary intertwiner 
\[
(C_{\Hsp}\otimes \check{J})U(C_{\overline{\Hsp}}\otimes \check{J}): D^{\G}(\overline{\Hsp}) \rightarrow \overline{D^{\G}(\Hsp)}, 
\]
using relation \eqref{EqBraidedCommInt} for the intertwining properties of the unitary $\G$-representations. 
\end{Exa}

It is easy to see that the statement in the previous example also holds for the regular amplification.

\begin{Exa}
  We have $\overline{C^\G_{A,B}} \cong C_{B,A}^\G$. Indeed, implementing $C_{A,B}^\G$ using the anti-unitary $J_A\otimes J \otimes J_B$, one verifies that the flip
  $$L^2(A)\otimes L^2(\G)\otimes L^2(B)\cong L^2(B)\otimes L^2(\G)\otimes L^2(A)$$
  defines the desired unitary intertwiner.
\end{Exa}

\subsection{Composition of \texorpdfstring{$\G$}{G}-equivariant correspondences}

Let $A,B$ be von Neumann algebras. Recall that $X_B(\Hsp)$ was introduced in \eqref{EqInterStand} as the self-dual Hilbert $A$-$B$-module of right $B$-intertwiners from $L^2(B)$ to $\Hsp$. If then $\Gsp$ is a $B$-$C$-correspondence, we can compose $\Hsp$ and $\Gsp$ into an $A$-$C$-correspondence by the \emph{Connes fusion product}
\[
{}_A(\Hsp \boxtimes \Gsp)_C := {}_A(X_B(\Hsp) \underset{B}{\otimes} \Gsp)_C,
\]
which is the separation-completion of $X_B(\Hsp) \odot_B \Gsp$ with respect to the semi-inner product 
\[
\langle x\otimes_B \xi,y\otimes_B \eta\rangle = \langle \xi,x^*y\eta\rangle,\qquad x,y\in X_B(\Hsp),\xi,\eta\in \Gsp,
\]
(where we identified $B \cong \pi_B(B)$). The $A$-$C$-bimodule structure is determined by 
\[
\pi_{\boxtimes}(a)(x\otimes_B \xi) = \pi_{\mathcal{H}}(a)x\otimes_B \xi,\qquad \rho_{\boxtimes}(c)(x\otimes_B \xi) = x\otimes_B \rho_{\mathcal{G}}(c)\xi,\qquad a\in A,c\in C, x\in X_B(\Hsp),\xi\in \Gsp.
\]
We write $\Hsp\boxtimes_B \Gsp$ if we want to emphasize that we take the fusion product over $B$.

If $\Ksp$ is a $C$-$D$-correspondence, we have a natural intertwining associativity unitary 
\begin{equation}\label{EqAssInt}
(\Hsp \boxtimes \Gsp)\boxtimes \Ksp \cong \Hsp \boxtimes (\Gsp \boxtimes \Ksp),
\end{equation}
using e.g.\ that $X_C(\Hsp \boxtimes \Gsp)$ is the $\sigma$-weak closure of the image of $X_B(\Hsp)\odot_B X_C(\Gsp)$ in $B(L^2(C),\Hsp \boxtimes \Gsp)$ by 
\[
(x\otimes_B y)\eta = x \otimes_B y\eta,\qquad x\in X_B(\Hsp),y\in X_C(\Gsp),\eta\in L^2(C).
\]
The identity correspondence $E_A = {}_AL^2(A)_A$ acts as the identity for this composition. See e.g.\ \cites{Lan01,Bro03} for a detailed treatment.

We can upgrade the above construction to the $\G$-equivariant setting. Recall that if $(A,\alpha)$ and $(B,\beta)$ are $\G$-W$^*$-algebras and $\Hsp$ a $\G$-equivariant $A$-$B$-correspondence, then $X = X_B(\Hsp)$ comes equipped with the coaction 
\[
\alpha_X = \alpha_{X_B(\Hsp)}: X \rightarrow X \bar{\otimes} M,\qquad x \mapsto U(x\otimes 1)U_{\beta}^*.
\]
\begin{Prop}\label{PropCrossRep}
Assume $(A,\alpha),(B,\beta),(C,\gamma)$ are $\G$-dynamical von Neumann algebras. If $\Hsp$ is a $\G$-equivariant $A$-$B$-correspondence and $\Gsp$ a $\G$-equivariant $B$-$C$-correspondence, we can endow $\Hsp \boxtimes \Gsp$ with a unique unitary $\G$-representation $U_{\boxtimes}$ such that 
\begin{equation}\label{EqUnitCorep}
U_{\boxtimes} ((x \otimes_B \xi)\otimes \eta) := \alpha_X(x) \otimes_{B \bar{\otimes} M} U_{\Gsp}(\xi\otimes \eta), \qquad x\in X_B(\Hsp),\xi\in \Hsp,\eta\in \Gsp.
\end{equation}
Then $\Hsp\boxtimes \Gsp$  becomes a $\G$-equivariant $A$-$C$-correspondence.
\end{Prop}
\begin{proof}
First of all, note that in \eqref{EqUnitCorep} one rather defines $U_{\boxtimes}$ as a map 
\[
(X_B(\Hsp) \odot_B \Gsp)\odot L^2(M) \rightarrow (X_B(\Hsp)\bar{\otimes}M)\underset{B\bar{\otimes}M}{\otimes} (\Gsp\otimes L^2(M)) \cong (X_B(\Hsp)\underset{B}{\otimes} \Gsp) \otimes L^2(M) = (\Hsp \boxtimes \Gsp)\otimes L^2(M).
\]
A direct computation shows that this indeed descends to an isometry 
\[
U_{\boxtimes}: (\Hsp\boxtimes\Gsp)\otimes L^2(M) \rightarrow (\Hsp\boxtimes\Gsp)\otimes L^2(M).
\]
By the bicommutant theorem, it follows straightforwardly that $U_{\boxtimes} \in B(\Hsp\boxtimes \Gsp)\bar{\otimes} M$. 

To check the corepresentation property, various identifications of Hilbert spaces have to be made, which we will not comment on explicitly. The main goal is to prove the identity 
\begin{equation}\label{EqCorepToProve}
U_{\boxtimes,13}U_{\boxtimes,12} = W_{32}^*U_{\boxtimes,12}W_{32},
\end{equation}
which we will show by applying both sides to a vector $(x\otimes_B \xi)\otimes \eta\otimes \zeta$ with $x \in X_B(\Hsp),\xi\in \Gsp$ and $\eta,\eta\in L^2(\G)$. Then on the one hand, the left hand side gives 
\begin{eqnarray*}
U_{\boxtimes,13}U_{\boxtimes,12}((x\otimes_B \xi)\otimes \eta\otimes \zeta)&\cong& U_{\boxtimes,13}((\alpha_X(x)\otimes_{B\bar{\otimes}M} U_{\Gsp}(\xi\otimes \eta))\otimes \zeta)\\
&\cong & (\alpha_X\otimes \id)\alpha_X(x)_{132} \otimes_{B\bar{\otimes}M\bar{\otimes}M}(U_{\Gsp,13}U_{\Gsp,12}(\xi\otimes \eta\otimes \zeta)),
\end{eqnarray*}
while for the right hand side we get
\begin{eqnarray*}
W_{32}^*U_{\boxtimes,12}W_{32}((x\otimes_B \xi)\otimes \eta\otimes \zeta)& = &W_{32}^*U_{\boxtimes,12}((x\otimes_B \xi)\otimes W_{21}(\eta\otimes \zeta)) \\
&\cong& (\id\otimes \Delta^{\opp})\alpha_X(x)\otimes_{B\bar{\otimes}M\bar{\otimes}M} W_{32}^*U_{\Gsp,12}W_{32}(\xi\otimes \eta\otimes \zeta))\\
&=& (\id\otimes \Delta^{\opp})\alpha_X(x)\otimes_{B\bar{\otimes}M\bar{\otimes}M} (\id\otimes \Delta^{\opp})(U_{\Gsp})(\xi\otimes \eta\otimes \zeta)),
\end{eqnarray*}
so that \eqref{EqCorepToProve} now follows from the coaction property for $\alpha_X$ and the corepresentation property of $U_{\Gsp}$. 

One can now see that $U_{\boxtimes}$ must be unitary since \emph{any} isometric corepresentation of a locally compact quantum group is unitary, by \cite{BDS13}*{Corollary 4.16}. Alternatively, one could more directly observe that the range of $U_{\boxtimes}$ will contain all vectors of the form $\alpha_X(x)(1\otimes m)\otimes_{B\bar{\otimes}M} (\xi\otimes \eta)$ for $x\in X_B(\Hsp),m\in M$ and $\xi\in \Hsp,\eta\in \Gsp$, and then conclude from the fact that the linear span of the $\alpha_X(x)(1\otimes m)$ is $\sigma$-strongly dense in $X_B(\Hsp)\bar{\otimes}M$. 

The equivariance conditions
\[
U_{\boxtimes}(\pi_{\boxtimes}(a)\otimes 1) = (\pi_{\boxtimes}\otimes \pi_M)\alpha(a) U_{\boxtimes},\qquad (\rho_{\boxtimes}(c)\otimes 1)U_{\boxtimes} = U_{\boxtimes}(\rho_{\boxtimes}\otimes \pi_MR)\gamma(c),\qquad a\in A,c\in C
\]
are easily checked.
\end{proof}

By using similar identifications as in the proof of Proposition \ref{PropCrossRep}, one can check that the associativity intertwiner \eqref{EqAssInt} is an intertwiner of $\G$-equivariant correspondences. We refrain from spelling out the details here. 

The following lemma is convenient for computations. 

\begin{Lem}\label{LemSimplForm}
Assume $(\Hsp\otimes L^2(B),U,\pi,\rho)$ is endowed with the structure of a $\G$-$A$-$B$-correspondence such that 
\[
\rho(b)(\xi\otimes \eta) = \xi \otimes \rho_B(b)\eta,\qquad U = U_{\Hsp,13}U_{\beta,23}
\]
for some unitary $\G$-representation $U_{\Hsp}$ on $\Hsp$. If $(\Gsp,U_{\Gsp},\pi_{\Gsp},\rho_{\Gsp})$ is a $\G$-equivariant $B$-$C$-correspondence, we have a unitary intertwiner 
\[
(\Hsp\otimes L^2(B))\boxtimes \Gsp \cong \Hsp \otimes \Gsp,
\]
where the right hand side is endowed with the $\G$-$A$-$C$-correspondence
\begin{equation}\label{EqSimplForm}
\pi_{\boxtimes}(a) = (\id\otimes \pi_{\Gsp})\pi(a),\qquad \rho_{\boxtimes}(c) = 1\otimes \rho_{\Gsp}(c),\qquad a\in A,c\in C,
\end{equation}
\[
U_{\boxtimes} = U_{\Hsp,13}U_{\Gsp,23}.
\]
\end{Lem}
\begin{proof}
Note that the assumption on $\rho$ implies that $\pi(a) \in B(\Hsp) \bar{\otimes} \pi_{B}(B)$ for $a\in A$, so under the identification $\pi_B(B) \cong B$ we see that $\pi_{\boxtimes}$ in \eqref{EqSimplForm} is at least a well-defined unital normal $*$-homomorphism $A \rightarrow B(\Hsp\otimes \Gsp)$.  

Now it is immediate that $X_B(\Hsp \otimes L^2(B)) = \Hsp \overline{\otimes} B$ with its obvious self-dual Hilbert $A$-$B$-bimodule structure, and with the coaction given by 
\[
\alpha_X(\xi\otimes b) = U_{\Hsp,13}(\xi\otimes \beta(b)). 
\]
We hence obtain an isometric map with dense image
\[
(\Hsp \odot B) \odot_B \Gsp \rightarrow \Hsp \otimes \Gsp,\qquad (\xi\otimes b) \otimes \zeta \mapsto \xi\otimes b\zeta.
\]
It is clear that this induces \eqref{EqSimplForm} on $\Hsp \otimes \Gsp$. Also the form for $U_{\boxtimes}$ is immediate.
\end{proof}

\begin{Def}
Let $A,B$ be $\G$-W$^*$-algebras. We say that $A$ is \emph{$\G$-equivariantly W$^*$-Morita equivalent with $B$} if there exists a $\G$-equivariant $A$-$B$-correspondence $(\Hsp,U,\pi,\rho)$ which is full (or equivalently, faithful) as a right $B$-module\footnote{Concretely, this means that $\pi_B(B) = [X_B(\mathcal{H})^*X_B(\mathcal{H})]^{\sigma\text{-weak}}.$} and faithful as a left $A$-module, and such that 
\[
\pi(A) = \rho(B)'.
\]
In this case, we call $\Hsp$ a $\G$-W$^*$-Morita correspondence (between $A$ and $B$).
\end{Def}

If $\Hsp$ is a $\G$-W$^*$-Morita correspondence between $A$ and $B$, we put (recalling the notation \eqref{EqHilbertMod})
\[
Q_{\Hsp} := X_B(L^2(B) \oplus \Hsp,L^2(B) \oplus \Hsp). 
\]
It follows from basic theory of correspondences \cite{Tak03} that we have a decomposition
\[
Q_{\Hsp} \cong \begin{pmatrix} A & X_B(\Hsp)\\ X_B(\Hsp)^* & B \end{pmatrix},
\]
and that we can then decompose 
\begin{equation}\label{EqUnitCorrDec}
L^2(Q_{\Hsp})  = \begin{pmatrix} L^2(A) & L^2(X_B(\Hsp))\\ L^2(X_B(\Hsp)^*) & L^2(B) \end{pmatrix} \cong  \begin{pmatrix} L^2(A) &\Hsp \\ \overline{\Hsp}& L^2(B)\end{pmatrix}
\end{equation}
as an $L^2$-direct sum. The isomorphism $L^2(X_B(\Hsp)) \cong \Hsp$ is uniquely determined by $\pi_Q(x) \xi \cong x\xi$ if $x\in X_B(\Hsp)$ and $\xi\in L^2(B)$. The isomorphism $L^2(X_B(\Hsp)^*) \cong \overline{\Hsp}$ is uniquely determined by $J_Q\xi \cong \overline{\xi}$ for $\xi\in \Hsp$. 

In the sequel, we simply write $Q= Q_{\mathcal{H}}$.
From the fact that we are in the equivariant setting, $Q$ obtains a natural $\G$-action $\alpha_Q$ on $Q$ by conjugating with the unitary representation $U_{\beta} \oplus U$ on $L^2(B) \oplus \Hsp$. Let $U_Q$ be the associated canonical unitary implementing $\alpha_Q$.

\begin{Lem}
Under the identification \eqref{EqUnitCorrDec}, we have 
\begin{equation}\label{EqIdentUniCorep}
U_Q \cong \begin{pmatrix} U_{\alpha} & U \\ \overline{U} & U_{\beta}\end{pmatrix}.
\end{equation}
\end{Lem}
\begin{proof}
It is clear that there will be a decomposition
\[
U_Q \cong \begin{pmatrix} U_{11}& U_{12} \\ U_{21} & U_{22}\end{pmatrix},
\]
as $\begin{pmatrix} 1  & 0 \\0 & 0\end{pmatrix}$ and $\begin{pmatrix} 0  & 0 \\0 & 1\end{pmatrix}$ are invariant under $\alpha_Q$. The identifications on the diagonal of \eqref{EqIdentUniCorep} then follow by noticing that the construction of the unitary $\G$-representation in \cite{Vae05} is compatible with the ($\G$-invariant) projection maps $L^2(Q) \rightarrow L^2(B)$ and $L^2(Q) \rightarrow L^2(A)$.

As $U_Q$ implements $\alpha_Q$, we must have for $x \in X_B(\Hsp)$ and $\xi\in L^2(B),\eta\in L^2(\G)$ that
\[
U_Q(x\xi\otimes \eta) = \alpha_X(x)U_{\beta}(\xi\otimes \eta) \cong U(x\xi\otimes \eta),
\]
showing the isomorphism in the top right corner of \eqref{EqIdentUniCorep}. The isomorphism on the lower left corner is deduced from the observation that 
\[
U_Q(J_Q\xi\otimes \check{J}\eta) = (J_Q\otimes \check{J})U^*(\xi\otimes \eta) \cong (C_{\Hsp}\otimes \check{J})U^*(\xi\otimes \eta) = \overline{U}(C_{\Hsp}\xi\otimes \check{J}\eta),\qquad \xi\in \Hsp,\eta\in L^2(\G).
\]\end{proof}

\begin{Prop}\label{PropMoritaGood}
If $\Hsp$ is a $\G$-W$^*$-Morita correspondence between $A$ and $B$, then 
\begin{equation}\label{EqIdentTensBla}
\Hsp \boxtimes \overline{\Hsp}\cong L^2(A),\qquad \overline{\Hsp}\boxtimes \Hsp \cong L^2(B).
\end{equation}
\end{Prop}
\begin{proof}
We can construct the maps in \eqref{EqIdentTensBla} as equivariant isometries simply by restriction of the natural unitary $\G$-equivariant intertwiners 
\[
L^2(Q) \boxtimes_Q L^2(Q) \cong L^2(Q) \cong L^2(Q)\boxtimes_Q L^2(Q).
\]
It is then well-known from the theory of (non-equivariant) correspondences that these restrictions are still unitary (cf.\ \cite{Sau83}*{Proposition 3.1}).
\end{proof}

\begin{Cor}\label{CorCompOpp}
If $A,B,C$ are $\G$-W$^*$-algebras, $\Hsp$ is a $\G$-equivariant $A$-$B$-correspondence and $\Gsp$ is a $\G$-equivariant $B$-$C$-correspondence, then 
\begin{equation}\label{EqIsoConjTens}
\overline{\Hsp \boxtimes \Gsp} \cong \overline{\Gsp} \boxtimes \overline{\Hsp},
\end{equation}
\end{Cor} 
\begin{proof}
We may replace $\Gsp$ by $\overline{\Gsp}$, for $\Gsp$ a $C$-$B$-correspondence. By considering $A\oplus C$ and $\Hsp \oplus \Gsp$, we may moreover assume that $A= C$ and $\Hsp = \Gsp$. Furthermore, we can then as well assume that $\rho_{\Hsp}$ is faithful, and $A = \rho_{\Hsp}(B)'$. In other words, $\Hsp$ is now a $\G$-W$^*$-Morita correspondence between $A$ and $B$, and the result follows from Proposition \ref{PropMoritaGood}.
\end{proof}

\begin{Cor}\label{CorInclIsStandard}
Let $A \subseteq B$ be an equivariant inclusion, and put $C = \rho_B(A)'$. Then $C$ is $\G$-equivariant by 
\[
\gamma(c) := U_{\beta}(c\otimes1)U_{\beta}^*,\qquad c\in C,
\]
and we obtain a $\G$-W$^*$-Morita $C$-$A$-correspondence ${}_CL^2(B)_A$. We then have
\[
{}_BL^2(B)_A \boxtimes {}_AL^2(B)_B \cong {}_BL^2(C)_B. 
\]
\end{Cor}
\begin{proof}
The result is just a particular case of Proposition \ref{PropMoritaGood}.
\end{proof}

\begin{Rem}\label{RemSauv}
More generally, we have that if $\Hsp$ is a $\G$-equivariant $A$-$B$-correspondence and $C = \rho(B)'$, we can identify ${}_A\Hsp \boxtimes \overline{\Hsp}_{A} \cong {}_AL^2(C)_A$. The proof easily reduces to the above case.
\end{Rem}

Let us look at some more particular examples. We again refer back to Section \ref{SecIntroAndExa} for the notations we use.

\begin{Exa}\label{EqUnitUnderTens}
Using Lemma \ref{LemSimplForm} and Corollary \ref{CorCompOpp}, it is easily seen that $E_A^{\G}$ acts as an identity under $\boxtimes$: for any $\Hsp \in \Corr^{\G}(A,B)$ we obtain canonical isomorphisms
\[
E_A^{\G}\boxtimes \Hsp \cong \Hsp \cong \Hsp \boxtimes E_B^{\G}.
\]
\end{Exa}

\begin{Exa}\label{ExaRestrComp}
If $A \subseteq B$ is an equivariant inclusion and $\Hsp$ is a $\G$-equivariant $B$-$C$-correspondence, we have 
\[
{}_AL^2(B)_B \boxtimes {}_B\Hsp_C \cong {}_A\Hsp_C. 
\] 
Similarly, if $\Hsp$ is a $\G$-equivariant $C$-$B$-correspondence, then 
\[
{}_C\Hsp_B \boxtimes {}_BL^2(B)_A \cong {}_C\Hsp_A. 
\]
\end{Exa}

\begin{Exa}\label{ExaW*semiCoars}
We can identify
\[
S_{A,B}^{\G} \cong {}_AL^2(A)_{\C}\boxtimes {}_{\C}L^2(B)_B,
\]
as follows directly by definition of the Connes fusion product.
\end{Exa}

\begin{Exa}\label{ExaCompGivesCoars}
Combining Example \ref{EqUnitUnderTens}, Example \ref{ExaRestrComp} and Example \ref{ExaW*semiCoars} with the associativity of $\boxtimes$, we see that we can identify
\[
S_{A,B}^{\G}\boxtimes S_B^{\G} \cong L^2(A)\otimes L^2(B)\otimes L^2(\G),
\]
where the latter is endowed with the normal (anti-)representations
\[
\pi(a) = \pi_A(a)\otimes 1\otimes 1,\qquad \rho(b) = 1\otimes (\rho_B\otimes \rho_M)\beta(b),\qquad a\in A,b\in B,
\]
and the unitary $\G$-representation $U_{\alpha,14}V_{34}$. Conjugating with $\Sigma_{23}U_{\alpha,13}$, we find that 
\[
S_{A,B}^{\G}\boxtimes S_B^{\G} \cong C_{A,B}^{\G}.
\]
It follows from the computations in Example \ref{ExaOppExasCan} and Corollary \ref{CorCompOpp} that also 
\[
S_{A}^{\G}\boxtimes S_{A,B}^{\G} \cong C_{A,B}^{\G}.
\]
\end{Exa}

\begin{Exa}\label{ExaCompoEasy}
Similar as in the previous Example, we obtain that
\[
S_{A,B}^{\G}\boxtimes S_{B,C}^{\G} \cong L^2(A)\otimes L^2(B)\otimes L^2(C)
\]
endowed with the $A$-$C$-correspondence by 
\[
\pi(a) = \pi_A(a)\otimes 1\otimes 1,\qquad \rho(c) = 1\otimes 1 \otimes \rho_C(c),\qquad a\in A,c\in C,
\]
and the unitary $\G$-representation $U_{\alpha,14}U_{\beta,24}U_{\gamma,34}$. In particular, if we take $B = M$ with $\beta = \Delta$, we see that 
\[
S_{A,M}^{\G}\boxtimes S_{M,C}^{\G} \cong C_{A,C}^{\G}.
\] 
Combining this with Example \ref{ExaResRestr}, we can write this as
\[
S_{A,M}^{\G}\boxtimes \overline{S_{C,M}^{\G}} \cong C_{A,C}^{\G}.
\]
\end{Exa}
\begin{Exa}\label{ExaRegAndAdjFusion}
Consider the $\G$-regular $\G$-$A$-$A$-correspondence $S^A_{\G}$. Since $S_{A}^{\G}$ is the restriction of the identity correspondence $E_{A\bar{\otimes} M}$ to $A\cong \alpha(A)$, we see from Example \ref{ExaRestrComp} and Example \ref{ExaResRestr} that 
\[
S_A^{\G} \cong \Hsp \boxtimes\overline{\Hsp},\qquad\Hsp  = {}_{\alpha(A)}(L^2(A)\otimes L^2(\G))_{A\bar{\otimes}M}.
\]
A similar reasoning applies to the adjoint $\G$-$A$-$A$-correspondence $D^A_{\G}$ as introduced in Example \ref{Exa5}: since it is the restriction of the identity correspondence $E_{A\rtimes \G}$ to $A\cong \alpha(A)$, where $A\rtimes \G$ is endowed with the adjoint $\G$-action, we find that 
\[
D_A^{\G} \cong \Hsp \boxtimes \overline{\Hsp},\qquad\Hsp  = {}_{\alpha(A)}(L^2(A)\otimes L^2(\G))_{A\rtimes \G}.
\]
\end{Exa}

\begin{Exa}\label{ExaSemiAdjTensss}
It is easily computed that the regular amplification in Example \ref{ExaSemiCoarsG} is obtained by tensoring with the $\G$-semi-coarse correspondence: if $\Hsp$ is a $\G$-equivariant $A$-$B$-correspondence, then 
\[
S^{\G}(\Hsp) \cong \Hsp \boxtimes S_{B}^{\G} \cong S_A^{\G}\boxtimes \Hsp.
\]
Similarly, the adjoint amplification of $\Hsp$ is gotten by tensoring with the adjoint $\G$-correspondence, 
\[
D^{\G}(\Hsp) \cong \Hsp \boxtimes D_B^{\G}\cong D_A^{\G}\boxtimes \Hsp. 
\]
\end{Exa}

\subsection{Crossed products}

We lift the crossed product construction from actions to the level of equivariant correspondences. This is an operation which is non-existent (or more precisely trivial) in the case of non-equivariant correspondences.

\begin{Prop}\label{PropDualities}
Let $(\Hsp,U,\pi,\rho)$ be a $\G$-$A$-$B$-correspondence. Then 
\[
\Hsp^{\rtimes} = \Hsp \rtimes \G := \Hsp \otimes L^2(\G)
\]
becomes a $\check{\G}$-$A^{\rtimes}$-$B^{\rtimes}$-correspondence for the following structures:  
\begin{itemize}
\item the normal $A^{\rtimes}$-representation $\pi^{\rtimes}$ determined by
\[
\pi^{\rtimes}(z) = (\pi\otimes \id)(z),\qquad z\in A^{\rtimes},
\]
\item the normal anti-$*$-representation $\rho^{\rtimes}$ of $B^{\rtimes}$ given by  
\[
\rho^{\rtimes}(z) = U(\rho\otimes \check{J}(-)^*\check{J})(z)U^*,\qquad z\in B^{\rtimes},
\]
\item the unitary $\check{\G}$-representation $U_{\rtimes} := \wW_{23}$. 
\end{itemize}
\end{Prop}
\begin{proof}
It is easily seen that $\pi^{\rtimes}$ and $\rho^{\rtimes}$ indeed define respectively a normal $*$-representation and a normal anti-$*$-representation of $A^{\rtimes}$ and $B^{\rtimes}$. Note that the ranges of $\pi^{\rtimes}$ and $\rho^{\rtimes}$ commute: indeed, we have
\begin{equation}
\pi^{\rtimes}(\alpha(a)) = (\pi\otimes \pi_M)\alpha(a),\qquad \pi^{\rtimes}(1\otimes x) = 1\otimes \pi_{\check{M}}(x),\qquad a \in A,x\in \check{M},
\end{equation}
\begin{equation}
\rho^{\rtimes}(\beta(b)) = \rho(b)\otimes 1,\qquad \rho^{\rtimes}(1\otimes x) = U(1\otimes \rho_{\check{M}}(x))U^*,\qquad b \in B,x\in \check{M},
\end{equation}
so the only non-trivial commutation to check is the one between $U(1\otimes \rho_{\check{M}}(\check{M}))U^*$ and $1\otimes \pi_{\check{M}}(\check{M})$. But this follows easily from the identity $U_{12}^*V_{23}U_{12} = U_{13}V_{23}$. 

It is clear that $U_{\rtimes}$ interacts in the correct way with $\pi^{\rtimes}$. To see its $\rho^{\rtimes}$-compatibility, simply note that 
\[
\wW_{23}^*U_{12} = U_{12}\wW_{23}^*. 
\]\end{proof}

We now determine how the above crossed product constructions interact with opposites and composition.
\begin{Prop}\label{PropOpposit}
Let $\Hsp$ be a $\G$-$A$-$B$-correspondence. We have 
\[
\overline{\Hsp}^{\rtimes} \cong  \overline{\Hsp^{\rtimes}}.
\]
\end{Prop} 
\begin{proof}
The unitary
\[
\overline{\Hsp}\otimes L^2(\G)\rightarrow \overline{\Hsp\otimes L^2(\G)},\qquad \overline{\xi}\otimes \eta\mapsto \overline{U(\xi\otimes \check{J}\eta)}
\]
does the job.
\end{proof}

In the proof of the following theorem, it will be convenient to have a different model for the crossed product correspondence. If $\Hsp$ is a $\G$-$A$-$B$-correspondence, we will denote by $\Hsp^{\blacktriangleleft}$ the $\check{\G}$-$A^{\rtimes}$-$B^{\rtimes}$-correspondence $\Hsp\otimes L^2(\G)$ obtained by conjugating $\Hsp^{\rtimes}$ with $U^*$, so in particular 
\[
\pi^{\blacktriangleleft}(\alpha(a)) = \pi(a)\otimes 1,\quad \pi^{\blacktriangleleft}(1\otimes x) = U^*(1\otimes \pi_{\check{M}}(x))U,\qquad a\in A,x\in \check{M},
\]
\[
\rho^{\blacktriangleleft}(\beta(b)) = (\rho\otimes \pi_MR)\beta(b),\quad \rho^{\blacktriangleleft}(1\otimes x) = 1\otimes \rho_{\check{M}}(x),\qquad b\in B,x\in \check{M}, 
\]
while $U_{\blacktriangleleft} = \wW_{23}$ is still the associated unitary $\check{\G}$-representation. Note that under the above representation, we have for example that 
\begin{equation}\label{EqCompuUUUU}
(\pi^{\blacktriangleleft}\otimes \id)(V_{23})= U_{12}^*V_{23}U_{12} = U_{13}V_{23}.
\end{equation}

\begin{Theorem}\label{TheoCompTTDual}
Let $(\Hsp,U_{\Hsp},\pi_{\Hsp},\rho_{\Hsp})$ be a $\G$-$A$-$B$-correspondence, and let $(\Gsp,U_{\Gsp},\pi_{\Gsp},\rho_{\Gsp})$ be a $\G$-$B$-$C$-correspondence. Then we have 
\[
(\Hsp\boxtimes \Gsp)^{\rtimes} \cong  \Hsp^{\rtimes}\boxtimes \Gsp^{\rtimes}.
\]
\end{Theorem}
\begin{proof}
We will show rather that 
\[
(\Hsp\boxtimes \Gsp)^{\blacktriangleleft}\cong \Hsp^{\rtimes}\boxtimes \Gsp^{\rtimes}
\]
Indeed, write $X = X_B(\Hsp)$ and consider the map
\begin{equation}\label{EqEquivToEst}
(X \odot_B \Gsp) \odot L^2(\G) \rightarrow \Hsp^{\rtimes}\boxtimes \Gsp^{\rtimes},\quad (x \otimes_B \xi)\otimes \eta \mapsto \alpha_X(x) \otimes_{B^{\rtimes}} U_{\Gsp}(\xi\otimes \eta),
\end{equation}
where we have used that $\alpha_X(X) \subseteq X_{B^{\rtimes}}(\Hsp^{\rtimes})$. A direct calculation shows that this is a well-defined isometric map, which then extends to an isometry 
\begin{equation}\label{EqEquivToEst2}
T:(\Hsp\boxtimes \Gsp)\otimes L^2(\G) \rightarrow \Hsp^{\rtimes}\boxtimes \Gsp^{\rtimes}.
\end{equation}
In fact, it is an easy observation that $X_{B^{\rtimes}}(\Hsp^{\rtimes})$ equals the $\sigma$-strong closure of $\alpha_X(X)(1\otimes \check{M})$, hence \eqref{EqEquivToEst2} will be unitary.

It is immediate that $T$ is an intertwiner of $A$-representations and of right $\check{M}$-representations. 

To see that $T$ is an intertwiner of $\check{\G}$-representations, we compute
\begin{eqnarray*}
U_{\rtimes \boxtimes}((\alpha_X(x) \otimes_{B^{\rtimes}}U_{\Gsp}(\xi\otimes \eta))\otimes \zeta) &=& (\alpha_X(x)\otimes 1) \otimes_{B^{\rtimes}\bar{\otimes} \check{M}} \widetilde{W}_{23}U_{\Gsp,12}(\xi\otimes \eta \otimes \zeta) \\
&=& (\alpha_X(x)\otimes 1) \otimes_{B^{\rtimes}\bar{\otimes} \check{M}}U_{\Gsp,12}(\xi\otimes \widetilde{W}(\eta \otimes \zeta)),
\end{eqnarray*}
from which we see that 
\[
U_{\rtimes\boxtimes}(\omega)(\alpha_X(x) \otimes_{B^{\rtimes}}U_{\Gsp}(\xi\otimes \eta)) =  \alpha_X(x) \otimes_{B^{\rtimes}}U_{\Gsp}(\xi\otimes (\id\otimes \omega)(\widetilde{W})\eta),\qquad \omega \in B(L^2(\G))_*.
\]
Hence $TU_{\boxtimes \blacktriangleleft}(\omega) =U_{\rtimes \boxtimes}(\omega)T$, and so $T$ intertwines the two $\check{\G}$-representations.  

The intertwining property for right $C$-representations can be obtained by observing first the identification 
\[
(\Hsp\boxtimes \Gsp)\otimes L^2(\G) \cong (\Hsp\otimes L^2(\G))\underset{B\bar{\otimes}M}{\boxtimes} (\Gsp\otimes L^2(\G)),\qquad (x\otimes_B \xi)\otimes \eta \mapsto (x\otimes 1)\otimes_{B\bar{\otimes} M} (\xi\otimes \eta). 
\]
Under this identification, the map $T$ can be written as 
\[
T: (\Hsp\otimes L^2(\G))\boxtimes (\Gsp\otimes L^2(\G)) \rightarrow \Hsp^{\rtimes}\boxtimes \Gsp^{\rtimes},\qquad (x\otimes 1_M) \otimes_{B\bar{\otimes}M} \xi \mapsto \alpha_X(x)\otimes_{B^{\rtimes}} U_{\Gsp}\xi, 
\]
where $x\in X_B(\Hsp)$ and $\xi\in \Gsp\otimes L^2(\G)$. With this description, it is immediate that $T$ intertwines the right $C$-representations as well. 

Finally, we need to check that $T$ is an intertwiner of left $\check{M}$-representations. For this, we note that we can identify 
\[
(\Hsp\boxtimes \Gsp)\otimes L^2(\G)\otimes L^2(\G) \cong (\Hsp\otimes L^2(\G)\otimes L^2(\G))\underset{B\bar{\otimes}\check{M}\bar{\otimes}M}{\boxtimes} (\Gsp\otimes L^2(\G)\otimes L^2(\G)),
\]
\[
(x \otimes \xi)\otimes \eta\otimes \zeta \mapsto (x\otimes 1\otimes 1)\otimes_{B\bar{\otimes}\check{M}\bar{\otimes}M} (\xi\otimes \eta\otimes \zeta),
\]
and similarly we can identify 
\[
(\Hsp^{\rtimes}\underset{B^{\rtimes}}{\boxtimes} \Gsp^{\rtimes}) \otimes L^2(\G) \cong (\Hsp^{\rtimes}\otimes L^2(\G))\underset{B^{\rtimes}\bar{\otimes}M}{\boxtimes} (\Gsp^{\rtimes} \otimes L^2(\G)),
\]
\[
(x \otimes_{B^{\rtimes}} \xi)\otimes \eta \mapsto (x\otimes 1)\otimes_{B^{\rtimes}\bar{\otimes}M} (\xi\otimes \eta).
\]
Then $T\otimes \id_{L^2(\G)}$ can be viewed as the unitary 
\[
(\Hsp\otimes L^2(\G)\otimes L^2(\G))\underset{B\bar{\otimes}\check{M}\bar{\otimes}M}{\boxtimes} (\Gsp\otimes L^2(\G)\otimes L^2(\G)) \rightarrow (\Hsp^{\rtimes}\otimes L^2(\G))\underset{B^{\rtimes}\bar{\otimes}M}{\boxtimes} (\Gsp^{\rtimes} \otimes L^2(\G)),
\]
\[
z_{13}\otimes_{B\bar{\otimes}\check{M}\bar{\otimes}M} \xi \mapsto (\alpha\otimes \id)(z) \bar{\otimes}_{B^{\rtimes}\bar{\otimes} M} U_{\Gsp,12}\xi,
\]
where $z \in X_B(\Hsp)\bar{\otimes}M$ and $\xi\in \Gsp\otimes L^2(\G)\otimes L^2(\G)$.

It is now sufficient to prove that, under these models,
\begin{equation}\label{EqHelp}
(T\otimes 1)(\pi^{\boxtimes\blacktriangleleft}\otimes \id)(V) = (\pi^{\rtimes \boxtimes}\otimes \id)(V)(T\otimes 1).
\end{equation}
But since by \eqref{EqCompuUUUU} we have
\[
(\pi^{\boxtimes\blacktriangleleft}\otimes \id)(V) = U_{\boxtimes,13}V_{23},
\]
we see that \eqref{EqHelp} follows from the following computation for $x\in X_B(\Hsp)$ and $\xi \in \Gsp\otimes L^2(\G)\otimes L^2(\G)$: 
\begin{eqnarray*} 
&& \hspace{-2cm} (T\otimes 1)(\pi^{\boxtimes\blacktriangleleft}\otimes \id)(V)((x\otimes 1\otimes 1)\otimes_{B\bar{\otimes}\check{M}\bar{\otimes}M}\xi) \\
&=& (T\otimes1) U_{\boxtimes,13}((x\otimes 1 \otimes 1)\otimes_{B\bar{\otimes}\check{M}\bar{\otimes}M} V_{23}\xi) \\
&=& (T\otimes1) (\alpha_X(x)_{13}\otimes_{B\bar{\otimes}\check{M}\bar{\otimes}M} U_{\Gsp,13}V_{23}\xi) \\
&=& (\alpha\otimes \id)\alpha(x)\otimes_{B^{\rtimes}\bar{\otimes}M} U_{\Gsp,12}U_{\Gsp,13}V_{23}\xi\\
&=& V_{23}\alpha(x)_{12}\otimes_{B^{\rtimes}\bar{\otimes}M} V_{23}^*U_{\Gsp,12}U_{\Gsp,13}V_{23}\xi\\
&=& V_{23}\alpha(x)_{12} \otimes_{B^{\rtimes}\bar{\otimes}M}U_{\Gsp,12}\xi\\
&=& (\pi^{\rtimes \boxtimes}\otimes \id)(V)(\alpha(x)_{12} \otimes_{B^{\rtimes}\bar{\otimes}M}U_{\Gsp,12}\xi)\\
&=& (\pi^{\rtimes \boxtimes}\otimes \id)(V)(T\otimes 1)((x\otimes 1\otimes 1)\otimes_{B\bar{\otimes}\check{M}\bar{\otimes}M}\xi).
\end{eqnarray*}\end{proof}

\subsection{Takesaki-Takai duality}

Recall that $M^{\vee\vee}\cong M$ via $\Ad(u_{\G})$. Hence any $\G$-W$^*$-algebra $A$ can be seen as a $\G^{\vee\vee}$-W$^*$-algebra by a coaction that we write 
\[
\widetilde{\alpha}: A \rightarrow A \bar{\otimes} M^{\vee\vee},\qquad a \mapsto (1\otimes u_{\G})\alpha(a)(1\otimes u_{\G}). 
\]
The \emph{Takesaki-Takai duality} for actions of locally compact quantum groups on von Neumann algebras \cites{ES80,Vae01} can now be stated as follows. 
\begin{Theorem}[Takesaki-Takai duality]\label{TheoTakTakEq}
Endow $A\bar{\otimes} B(L^2(\G))$ with the $\G$-action
\[
\alpha_+ = \Ad(\Sigma_{23}W_{23})\circ (\alpha\otimes \id).
\] 
Then there is a $\G^{\vee\vee}$-equivariant $*$-isomorphism 
\begin{equation}\label{EqTakTakIso}
\Phi: (A^{\rtimes \rtimes},\alpha^{\rtimes\rtimes}) \cong (A\bar{\otimes} B(L^2(\G)),\widetilde{\alpha_+}),
\end{equation}
whose inverse is concretely implemented by 
\[
\Phi^{-1}: A\bar{\otimes}B(L^2(\G))\rightarrow A^{\rtimes\rtimes}\subseteq B(L^2(A)\otimes L^2(\G)\otimes L^2(\G)),\quad z\mapsto V_{23}^*(\alpha\otimes \id)(z)V_{23}. 
\]
On generating elements, $\Phi$ is determined by 
\[
\alpha(a)\otimes 1 \mapsto \alpha(a),\quad 1\otimes \check{\Delta}(x)\mapsto 1\otimes x,\quad 1\otimes 1 \otimes y\mapsto 1\otimes y,\qquad a\in A,x\in \check{M},y\in M', 
\]
\end{Theorem}

We reformulate this result as follows. 

\begin{Theorem}\label{TheoTakTakTak}
Endow the Hilbert space ${}_{A^{\rtimes\rtimes}}\msM_A := L^2(A)\otimes L^2(\G)$ with 
\begin{itemize}
\item the normal unital $*$-representation $(\pi_A\otimes \id)\circ \Phi$ of $A^{\rtimes\rtimes}$, 
\item the normal unital anti-$*$-representation $\rho_A\otimes 1$ of $A$, and
\item the $\G^{\vee\vee}$-representation $u_{\G,3}W_{32}U_{\alpha,13}u_{\G,3}$. 
\end{itemize}
Then ${}_{A^{\rtimes\rtimes}}\msM_A$ becomes a $\G^{\vee\vee}$-W$^*$-Morita correspondence between $(A^{\rtimes\rtimes},\alpha^{\rtimes\rtimes})$ and $(A,\widetilde{\alpha})$.
\end{Theorem} 

Note now that if $(A,\alpha),(B,\beta)$ are $\G$-W$^*$-algebras, any $\G$-$A$-$B$-correspondence $(\Hsp,\pi,\rho,U)$ becomes a $\G^{\vee\vee}$-equivariant $(A,\widetilde{\alpha})$-$(B,\widetilde{\beta})$-correspondence by means of 
\[
\widetilde{\Hsp} = (\Hsp,\pi,\rho,(1\otimes u_{\G})U(1\otimes u_{\G})). 
\]
We can then extend the Takesaki-Takai duality to general $\G$-equivariant correspondences.

\begin{Theorem}\label{TheoTakTakTakTak}
Let $(\Hsp,U_{\Hsp},\pi_{\Hsp},\rho_{\Hsp})$ be a $\G$-equivariant $A$-$B$-correspondence. Then we have a unitary intertwiner of $\G^{\vee\vee}$-$A$-$B$-correspondences  
\[
\Hsp^{\rtimes\rtimes} \cong {}_{A^{\rtimes\rtimes}}\msM_A\boxtimes \widetilde{\Hsp} \boxtimes {}_B\overline{\msM}_{B^{\rtimes\rtimes}}. 
\]
\end{Theorem}
\begin{proof}
Since $\msM$ is a $\G^{\vee\vee}$-W$^*$-Morita correspondence, it is (by associativity of $\boxtimes$ and Proposition \ref{PropMoritaGood}) enough to show that 
\begin{equation}\label{EqEasierMade}
{}_{A^{\rtimes\rtimes}}\Hsp^{\rtimes\rtimes}_{B^{\rtimes\rtimes}} \boxtimes {}_{B^{\rtimes\rtimes}}\msM_{B} \cong {}_{A^{\rtimes\rtimes}}\msM_A\boxtimes {}_A\widetilde{\Hsp}_B.
\end{equation}
Now the representation of $A\rtimes \G\rtimes \check{\G}$ and anti-representation of $B\rtimes \G\rtimes \check{\G}$ on $\Hsp^{\rtimes\rtimes}= \Hsp\otimes L^2(\G)\otimes L^2(\G)$ are given by 
\[
\pi^{\rtimes\rtimes}(\alpha(a)\otimes 1) = (\pi_{\Hsp}\otimes \pi_M)(\alpha(a))_{12},\quad \rho^{\rtimes\rtimes}(\beta(b)\otimes 1) = \rho_{\Hsp}(b)_1,\qquad a\in A,b\in B,
\]
\[
\pi^{\rtimes\rtimes}(1\otimes \check{\Delta}(y)) = \check{\Delta}(y)_{23},\quad \rho^{\rtimes\rtimes}(1\otimes \check{\Delta}(y)) = U_{\Hsp,12}\rho_{\check{M}}(y)_2U_{\Hsp,12}^*,\qquad y\in \check{M},
\]
\[
\pi^{\rtimes\rtimes}(1\otimes 1 \otimes z) = \pi_{M^{\vee\vee}}(z)_3,\quad \rho^{\rtimes\rtimes}(1\otimes 1 \otimes z) = \widetilde{W}_{23}\rho_{M^{\vee\vee}}(z)_3\widetilde{W}_{23}^*,\qquad z\in M^{\vee\vee},
\]
with associated unitary $\G^{\vee\vee}$-representation given by $\widetilde{V}_{34}$, where we recall that $\widetilde{V}= \Ad(u_{\G}\otimes u_{\G})(V)$. 

It is then easily verified that 
\[
X_{B^{\rtimes\rtimes}}(\Hsp^{\rtimes\rtimes}) = \sigma\textrm{-closure of the span of } \{\Delta_X(x)_{12}\check{\Delta}(y)_{23}\pi_{M^{\vee\vee}}(z)_3\mid x\in X_B(\Hsp), y\in \check{M}, z \in M^{\vee\vee}\}.
\]
From this, we find a unitary map 
\begin{equation}\label{EqTheUnitTransfooo}
\Hsp^{\rtimes\rtimes}\underset{B^{\rtimes\rtimes}}{\boxtimes} (L^2(B)\otimes L^2(\G))\rightarrow \Hsp \otimes L^2(\G),\quad \left(\Delta_X(x)_{12}\check{\Delta}(y)_{23}\pi_{M^{\vee\vee}}(z)_3\right)\otimes_{B^{\rtimes\rtimes}} \xi \mapsto \Delta_X(x)(1\otimes yz)\xi,
\end{equation}
which hence by transport of structure endows  $\Hsp\otimes L^2(\G)$ with the structure of a $\G^{\vee\vee}$-equivariant $A^{\rtimes\rtimes}$-$B$-correspondence $(\Hsp\otimes L^2(\G),U,\pi,\rho)$. 

Clearly $\pi,\rho$ are given by 
\[
\pi(z) = (\pi_{\Hsp}\otimes \id_{B(L^2(\G))})\Phi^{-1}(z),\qquad \rho(b) = \rho_{\Hsp}(b)\otimes 1,\qquad z\in A^{\rtimes\rtimes},b\in B.
\]
The resulting unitary $\G^{\vee\vee}$-representation $U$ on $\Hsp\otimes L^2(\G)$ is given by the following transformation (dropping the notations $\pi_{\check{M}}$ and $\pi_{M^{\vee\vee}}$): for $x\in X_B(\Hsp),y\in \check{M},z\in M^{\vee\vee}$ and $\xi,\eta\in L^2(\G)$, we have
\begin{eqnarray*}
(\Delta_X(x)_{12}(1\otimes yz)\xi)\otimes \eta &\mapsto& \Delta_X(x)_{12}y_2\Delta^{\vee\vee}(z)_{23}u_{\G,3}W_{32}U_{\beta,13}u_{\G,3}(\xi\otimes \eta) \\
&=& \Delta_X(x)_{12}y_2\Delta^{\vee\vee}(z)_{23}\widetilde{V}_{23}u_{\G,3}U_{\beta,13}u_{\G,3}(\xi\otimes \eta) \\ &=& \Delta_X(x)_{12}y_2\widetilde{V}_{23}z_2u_{\G,3}U_{\beta,13}u_{\G,3}(\xi\otimes \eta)  \\
&=& \Delta_X(x)_{12}\widetilde{V}_{23}u_{\G,3}U_{\beta,13}u_{\G,3}y_2z_2(\xi\otimes \eta)  \\
&=&u_{\G,3} \Delta_X(x)_{12}W_{32}U_{\beta,13}u_{\G,3}y_2z_2(\xi\otimes \eta) \\
&=&u_{\G,3}W_{32} (\id\otimes \Delta^{\opp})\Delta_X(x)U_{\beta,13}u_{\G,3}y_2z_2(\xi\otimes \eta)  \\
&=&u_{\G,3}W_{32} (\Delta_X\otimes \id)\Delta_X(x)_{132}U_{\beta,13}u_{\G,3}y_2z_2(\xi\otimes \eta) \\ 
&=& u_{\G,3}W_{32}U_{\Hsp,13}u_{\G,3}((\Delta_X(x)_{12}y_2z_2\xi)\otimes \eta).
\end{eqnarray*}
We hence find that \eqref{EqTheUnitTransfooo} endows $\Hsp\otimes L^2(\G)$ with the unitary $\G^{\vee\vee}$-representation $U = u_{\G,3}W_{32}U_{\Hsp,13}u_{\G,3}$. But this gives precisely the right hand side of \eqref{EqEasierMade} under the direct identification 
\[
{}_{A^{\rtimes\rtimes}}(L^2(A)\otimes L^2(\G))_A\boxtimes {}_A\widetilde{\Hsp}_B \cong {}_{A^{\rtimes\rtimes}}(\widetilde{\Hsp}\otimes L^2(\G))_{B\otimes 1}.
\]\end{proof}

\begin{Rem}
The above results can be upgraded slightly by noting that we obtain a W$^*$-bicategory $\Corr^{\G}$ whose objects are $\G$-equivariant W$^*$-algebras, whose $1$-cells are $\G$-equivariant correspondences, and whose $2$-cells are unitary intertwiners. Composition is determined by the Connes fusion product. Takesaki-Takai duality provides then an equivalence of W$^*$-bicategories between $\Corr^{\G}$ and $\Corr^{\check{\G}}$. We do not enter into the details here. 
\end{Rem}

\subsection{Examples}\label{SubSecExa}

In our next few results, we look at how the above duality interacts with our canonical $\G$-equivariant correspondences, see Example \ref{Exa1} and further. 

\begin{Prop}\label{PropIdRestr}
We have 
\begin{equation}\label{EqDualUnit}(E_A^{\G})^{\rtimes} \cong E_{A^{\rtimes}}^{\check{\G}}
\end{equation}
\end{Prop}
\begin{proof}
This follows immediately from Theorem \ref{TheoDualCoact}. 
\end{proof}

\begin{Prop}\label{PropRestrBothBoth}
Let $A \subseteq C$ and $B \subseteq D$ be $\G$-equivariant inclusions. Then for any $\G$-$C$-$D$-correspondence $\Hsp$ we have 
\begin{equation}\label{EqRestriction}
({}_A\Res_B(\Hsp))^{\rtimes} = {}_{A^{\rtimes}}\Res_{B^{\rtimes}}(\Hsp^{\rtimes}).
\end{equation}
\end{Prop} 
Note that e.g.\ the inclusion $A\bar{\otimes}B(L^2(\G)) \subseteq C\bar{\otimes} B(L^2(\G))$ descends to an inclusion $A^{\rtimes} \subseteq C^{\rtimes}$, so the restrictions in the right hand side of \eqref{EqRestriction} are meaningful.
\begin{proof}
Immediate.
\end{proof}

\begin{Prop}\label{PropRegAdj}
Duality interchanges the regular amplification and the adjoint amplification: if $\Hsp$ is a $\G$-equivariant $A$-$B$-correspondence, then 
\begin{equation}\label{EqDualSemiG}
(D^{\G}(\Hsp))^{\rtimes}  \cong S^{\check{\G}}(\Hsp^{\rtimes}), \qquad S^{\G}(\Hsp)^{\rtimes}  \cong D^{\check{\G}}(\Hsp^{\rtimes}).
\end{equation}
\end{Prop}
\begin{proof}
By Example \ref{ExaSemiAdjTensss} and Theorem \ref{TheoCompTTDual}, it is sufficient to show that $(D^{\G}_A)^{\rtimes} \cong S_{A^{\rtimes}}^{\check{\G}}$ and $(S^{\G}_A)^{\rtimes} \cong D_{A^{\rtimes}}^{\check{\G}}$. Furthermore, by the observations in Example \ref{Exa3} and Example \ref{Exa5} that $S_A^{\G}$ and $D_A^{\G}$ are restrictions of identity correspondences, it is by Proposition \ref{PropIdRestr} and Proposition \ref{PropRestrBothBoth} sufficient to show that there are $\check{\G}$-equivariant $*$-isomorphisms 
\begin{equation}\label{EqToProveIso}
A^{\rtimes} \rtimes_{\alpha_{ad}} \G \cong A^{\rtimes} \bar{\otimes} \check{M},\qquad 
(A\bar{\otimes} M)\rtimes_{\id\otimes \Delta} \G\cong A^{\rtimes}\rtimes_{\alpha^{\rtimes}} \check{\G},
\end{equation}
with in each case $\alpha(A) \rtimes \G$ carried onto $\alpha^{\rtimes}(A^{\rtimes})$. Note that all of the above von Neumann algebras are (by their canonical constructions) implemented on $L^2(A)\otimes L^2(\G)\otimes L^2(\G)$. 

Now $A^{\rtimes} \rtimes_{\alpha_{ad}} \G$ is generated by $V_{23}(\alpha(A)\otimes1)V_{23}^*$, $V_{23}(1\otimes \check{M}\otimes 1)V_{23}^*$ and $1\otimes 1\otimes \check{M}$, with coaction implemented by the unitary representation $\widetilde{W}_{34}$. Upon conjugating with $V_{23}^*$, the unitary representation remains unchanged, while we now get the von Neumann algebra generated by $\alpha(A)\otimes 1$, $1\otimes \check{M}\otimes 1$ and $1\otimes \check{\Delta}(\check{M})$. The second and third term generate $\check{M}\bar{\otimes}\check{M}$, so all the terms together generate $A^{\rtimes}\bar{\otimes} \check{M}$. Meanwhile, the first and third term together generate $\alpha^{\rtimes}(A^{\rtimes})$. Since $\widetilde{W}_{34}$ is still an implementing unitary for the $\check{\G}$-action on $A^{\rtimes} \bar{\otimes} \check{M}$, this proves the existence of the first isomorphism in \eqref{EqToProveIso}. 

Turning now to the second isomorphism in \eqref{EqToProveIso}, note that the right hand side can be identified with $A\bar{\otimes} B(L^2(\G))$ under the Takesaki-Takai-duality isomorphism \eqref{EqTakTakIso}. Recalling that $\check{V}=\widetilde{W}$, it is immediately verified that the $\check{\G}$-adjoint action on $A^{\rtimes}\rtimes \check{\G}$ gets turned into the coaction 
\begin{equation}\label{EqIsoToAdj}
A\bar{\otimes} B(L^2(\G)) \rightarrow A\bar{\otimes} B(L^2(\G))\bar{\otimes} \check{M},\quad z \mapsto \widetilde{W}_{23}z_{12}\widetilde{W}_{23}^*.
\end{equation}
On the other hand, the left hand side of the second isomorphism in \eqref{EqToProveIso} is generated by $A\otimes 1 \otimes 1$, $1\otimes \Delta(M)$ and $1\otimes 1\otimes \check{M}$, which under conjugation with $W_{23}$, and the observation that $M$ and $\check{M}$ generate $B(L^2(\G))$, becomes $A\bar{\otimes}1\bar{\otimes}B(L^2(\G)) \cong A\bar{\otimes}B(L^2(\G))$. Under this identification, the dual coaction gets exactly identified with \eqref{EqIsoToAdj}. It is also immediate that under these identifications, the von Neumann algebras generated by the copies of $A$ and $\check{M}$ in respectively $(A\bar{\otimes} M)\rtimes_{\id\otimes \Delta} \G$ and $A^{\rtimes}\rtimes_{\alpha^{\rtimes}} \check{\G}$ get sent to each other.
\end{proof}

\begin{Prop}\label{PropDualSquareRootCoars}
We have 
\[
(S_{A,M}^{\G})^{\rtimes} \cong S_{A^{\rtimes},\C}^{\check{\G}} \boxtimes {}_{\C}L^2(\G)_{M^{\rtimes}},
\]
where ${}_{\C}L^2(\G)_{M^{\rtimes}}$ is a $\check{\G}$-W$^*$-Morita correspondence between $\C$ and $M^{\rtimes}$.  
\end{Prop}

\begin{proof}
By Theorem \ref{TheoCompTTDual}, Example \ref{ExaW*semiCoars}, Proposition \ref{PropIdRestr} and Proposition \ref{PropRestrBothBoth}, we have that 
\[
(S_{A,M}^{\G})^{\rtimes} \cong ({}_AL^2(A)_{\C}\boxtimes {}_{\C}L^2(\G)_M)^{\rtimes}\cong {}_{A^{\rtimes}}L^2(A^{\rtimes})_{\check{M}} \boxtimes {}_{\check{M}} L^2(M^{\rtimes})_{M^{\rtimes}}.
\]
So by Example \ref{ExaRestrComp}, it is sufficient to show that 
\[
{}_{\check{M}} L^2(M^{\rtimes})_{M^{\rtimes}} \cong {}_{\check{M}}L^2(\G)_{\C}\boxtimes {}_{\C}L^2(\G)_{M^{\rtimes}}.  
\]
But we have in fact that 
\[
{}_{M^{\rtimes}} L^2(M^{\rtimes})_{M^{\rtimes}} \cong {}_{\check{M}}L^2(\G)_{\C}\boxtimes {}_{\C}L^2(\G)_{M^{\rtimes}}.  
\]
Indeed, we may as well prove this with $M$ replaced by $\check{M}$, in which case we can write the isomorphism as 
\[
{}_{\C^{\rtimes\rtimes}}L^2(\C^{\rtimes\rtimes})_{\C^{\rtimes\rtimes}} \cong {}_{\C^{\rtimes\rtimes}}L^2(\check{\G})_{\C}\boxtimes {}_{\C}L^2(\check{\G})_{\C^{\rtimes\rtimes}},
\]
which is just a special case of the Takesaki-Takai duality of Theorem \ref{TheoTakTakTakTak}. 
\end{proof}

\begin{Cor}\label{CorDualCoars}
We have
\begin{equation}\label{EqDualSemiAdj2}
(C_{A,B}^{\G})^{\rtimes} \cong S_{A^{\rtimes},B^{\rtimes}}^{\check{\G}}
\end{equation}
\end{Cor}
\begin{proof}
Using the Proposition \ref{PropDualSquareRootCoars} and Theorem \ref{TheoCompTTDual}, we find 
\[
(C_{A,B}^{\G})^{\rtimes} \cong  (S_{A,M}^{\G}\boxtimes \overline{S_{B,M}^{\G}})^{\rtimes} \cong S_{A^{\rtimes},\C}^{\check{\G}} \boxtimes {}_{\C}L^2(\G)_{M^{\rtimes}} \boxtimes {}_{M^{\rtimes}}L^2(\G)_{\C} \boxtimes S_{\C,B^{\rtimes}}^{\check{\G}} \cong S_{A^{\rtimes},\C}^{\check{\G}} \boxtimes S_{\C,B^{\rtimes}}^{\check{\G}} \cong S_{A^{\rtimes},B^{\rtimes}}^{\check{\G}},
\]
where in the penultimate step we used that ${}_{\C}L^2(\G)_{M^{\rtimes}}$ is a $\check{\G}$-W$^*$-Morita correspondence between $\C$ and $M^{\rtimes}$ and Proposition \ref{PropMoritaGood}, and where in the last step we used Example \ref{ExaW*semiCoars}.
\end{proof}

We also have the following observation.

\begin{Prop}\label{PropDualSquareRootCoars2}
We have 
\[
({}_AL^2(A)_{\C} \boxtimes {}_{\C}L^2(\G)_{\C})^{\rtimes} \cong S_{A^{\rtimes},\check{M}}^{\check{\G}}, 
\]
where ${}_AL^2(A)_{\C} \boxtimes {}_{\C}L^2(\G)_{\C} = {}_A\Res_{\C}(S_{A,M}^{\G})$. 
\end{Prop}
\begin{proof}
This follows directly from Proposition \ref{PropDualSquareRootCoars} upon using Proposition \ref{PropRestrBothBoth}. 
\end{proof}

This leads to the following, maybe more surprising corollary.
\begin{Cor}\label{CorDualCoars2}
We have
\begin{equation}
\label{EqCoarseAmpl}(C_{A,B}^{\G}\otimes L^2(\G))^{\rtimes}  \cong C_{A^{\rtimes},B^{\rtimes}}^{\check{\G}}
\end{equation}
where the $\G$-$A$-$B$-correspondence on the left is $C_{A,B}^{\G}$, amplified with a multiplicity Hilbert space $L^2(\G)$ (with which the structures have no interaction, indicated by the fact that we are using $\otimes$ and not $\boxtimes$). 
\end{Cor} 
\begin{proof}
By Fell absorption, we have 
\[
{}_{\C}L^2(\G)_{\C}\otimes L^2(\G) \cong {}_{\C}L^2(\G)_{\C}\boxtimes {}_{\C}L^2(\G)_{\C},
\]
where on the left $L^2(\G)$ is just a multiplicity Hilbert space, and where on the right we take the Connes fusion product of equivariant correspondences (which in this case is just the usual tensor product of the right regular $\G$-representation with itself). We thus obtain 
\[
C_{A,B}^{\G}\otimes L^2(\G) = ({}_AL^2(A)_{\C} \boxtimes {}_{\C}L^2(\G)_{\C})\boxtimes ({}_{\C}L^2(\G)_{\C}\boxtimes {}_{\C}L^2(B)_B),
\]
and the corollary now follows directly from Proposition \ref{PropDualSquareRootCoars2}. 
\end{proof}

As a guide to the reader, we offer the following table which summarizes some of the results of this section. 

\begin{center}
\begin{table}[ht]
\begin{NiceTabular}{|c | c | c|c|}[vlines]
 \hline \\[-0.2cm]
 $\Hsp$ & $\Hsp\boxtimes \overline{\Hsp}$ & $\Hsp^{\rtimes}$& $(\Hsp\boxtimes\overline{\Hsp})^{\rtimes} \cong \Hsp^{\rtimes}\boxtimes \overline{\Hsp^{\rtimes}}$ \\[0.2cm] 
 \hline\hline\\[-0.2cm]   ${}_{\alpha(A)}L^2(A\bar{\otimes}M)_{A\bar{\otimes}M}$& $S_A^{\G}$ & ${}_{\alpha^{\rtimes}(A^{\rtimes})}L^2(A^{\rtimes}\rtimes_{\alpha^{\rtimes}} \check{\G})_{A^{\rtimes}\rtimes_{\alpha^{\rtimes}} \check{\G}}$& $D_{A^{\rtimes}}^{\check{\G}}$
  \\[0.2cm]  \hline\\[-0.2cm] 
    ${}_{\alpha(A)}L^2(A\rtimes \G)_{A\rtimes M}$ & $D_A^{\G}$&  ${}_{\alpha^{\rtimes}(A^{\rtimes})}L^2(A^{\rtimes}\overline{\otimes} \check{M})_{A^{\rtimes}\overline{\otimes} \check{M}} $ & $S_{A^{\rtimes}}^{\check{\G}}$
   \\[0.2cm] 
 \hline\\[-0.2cm]
  $S_{A,\C}^{\G} = {}_AL^2(A)_{\C}$ & $S_{A,A}^{\G}$& ${}_{A^{\rtimes}}L^2(A^{\rtimes})_{\check{M}}$ & ${}_{A^{\rtimes}}L^2(A^{\rtimes})_{\check{M}} \boxtimes {}_{\check{M}}L^2(A^{\rtimes})_{A^{\rtimes}}$ 
  \\[0.2cm]
 \hline\\[-0.2cm] 
 $S_{A,M}^{\G}$ & $C_{A,A}^{\G}$ & $S_{A^{\rtimes},\C}^{\check{\G}} \boxtimes {}_{\C}L^2(\G)_{M^{\rtimes}}$ & $S_{A^{\rtimes},A^{\rtimes}}^{\check{\G}}$
 \\[0.2cm] 
 \hline\\[-0.2cm]
 ${}_AL^2(A)_{\C} \boxtimes {}_{\C}L^2(\G)_{\C}$ & $C_{A,A}^{\G}\otimes L^2(\G)$ & $S_{A^{\rtimes},\check{M}}^{\check{\G}}$ & $C_{A^{\rtimes},A^{\rtimes}}^{\check{\G}}$\\[0.2cm]\hline
\end{NiceTabular}
\caption{Operations for basic correspondences}
\label{TableCorrConn}
\end{table}
\end{center}

\section{Preservation of weak containment and left amenability}

\subsection{Fell continuity for operations on equivariant correspondences}

We assume by default that W$^*$-categories admit subobjects and arbitrary (small) direct sums. Then a functor $F: \mcC\rightarrow \mcD$ between W$^*$-categories is called \emph{normal} if all the induced maps $F_X: \End_{\mcC}(X) \rightarrow \End_{\mcD}(F(X))$ are normal for any object $X \in \mcC$. This is the same \cite{GLR85} as each $F_{X,Y}: \Hom_{\mcC}(X,Y) \rightarrow \Hom_{\mcD}(F(X),F(Y))$ being $\sigma$-weakly continuous for each pair of objects $X,Y \in \mcC$. Moreover, this is also the same as $F$ preserving arbitrary (small) direct sums, for a unital $*$-homomorphism between von Neumann algebras is normal if and only if it is completely additive (i.e.\ preserves arbitrary sums of orthogonal projections). 

Normality takes care of continuity on the level of morphisms. Assume now in particular that $\mcC  \subseteq \Rep_*(C)$ and $\mcD \subseteq \Rep_*(D)$ for $C,D$ C$^*$-algebras. Then on the class of objects we have the Fell topology, and we can ask when a normal functor $F: \mcC\rightarrow \mcD$ is continuous with respect to the Fell topology. This will be the case if and only if $F$ preserves weak containment: 
\[
\pi \preccurlyeq \theta \quad \Rightarrow \quad F(\pi) \preccurlyeq F(\theta). 
\]
Indeed, it is immediate that Fell continuity of $F$ implies that $F$ preserves weak containment. For the converse direction, one uses the well-known fact that if $\pi\in \Rep_*(C)$ and $(\pi_{\alpha})_{\alpha\in I}$ a net in $\Rep_*(C)$, then $\pi_{\alpha}\rightarrow \pi$ in the Fell topology if and only if $\pi \preccurlyeq \oplus_{\beta\in J}\pi_{\beta}$ for each subnet $(\pi_{\beta})_{\beta\in J}$ of $(\pi_{\alpha})_{\alpha\in I}$ (see the introduction of \cite{Fel63} for the explicit statement, and for example \cite{BdlHV08b}*{Proposition F.2.2} for the elements of the proof). 

In light of the above, we introduce the following terminology.
\begin{Def}
We call a $*$-functor $F: \Rep_*(C) \rightarrow \Rep_*(D)$ \emph{W$^*$-continuous} if it preserves arbitrary direct sums and weak containment. 
\end{Def}
We can extend this terminology to functors $F: \mcC \rightarrow \mcD$ of W$^*$-categories $\mcC \subseteq \Rep_*(C)$ and $\mcD \subseteq \Rep_*(D)$ which are closed under direct sums and subobjects. In particular, we can use this terminology for functors between categories of equivariant correspondences.  

We now show the W$^*$-continuity of the operations introduced in Section \ref{SecOpEquiCorr}. It is clear that the operations introduced there are compatible with arbitrary direct sums. We hence focus below on the preservation of weak containment.

\begin{Prop}\label{PropWeakContOpp}
Let $\Gsp,\Hsp$ be $\G$-$A$-$B$-correspondences. If $\Gsp \preccurlyeq \Hsp$, then $\overline{\Gsp} \preccurlyeq \overline{\Hsp}$.
\end{Prop} 
\begin{proof}
This follows immediately from the original definition, since for $a\in A,b\in B,\omega \in M_*$ and $\xi,\eta\in \Gsp$, we have
\[
\langle \overline{\xi},\overline{\pi}(b)\overline{U}(\omega)\overline{\rho}(a)\overline{\eta}\rangle = \langle \overline{\xi}, \overline{\rho(b^*)U(\omega\circ R)^*\pi(a^*) \eta}\rangle  = \langle \eta,\pi(a)U(\omega\circ R)\rho(b)\xi\rangle.
\]\end{proof}

\begin{Prop}\label{PropContCompJa}
Let $\Gsp,\Hsp$ be $\G$-$A$-$B$-correspondences, and  $\Ksp$ a $\G$-$B$-$C$-correspondence. Then 
\begin{equation}\label{EqLeftCont}
\Gsp \preccurlyeq \Hsp \quad \Rightarrow \quad \Gsp \boxtimes \Ksp \preccurlyeq \Hsp \boxtimes \Ksp,
\end{equation}
and this becomes $\Leftrightarrow$ if $\Ksp$ is a $\G$-W$^*$-Morita correspondence.
 
Similarly for taking the fusion product on the left. 
\end{Prop}
\begin{proof}
By Theorem \ref{TheoMain1}, we can choose a contractive net $z_{\alpha} \in X = X_B(\Gsp,\Hsp^{\infty})$ such that for all $a\in A$
\[
(z_{\alpha}^*\pi_{\Hsp^{\infty}}(a)\otimes 1)\alpha_X(z_{\alpha}) \rightarrow \pi_{\Gsp}(a)\otimes 1,\qquad \sigma\textrm{-weakly}.
\]
If then $x,y \in X_B(\Gsp)$ and $\xi,\eta\in \Ksp$, we have $z_{\alpha}x,z_{\alpha}y \in X_B(\Hsp^{\infty})$ and, for $a\in A,c\in C$ and $\omega = \omega_{\zeta,\zeta'}\in M_*$, we compute that
\begin{eqnarray*}
&& \hspace{-2cm}\langle (z_{\alpha}x\otimes_B \xi),\pi_{\boxtimes}(a)U_{\boxtimes}(\omega)\rho_{\boxtimes}(c)(z_{\alpha}y\otimes_B \eta)\rangle  \\  &=& \langle (z_{\alpha}x\otimes 1)\otimes_{B\bar{\otimes}M}(\xi \otimes \zeta),(\pi_{\Hsp^{\infty}}(a)\otimes 1)\alpha_{X_B(\Hsp^{\infty})}(z_{\alpha}y)\otimes_{B\bar{\otimes}M} U_{\Ksp}(\rho_{\Ksp}(c)\eta\otimes \zeta')\rangle\\
&=&  \langle \xi\otimes \zeta,(x^*z_{\alpha}^*\pi_{\Hsp^{\infty}}(a) \otimes 1)\alpha_X(z_{\alpha})\alpha_{X_B(\Gsp)}(y)U_{\Ksp}(\rho_{\Ksp}(c)\eta\otimes \zeta')\rangle,
\end{eqnarray*}
which converges to 
\[
\langle \xi\otimes \zeta,(x^*\pi_{\Gsp}(a)\otimes 1) \alpha_X(y)U_{\Ksp}(\rho_{\Ksp}(c)\eta\otimes \zeta')\rangle = \langle (x\otimes_B \xi),\pi_{\boxtimes}(a)U_{\boxtimes}(\omega)\rho_{\boxtimes}(c)(y\otimes_B \eta)\rangle
\]
by assumption. 

Together with $\Hsp^{\infty} \boxtimes \Ksp \cong (\Hsp\boxtimes \Ksp)^{\infty}$, this immediately implies \eqref{EqLeftCont}. If moreover $\Ksp$ is a $\G$-W$^*$-Morita equivalence, we use the above result in combination with 
\[
(\Gsp\boxtimes \Ksp)\boxtimes \overline{\Ksp} \cong \Gsp\boxtimes (\Ksp\boxtimes \overline{\Ksp}) \cong \Gsp \boxtimes E_C^{\G} \cong \Gsp,
\]
using associativity of $\boxtimes$ and Proposition \ref{PropMoritaGood}, and similarly for $\Hsp$.

For composition on the left, we can use Proposition \ref{PropWeakContOpp} and Corollary \ref{CorCompOpp}.
\end{proof}

We now consider the crossed product construction. 

\begin{Theorem}\label{TheoAsMult}
If $\Hsp,\Gsp$ are $\G$-$A$-$B$-correspondences, then 
\[
\Gsp \preccurlyeq \Hsp \quad \Leftrightarrow \quad \Gsp^{\rtimes}\preccurlyeq \Hsp^{\rtimes}.
\] 
\end{Theorem}
\begin{proof}
By Theorem \ref{TheoTakTakEq} and Proposition \ref{PropContCompJa}, we only need to prove that $\Gsp \preccurlyeq \Hsp$ implies $\Gsp^{\rtimes}\preccurlyeq \Hsp^{\rtimes}$. By replacing $\Hsp$ with $\Hsp^{\infty}$, we may assume that $\Hsp^{\infty} \cong \Hsp$.

By Theorem \ref{TheoMain1}, we may then pick a contractive net $y_i \in X = X_B(\Gsp,\Hsp)$ such that 
\[
(y_i^*\pi_{\Hsp}(a)\otimes 1)\alpha_{X}(y_i) \overset{\sigma\textrm{-weakly}}{\underset{i\rightarrow \infty}{\longrightarrow}} \pi_{\Gsp}(a)\otimes 1,\qquad a\in A.
\]
Put $z_i = \alpha_X(y_i)$. Then $z_i \in X_{B^{\rtimes}}(\Gsp^{\rtimes},\Hsp^{\rtimes})$ is a contractive net, and for all $a\in A$ one has
\begin{equation}\label{EqConvStart}
(z_i^*(\pi_{\Hsp}\otimes \pi_M)\alpha(a))_{12}V_{23}z_{i,12} \overset{\sigma\textrm{-weakly}}{\underset{i\rightarrow \infty}{\longrightarrow}}
 (\pi_{\Gsp}\otimes \pi_M)(\alpha(a))_{12}V_{23}. 
\end{equation}
In particular, we see that for all $a\in A$ and all $x\in C_0^{\red}(\check{\G})$ one has
\begin{equation}\label{EqConvStart2}
z_i^*(\pi_{\Hsp}\otimes \pi_M)(\alpha(a))(1\otimes \pi_{\check{M}}(x))z_i \overset{\sigma\textrm{-weakly}}{\underset{i\rightarrow \infty}{\longrightarrow}}
(\pi_{\Gsp}\otimes \pi_M)(\alpha(a))(1\otimes \pi_{\check{M}}(x)).
\end{equation}
By passing to a subnet, we may assume that $\Ad_{z_i^*}$ pointwise $\sigma$-weakly converges to a ccp map
\[
E: B(\Hsp\otimes L^2(\G))\rightarrow B(\Gsp\otimes L^2(\G)).
\]
Then $1\otimes \pi_{\check{M}}(C_0^{\red}(\check{\G}))$ lies in the multiplicative domain of $E$. In particular, since $C_0^{\red}(\G)$ acts non-degenerately on $L^2(\G)$, it follows that $E$ is ucp. It is also immediate from their definition that the $z_i$ intertwine $U_{\Gsp}(1\otimes x)U_{\Gsp}^*$ and $U_{\Hsp}(1\otimes x)U_{\Hsp}^*$ for $x\in B(L^2(\G))$. In particular, $U_{\Hsp}(1\otimes \mcK(L^2(\G)))U_{\Hsp}^*$ lies in the multiplicative domain of $E$. 

Now 
\begin{eqnarray*}
&& \hspace{-2cm}[U_{\Hsp}(1\otimes \mcK(L^2(\G)))U_{\Hsp}^*(1\otimes C_0^{\red}(\check{\G}))U_{\Hsp}(1\otimes \mcK(L^2(\G)))U_{\Hsp}^*] \\&=& 
[U_{\Hsp}(1\otimes x)(\id\otimes \id\otimes \omega)(U_{\Hsp,12}^*V_{23}U_{\Hsp,12}) (1\otimes y)U_{\Hsp}^*\mid x,y\in \mcK(L^2(\G)),\omega\in B(L^2(\G))_*]\\
&=& [U_{\Hsp}(1\otimes x)(\id\otimes \id\otimes \omega)(U_{\Hsp,13}V_{23}) (1\otimes y)U_{\Hsp}^*\mid x,y\in \mcK(L^2(\G)),\omega\in B(L^2(\G))_*]\\
&=& [U_{\Hsp}((\id\otimes \omega)(U_{\Hsp}(1\otimes (\chi\otimes \id)(V)))\otimes x)U_{\Hsp}^*\mid x\in \mcK(L^2(\G)),\omega,\chi\in B(L^2(\G))_*]\\
&=&[U_{\Hsp}(U_{\Hsp}(\omega)\otimes x)U_{\Hsp}^*\mid x\in \mcK(L^2(\G)),\omega\in B(L^2(\G))_*]\\
&=& U_{\Hsp}(C_{U_{\Hsp}}\otimes \mcK(L^2(\G)))U_{\Hsp}^*
\\ &=& C_{U_{\Hsp}}\otimes \mcK(L^2(\G)),
\end{eqnarray*}
where in the last step we used Proposition \ref{PropImRepC}. 

Letting $C^{\G}(A,B) = C^{G}_{\Hsp_u}(A,B)$ with respect to a universal model $(\Hsp_u,\pi_u,\rho_u,U_u)$, we find from the above computation that for all $x\in C_{U_u}\otimes \mcK(L^2(\G))$ we have
\begin{equation}\label{EqMultDomArgAgainAgainAgainAgain}
E((\theta_{\Hsp}\otimes \id)(x)) =  (\theta_{\Gsp}\otimes \id)(x).
\end{equation}
Since also $(\pi_{\Hsp}\otimes \pi_M)(\alpha(A))$ lies in the multiplicative domain of $E$ by \eqref{EqConvStart2}, and 
\[
[(w\otimes x)(\pi_{\Hsp_u}\otimes \pi_M)(\alpha(a))(z\otimes y)\mid a\in A,w,z\in C_{U_u},x,y\in \mcK(L^2(\G))] = \pi_{\Hsp_u}(C^{\G}_{\Hsp_u}(A)) \otimes \mcK(L^2(\G)),
\] 
we find that \eqref{EqMultDomArgAgainAgainAgainAgain} in fact holds for all  $x\in C^{\G}_{\Hsp_u}(A)\otimes \mcK(L^2(\G))$. 

But $(\pi_{\Hsp}\otimes \id)(A\rtimes \G) \subseteq QM( \theta_{\Hsp}(C^{\G}_{\Hsp_u}(A))\otimes \mcK(L^2(\G)))$, so we can finally conclude that for each $x\in A\rtimes \G$ we have
\begin{equation}\label{EqMultDomArgAgainAgainAgain}
z_i^*(\pi_{\Hsp}\otimes \id)(x)z_i\overset{\sigma\textrm{-weakly}}{\underset{i\rightarrow \infty}{\longrightarrow}} E((\pi_{\Hsp}\otimes \id)(x)) =  (\pi_{\Gsp}\otimes \id)(x).
\end{equation}
Since $\alpha_{X^{\rtimes}}(z_i) =z_i\otimes 1$, we can conclude from Theorem \ref{TheoMain1} that $\Gsp^\rtimes \preccurlyeq \Hsp^{\rtimes}$.
\end{proof}

\subsection{Left amenability revisited}

Recall the definition of strong left equivariant amenability for a $\G$-equivariant correspondence from Definition \ref{DefEquivariantlyStrong}. By Remark \ref{RemSauv}, we obtain the following equivalent characterisation. 
\begin{Prop}\label{PropEquivCharStrongLeftAm}
Let $A,B$ be $\G$-W$^*$-algebras. A $\G$-equivariant $A$-$B$-correspondence $\Hsp$ is  strongly left equivariantly amenable if and only if $\Hsp\boxtimes \overline{\Hsp}$ equivariantly weakly contains $E_A^{\G}$. 
\end{Prop}
For example, we can characterize strong equivariant amenability of an action of $\G$ on a von Neumann algebra $A$ by the property that ${}_{\alpha(A)}(L^2(A)\otimes L^2(\G))_{A\bar{\otimes}M}$ is strongly left $\G$-amenable, or by the fact that the $\G$-semi-coarse correspondence equivariantly weakly contains the trivial correspondence. 

In the remainder of this section, we want to explore some other notions of `equivariant amenability' and their interconnections. For completeness and ease of reference, we repeat the definition of (strong) $\G$-amenability within the next definition.

\begin{Def}\label{DefAmena}\label{DefPropApproxProps} Let $\G$ be a locally compact quantum group, and let $(A,\alpha)$ be a $\G$-dynamical von Neumann algebra. 

We say that the action of $\G$ on $A$ is
\begin{itemize}
\item \emph{(strongly) $\G$-amenable} if ${}_{\alpha(A)}(L^2(A)\otimes L^2(\G))_{A\bar{\otimes}M}$ is (strongly) left equivariantly amenable, where $A\bar{\otimes} M$ has the $\G$-action determined by $\id\otimes \Delta$,
\item \emph{(strongly) W$^*$-amenable} if $S_{A,\C}^{\G}$ is (strongly) left equivariantly amenable,
\item \emph{(strongly) $\G$-W$^*$-amenable} if $S_{A,M}^{\G}$ is (strongly) left equivariantly amenable,
\item \emph{(strongly) inner $\G$-amenable} if ${}_{\alpha(A)}(L^2(A)\otimes L^2(\G))_{A\rtimes \G}$ is (strongly) left equivariantly amenable, where $A\rtimes \G$ has the adjoint $\G$-action \eqref{EqAdjAction}.
\end{itemize}
\end{Def}

So, making use of the second column of Table \ref{TableCorrConn}, we see that $(A,\alpha)$ is
\begin{itemize}
\item strongly $\G$-amenable iff $E_A^{\G}$ is weakly contained in the $\G$-semi-coarse $\G$-$A$-$A$-correspondence $S_A^{\G}$,
\item strongly W$^*$-amenable iff $E_A^{\G}$ is weakly contained in the W$^*$-semi-coarse $\G$-$A$-$A$-correspondence $S_{A,A}^{\G}$,
\item strongly $\G$-W$^*$-amenable iff $E_A^{\G}$ is weakly contained in the coarse $\G$-$A$-$A$-correspondence $C_{A,A}^{\G}$, and
\item strongly inner $\G$-amenable iff $E_A^{\G}$ is weakly contained in the adjoint $\G$-$A$-$A$-correspondence $D_A^{\G}$.
\end{itemize}

On the other hand, referring back to Definition \ref{DefEquivariantlyStrong}, we see immediately from Definition \ref{DefAmena} that $(A,\alpha)$ is
\begin{itemize}
\item $\G$-amenable  iff there exists a $\G$-equivariant conditional expectation $A\bar{\otimes}M \rightarrow \alpha(A)$.
\item W$^*$-amenable iff there exists a $\G$-equivariant conditional expectation $B(L^2(A)) \rightarrow A$. Here $B(L^2(A))$ is endowed with the $\G$-action $x \mapsto U_{\alpha}(x\otimes 1)U_{\alpha}^*$. 
\item $\G$-W$^*$-amenable iff there exists a $\G$-equivariant conditional expectation 
\[
B(L^2(A))\bar{\otimes} M \rightarrow A\otimes \C,
\]
where the left hand side is endowed with the $\G$-action
\[
B(L^2(A))\bar{\otimes} M \rightarrow (B(L^2(A))\bar{\otimes} M)\bar{\otimes} M,\qquad x \mapsto U_{\alpha,13}V_{23}x_{12}V_{23}^*U_{\alpha,13}^*.
\]
\item inner $\G$-amenable iff there exists a $\G$-equivariant conditional expectation $E: A \rtimes \G\rightarrow \alpha(A)$, where $A\rtimes \G$ is endowed with the adjoint $\G$-action. 
\end{itemize}

Recall that from Theorem \ref{TheoAmenacorr}, it follows that the latter amenability properties follow from their respective strong version. 

For completeness, we include also the following obvious definition.

\begin{Def}\label{DefPropT} Let $\G$ be a locally compact quantum group, and let $(A,\alpha)$ be a $\G$-W$^*$-algebra. 

We say that $(A,\alpha)$ has \emph{$\G$-equivariant property $(T)$} if whenever a $\G$-$A$-$A$-correspondence $\Hsp$ contains $E_A$ weakly, then $E_A \subseteq \Hsp$ as a $\G$-equivariant subcorrespondence. 
\end{Def}

Note the following: 
\begin{itemize}
\item If $\G= \{e\}$ is trivial, then strong $\G$-W$^*$-amenability of $A$ is the same as strong $W^*$-amenability of $A$, and simply means that the trivial $A$-$A$-correspondence $L^2(A)$ is weakly contained in the coarse correspondence $L^2(A)\otimes L^2(A)$, i.e.\ $A$ is amenable (=semi-discrete) as a von Neumann algebra \cite{EL77}. 
\item If $A= \C$, we are considering the action of $\G$ on a point $\bullet$, and strong $\G$-amenability of the action is the same as strong $\G$-W$^*$-amenability. Moreover, the action $\bullet\curvearrowleft \G$ is strongly $\G$-amenable if and only if $\G$ is strongly amenable as a locally compact quantum group, i.e.\ the trivial representation is weakly contained in the regular representation (this is also known as \emph{co-amenability of $\hat{\G}$}), see \cites{Rua96,DQV02,Tom06}.
\item The trivial action $\bullet\curvearrowleft\G$ is strongly inner $\G$-amenable if and only if $\G$ is strongly inner amenable (in the sense of \emph{locally compact quantum groups} \cite{OOT17}). Similarly, the trivial action $\bullet\curvearrowleft\G$ is inner $\G$-amenable if $\G$ is inner amenable (in the sense of \cite{Cr19}).
\item If $A= \C$, then $\G$-equivariant property $(T)$ is the same as the usual property $(T)$ for locally compact quantum groups \cite{CN15} (see also \cites{Fim10,KyS12}).
\end{itemize}

\begin{Rem}
It is less clear to us how to introduce the \emph{Haagerup property} in the setting of equivariant correspondences. It should combine the notions in \cite{DFSW16} and \cite{OOT17}. 
\end{Rem}

We collect now some theorems illustrating the connections between these various approximation properties.

\begin{Theorem}\label{TheoGWStar}
Let $(A,\alpha)$ be a $\G$-W$^*$-algebra.
\begin{enumerate}
\item If $(A,\alpha)$ is strongly W$^*$-amenable and strongly $\G$-amenable, then $(A,\alpha)$ is strongly $\G$-$W^*$-amenable.
\item $(A,\alpha)$ is W$^*$-amenable and $\G$-amenable if and only if $(A,\alpha)$ is $\G$-$W^*$-amenable.
\end{enumerate}
\end{Theorem}
\begin{proof}
Under the assumptions in the first item, we have that $E_A^{\G}$ is weakly contained in $S_A^{\G}$ and in $S^{\G}_{A,A}$. By Example \ref{EqUnitUnderTens} and Example \ref{ExaCompGivesCoars}, together with  Proposition \ref{PropContCompJa}, we see that $E_A^{\G}$ is weakly contained in $C_{A,A}^{\G}$, i.e.\ $(A,\alpha)$ is strongly $\G$-W$^*$-amenable.

Let us now prove the second item. Assume first that $(A, \alpha)$ is both $W^*$-amenable and $\G$-amenable. This means that there exists a $\G$-equivariant ucp map $E: (A \bar{\otimes} L^\infty(\G), \id \otimes \Delta)\to (A, \alpha)$ such that $E\circ \alpha = \id_A$ and a $\G$-equivariant ucp map $\varphi: B(L^2(A))\to A$ such that $\varphi\circ \pi_A = \id_A$. Define then the ucp map
$$\psi: B(L^2(A))\bar{\otimes}L^\infty(\G) \to A: z \mapsto E((\varphi\otimes \id)(U_\alpha z U_\alpha^*)).$$
If $a\in A$, we have
$$\psi(\pi_A(a)\otimes 1)= E((\varphi\otimes \id)(\pi_A \otimes \id)(\alpha(a)))= E(\alpha(a)) = a$$
and if $z\in B(L^2(A))\bar{\otimes} L^\infty(\G)$, we have
\begin{align*}
    (\psi\otimes \id)(U_{\alpha,13}V_{23}z_{12}V_{23}^*U_{\alpha,13}^*)&= (E \otimes \id)(\varphi\otimes \id \otimes \id)(U_{\alpha,12}U_{\alpha,13}V_{23}z_{12}V_{23}^*U_{\alpha,13}^*U_{\alpha,12}^*)\\
    &= (E \otimes \id)(\varphi \otimes \id \otimes \id)(\id \otimes \Delta)(U_\alpha z U_\alpha^*)\\
    &= (E\otimes \id)(\id \otimes \Delta)(\varphi \otimes \id)(U_\alpha z U_\alpha^*)\\
    &= \alpha(E(\varphi \otimes \id)(U_\alpha z U_\alpha^*)) = \alpha(\psi(z))
\end{align*}
so $\psi$ is a $\G$-equivariant ucp map that witnesses the $\G$-$W^*$-amenability of $(A, \alpha)$. 

Conversely, assume that $(A, \alpha)$ is $\G$-$W^*$-amenable, i.e.\ there exists a ucp $\G$-equivariant map
\[
\psi: B(L^2(A))\bar{\otimes} L^\infty(\G)\to A
\]
that satisfies $\psi(\pi_A(a)\otimes 1)= a$ for all $a\in A$. Define 
$$E: A \bar{\otimes} L^\infty(\G)\to A: z \mapsto \psi(U_\alpha^*(\pi_A\otimes \id)(z)U_\alpha).$$
Then if $a\in A$, we have
$$E(\alpha(a)) = \psi(\pi_A(a)\otimes 1) = a$$
and if $z\in A \bar{\otimes} L^\infty(\G)$, we have
\begin{align*}
    (E \otimes \id)(\id \otimes \Delta)(z) &= (E \otimes \id)(U_{\alpha, 12}^*(\pi_A \otimes \id \otimes \id)(\id \otimes \Delta)(z)U_{\alpha,12})\\
    &= (\psi \otimes \id)((U_{\alpha, 12}^*(\id \otimes \Delta)(\pi_A\otimes \id)(z)U_{\alpha,12})\\
    &= (\psi \otimes \id)(U_{\alpha,13} (\id \otimes \Delta)(U_\alpha^*(\pi_A\otimes \id)(z)U_\alpha)U_{\alpha,13}^*)\\
    &= \alpha(\psi(U_\alpha^*(\pi_A\otimes \id)(z)U_\alpha))= \alpha(E(z)).
\end{align*}
Therefore, the map $E$ witnesses the $\G$-amenability of $(A, \alpha)$. On the other hand, define the ucp map
$$\varphi: B(L^2(A))\to A: x \mapsto \psi(x\otimes 1).$$
Then clearly $\varphi\circ \pi_A= \id_A$ and if $x\in B(L^2(A))$, we have
\begin{align*}
    (\varphi\otimes \id)(U_\alpha(x\otimes 1)U_\alpha^*)&= (\psi \otimes \id)(U_{\alpha,13}(x\otimes 1 \otimes 1)U_{\alpha,13}^*)\\
    &= (\psi \otimes \id)(U_{\alpha,13}V_{23}(x\otimes 1 \otimes 1)V_{23}^*U_{\alpha,13}^*)\\
    &= \alpha(\psi(x\otimes 1)) = \alpha(\varphi(x))
\end{align*}
so the map $\varphi$ witnesses the $W^*$-amenability of $(A, \alpha)$. 
\end{proof}

\begin{Rem}
We do not know if, in general, strong $\G$-W$^*$-amenability implies strong $\G$-amenability and strong W$^*$-amenability. See however statements $(3)$ and $(4)$ in Theorem \ref{TheoDualApprox}.
\end{Rem}

In what follows, we fix a locally compact quantum group $\G$ and we will use the notations
\begin{align*}
    &\Delta_r: B(L^2(\G))\to B(L^2(\G))\bar{\otimes} M: x \mapsto V(x\otimes 1)V^*,\\
    &\check{\Delta}_r: B(L^2(\G))\to B(L^2(\G))\bar{\otimes} \check{M}: x \mapsto \check{V}(x\otimes 1)\check{V}^*,\\
    &\Delta_l: B(L^2(\G))\to M \bar{\otimes} B(L^2(\G)): x \mapsto W^*(1\otimes x)W.
\end{align*}

We also recall that if $(A, \alpha)$ is a $\G$-von Neumann algebra, then the $\check{\G}$-fixed points of $A\rtimes_\alpha \G$ with respect to the $\check{\G}$-coaction $\alpha^\rtimes = \id \otimes \check{\Delta}_r$ is exactly $\alpha(A)$ and we have
$$A^\rtimes = A\rtimes_\alpha \G = \{z\in A \bar{\otimes} B(L^2(\G)): (\alpha\otimes \id)(z) = (\id \otimes \Delta_l)(z)\}.$$

\begin{Theorem}\label{TheoDualWeakContCross}
Let $\G$ be a locally compact quantum group, and let $A,B$ be $\G$-dynamical von Neumann algebras. Let $\Hsp = (\Hsp, U, \pi, \rho)$ be a $\G$-equivariant $A$-$B$-correspondence. 
\begin{enumerate}
    \item $\Hsp$ is strongly left $\G$-amenable if and only if $\Hsp^{\rtimes}$ is strongly left $\check{\G}$-amenable.
    \item $\Hsp$ is left $\G$-amenable if and only if $\Hsp^{\rtimes}$ is left $\check{\G}$-amenable.
\end{enumerate} 
\end{Theorem}
\begin{proof}
We first prove the first item. By Proposition \ref{PropOpposit} and Theorem \ref{TheoCompTTDual}, we have that
\[
(\Hsp \boxtimes \overline{\Hsp})^{\rtimes} \cong \Hsp^{\rtimes}\boxtimes \overline{\Hsp^{\rtimes}}. 
\]
So by Proposition \ref{PropIdRestr} and Theorem \ref{TheoAsMult}, we see that $E_A^{\G}$ is equivariantly weakly contained in $\Hsp\boxtimes \overline{\Hsp}$ if and only if $E_{A^{\rtimes}}^{\check{\G}}$ is equivariantly weakly contained in $\Hsp^{\rtimes}\boxtimes \overline{\Hsp^{\rtimes}}$, which is exactly the statement in the first item.

Next, we prove the second item. Note that the equivalence in $(2)$ means the equivalence of the following two items:
\begin{enumerate}
    \item[(a)] There exists a $\G$-equivariant ucp conditional expectation $\rho(B)'\to \pi(A).$
    \item[(b)] There exists a $\check{\G}$-equivariant ucp conditional expectation $\rho^{\rtimes}(B^\rtimes)' \to \pi^{\rtimes}(A^\rtimes).$
\end{enumerate}

Before proving the equivalence of these statements, we need to do some preparatory work. Let us write $\widetilde{R}$ for the anti-unitary on $B(L^2(\G))$ given by $\widetilde{R}(x) = \check{J}x^*\check{J}$.  Given a subset $S$ of a von Neumann algebra, we write $[S]^\sigma$ for the $\sigma$-weak closure of the linear span of $S$. Then note that
\begin{align*}
    \rho^{\rtimes}(B^{\rtimes}) &= U(\rho \otimes \widetilde{R})([\beta(B)(1\otimes \check{M})]^\sigma)U^*\\
    &= [U(1\otimes \hat{M})U^*U(\rho\otimes R)(\beta(B))U^*]^\sigma= [U(1\otimes \hat{M})U^*(\rho(B)\otimes 1)]^\sigma
\end{align*}
whence 
$\rho^{\rtimes}(B^{\rtimes})'= [U(B(\Hsp)\bar{\otimes} \check{M})U^*] \cap [\rho(B)'\bar{\otimes}B(L^2(\G))].$
Next, note that if $z\in \rho^{\rtimes}(B^\rtimes)'$, then \begin{equation}\label{hihi}
    W_{23}^* z_{13}W_{23}= U_{12}z_{13}U_{12}^*.
\end{equation}
Indeed, if we put $X:= U^*zU$, then \eqref{hihi} is equivalent with \begin{equation}\label{hihihi} 
W_{23}^* U_{13}X_{13}U_{13}^*W_{23}= U_{12}U_{13}X_{13}U_{13}^*U_{12}^*.\end{equation} But since $X \in B(\Hsp)\bar{\otimes} \check{M}$ and $W\in M \bar{\otimes}\hat{M}$, we find
\begin{align*}
    U_{12}U_{13}X_{13}U_{13}^*U_{12}^* = (\id \otimes \Delta)(U)X_{13}(\id \otimes \Delta)(U^*)= W_{23}^*U_{13}W_{23}X_{13}W_{23}^*U_{13}^* W_{23} = W_{23}^* U_{13}X_{13}U_{13}^* W_{23}
\end{align*}
and \eqref{hihihi} follows. 

Assume now that $(a)$ holds, and choose a $\G$-equivariant ucp conditional expectation $\phi: \rho(B)'\to \pi(A)$. Define the map
$$\psi:= \phi \otimes \id: \rho(B)' \bar{\otimes} B(L^2(\G))\to \pi(A)\bar{\otimes} B(L^2(\G)).$$
Let us write $\gamma: B(\Hsp)\to B(\Hsp)\bar{\otimes} M: x \mapsto U(x\otimes 1)U^*$. 
If $z\in \rho^{\rtimes}(B^{\rtimes})'\subseteq \rho(B)'\bar{\otimes} B(L^2(\G))$, then using \eqref{hihi} and $\G$-equivariance of $\phi$, we get
\begin{align*}
    (\id_{\pi(A)}\otimes \Delta_l)\psi(z)&= (\phi \otimes \id \otimes \id)(\id \otimes \Delta_l)(z) = (\phi \otimes \id \otimes \id)(U_{12}z_{13}U_{12}^*)= (\gamma \otimes \id)\psi(z)
\end{align*}
so $\psi(z) \in \pi(A)\rtimes_\gamma \G= \pi^{\rtimes}(A^{\rtimes})$. Thus, $\psi= \phi\otimes \id$ restricts to a map
$$\psi: \rho^\rtimes(B^\rtimes)' \to \pi^{\rtimes}(A^{\rtimes})$$ which trivially is a ucp $\check{\G}$-equivariant conditional expectation. Thus, $(b)$ holds.

Conversely, assume that $(b)$ holds, and choose a ucp $\check{\G}$-equivariant ucp conditional expectation $\psi: \rho^{\rtimes}(B^{\rtimes})' \to \pi^{\rtimes}(A^{\rtimes})$. Since 
$U(\rho(B)'\otimes 1)U^*\subseteq \rho^{\rtimes}(B^{\rtimes})'$, it makes sense to define
$$\widetilde{\psi}: \rho(B)'\to \pi^{\rtimes}(A^{\rtimes}): x \mapsto \psi(U(x\otimes 1)U^*).$$
Note that since $\check{V}_{23}$ and $U_{12}$ commute, we have for $x\in \rho(B)'$ that
\begin{align*}
    (\id \otimes \check{\Delta}_r)\widetilde{\psi}(x)&= (\id \otimes \check{\Delta}_r)\psi(U(x\otimes 1)U^*)\\
    &= (\psi\otimes \id)(\id \otimes \check{\Delta}_r)(U(x\otimes 1)U^*)\\
    &= (\psi \otimes \id)(\check{V}_{23}U_{12}(x\otimes 1\otimes 1)U_{12}^*\check{V}_{23}^*)\\
    &= \psi(U(x\otimes 1)U^*)\otimes 1 = \widetilde{\psi}(x)\otimes 1
\end{align*}
so that $\widetilde{\psi}(x)\in \gamma(\pi(A))$. Next, note that since $\psi(1\otimes y)= 1\otimes y$ for all $y\in \check{M}$, we have that $V_{23}\in M(1\otimes C_0^{\operatorname{red}}(\check{\G})\otimes C_0^{\operatorname{red}}(\G))$ is in the multiplicative domain of $\psi \otimes \id$, so that for $x\in \rho(B)'$,
\begin{align*}
    (\widetilde{\psi}\otimes \id)(U(x\otimes 1)U^*)&= (\psi\otimes \id)(U_{12}U_{13}(x\otimes 1 \otimes 1)U_{13}^* U_{12}^*)\\
    &= (\psi\otimes \id)(\id \otimes \Delta)(U(x\otimes 1)U^*)\\
    &= V_{23}((\psi\otimes \id)(U(x\otimes 1)U^*)\otimes 1)V_{23}^* = (\id \otimes \Delta)\widetilde{\psi}(x).
\end{align*}
Thus, the map $\widetilde{\psi}: \rho(B)'\to \gamma(\pi(A))$ is $\G$-equivariant. Composing with the $\G$-isomorphism $$\gamma^{-1}: (\gamma(\pi(A)), \id \otimes \Delta)\cong (\pi(A), \gamma),$$ we find that $\phi:= \gamma^{-1} \circ \widetilde{\psi}: \rho(B)'\to \pi(A)$ is $\G$-equivariant, and it satisfies 
$$\phi(\pi(a)) = \gamma^{-1}(\widetilde{\psi}(\pi(a)))= \gamma^{-1}(\psi((\pi\otimes \id)(\alpha(a)))) = \gamma^{-1}((\pi\otimes \id)(\alpha(a))) =\gamma^{-1}(\gamma(\pi(a))) = \pi(a)$$ for every $a\in A$. Thus, $(a)$ is satisfied.
\end{proof}

\begin{Theorem}\label{TheoDualApprox}
Let $\G$ be a locally compact quantum group, and $A$ a $\G$-dynamical von Neumann algebra. Then the following equivalences hold: 
\begin{enumerate}
\item $A$ is (strongly) $\G$-amenable if and only if $A\rtimes \G$ is (strongly) inner $\check{\G}$-amenable. 
\item $A$ is (strongly) inner $\G$-amenable if and only if $A\rtimes \G$ is (strongly) $\check{\G}$-amenable. 
\item $A$ is (strongly) $\G$-W$^*$-amenable if and only if $A\rtimes \G$ is (strongly) $W^*$-amenable (with respect to the dual $\check{\G}$-action).
\item $A$ is strongly $\G$-W$^*$-amenable if and only if $A\rtimes \G$ is strongly $\check{\G}$-W$^*$-amenable.
\end{enumerate}
\end{Theorem}
\begin{proof}
This follows directly from Table \ref{TableCorrConn} upon using Theorem \ref{TheoDualWeakContCross}. The only things that need a small comment, are
\begin{itemize}
\item that, for the last equivalence, the multiplicity factor $L^2(\G)$ does not matter, see Remark \ref{RemMulti}, and
\item that, for the third equivalence, left $\check{\G}$-amenability of $S_{A^{\rtimes},\C}^{\check{\G}} \boxtimes {}_{\C}L^2(\G)_{M^{\rtimes}}$, whose underlying Hilbert space is $L^2(A^{\rtimes})\otimes L^2(\check{\G})$, is the same as left $\check{\G}$-amenability of $S_{A^{\rtimes},\C}^{\check{\G}}$, since 
\[
(1\otimes M^{\rtimes})' = B(L^2(A))\otimes 1. 
\]
\end{itemize}\end{proof}

We complement the non-strong version of the above theorem with the following remark:

\begin{Rem}
    \begin{itemize}
     \item The last equivalence in Theorem \ref{TheoDualApprox} is also true in the non-strong setting, although we do not know a quick proof using the machinery of this paper. The proof will be given in follow-up work by one of the authors.
    \item Theorem \ref{TheoDualApprox} unifies and generalizes several results that are already spread throughout the literature. For instance, the equivalence $(1)$ in the above theorem was shown in \cite{Moa18}*{Theorem 7.5} for $\G$ a discrete quantum group, and was also proven in the setting of $\G$-operator systems in \cite{DH24}*{Theorem 4.2} where $\G$ is a discrete quantum group. A similar result also appeared in \cite{BC21}*{Theorem 5.2, Corollary 5.3} in the setting of locally compact groups (and where the authors note that the techniques extend to co-amenable locally compact quantum groups). 
    \item In \cite{MP22}, a notion of inner amenability for actions of locally compact groups on von Neumann algebras is introduced. However, the definition given in this paper is not the same one as in Definition \ref{DefAmena}. That both definitions nevertheless agree (in the setting of locally compact groups) follows by combining Theorem \ref{TheoDualApprox} (2) with \cite{MP22}*{Proposition 4.3}.
\end{itemize}
\end{Rem}

We end this paper by noting that it is also possible to prove the first three equivalences in the non-strong setting in Theorem \ref{TheoDualApprox} directly. This is interesting, as it gives us an explicit bridge to move between certain conditional expectations. For completeness, we briefly describe for each of the equivalences how to do this. We refrain from spelling out all details.
    \begin{enumerate}
        \item $(A, \alpha)$ being $\G$-amenable is equivalent to the existence of a $\G$-equivariant ucp conditional expectation
        $$\varphi: (A\bar{\otimes}L^\infty(\G), \id \otimes \Delta) \to (\alpha(A), \id \otimes \Delta).$$
        Then the map
        $$\psi:= \varphi\rtimes \G: A \bar{\otimes} B(L^2(\G))\cong (A\bar{\otimes}L^\infty(\G))\rtimes \G \to \alpha(A)\rtimes \G\cong A\rtimes \G$$
        is a $\check{\G}$-equivariant ucp conditional expectation. Conversely, if $\psi: A \bar{\otimes} B(L^2(\G)) \to A\rtimes \G$ is a $\check{\G}$-equivariant ucp conditional expectation, then it restricts to a $\G$-equivariant\footnote{Making use of the fact that $V\in M(C_0^{\operatorname{red}}(\check{\G})\otimes C_0^{\operatorname{red}}(\G))$, we see that $V_{23}$ is contained in the multiplicative domain of $\psi\otimes \id$ with $(\psi\otimes \id)(V_{23}) = V_{23}$. Hence, $\psi$ preserves the $\G$-action $\id \otimes \Delta_r$ and thus its restriction $\varphi$ preserves the $\G$-action $\id \otimes \Delta$.} conditional expectation
        $$\varphi: A\bar{\otimes}L^\infty(\G)\to \alpha(A)$$
        between the $\check{\G}$-fixed points. On the other hand, the Takesaki-Takai isomorphism $$\Phi: (A\rtimes \G)\rtimes \check{\G}\to A\bar{\otimes} B(L^2(\G))$$ is $\check{\G}$-equivariant (where the iterated crossed product has the adjoint $\check{\G}$-action - this follows by an explicit calculation on generators) with $\Phi^{-1}\vert_{A\rtimes \G}= \id\otimes \check{\Delta}_r$, so the existence of a $\check{\G}$-equivariant ucp conditional expectation $\psi: A \bar{\otimes} B(L^2(\G)) \to A\rtimes \G$ is equivalent with the existence of a $\check{\G}$-equivariant ucp conditional expectation
        $$(A\rtimes \G)\rtimes \check{\G} \to \alpha^\rtimes(A\rtimes \G)$$
        which is to say that $(A\rtimes \G, \alpha^\rtimes)$ is inner $\check{\G}$-amenable.
        \item $(A, \alpha)$ being inner $\G$-amenable is equivalent to the existence of a $\G$-equivariant ucp conditional expectation
        $$\varphi: A\rtimes \G \to \alpha(A).$$
        This is equivalent with the existence of a $\check{\G}$-equivariant ucp conditional expectation
        \begin{equation}\label{100}(A\rtimes \G)\rtimes \G \to \alpha(A)\rtimes \G\cong A \rtimes \G.\end{equation}
        Note however that the map
        $$\Psi: (A\rtimes \G)\bar{\otimes}L^\infty(\check{\G})\to (A\rtimes \G)\rtimes \G: z \mapsto V_{23}zV_{23}^*$$ is a $\check{\G}$-equivariant $*$-isomorphism such that $\Psi\circ \alpha^\rtimes = \alpha\otimes \id$ on $A\rtimes \G$. Therefore, the existence of a $\check{\G}$-equivariant ucp map as in \eqref{100} is equivalent with the existence of a $\check{\G}$-equivariant ucp map
        $$\psi: (A\rtimes \G)\bar{\otimes}L^\infty(\check{\G})\to A\rtimes \G$$
        such that $\psi\circ \alpha^\rtimes = \id$, which is to say that $(A\rtimes \G, \alpha^\rtimes)$ is $\check{\G}$-amenable.
        \item Consider the coaction 
        $$\beta: B(L^2(A))\bar{\otimes} L^\infty(\G)\to (B(L^2(A))\bar{\otimes}L^\infty(\G))\bar{\otimes} L^\infty(\G): z \mapsto U_{\alpha,13}V_{23}z_{12}V_{23}^* U_{\alpha,13}^*.$$ If $A$ is $\G$-$W^*$-amenable, there is a $\G$-equivariant ucp map $$\varphi: (B(L^2(A))\bar{\otimes}L^\infty(\G), \beta)\to (\alpha(A), \id \otimes \Delta)$$ such that $\varphi(\pi_A(a)\otimes 1) = \alpha(a)$ for all $a\in A$. Define then
$$\widetilde{\varphi}: B(L^2(A))\bar{\otimes}L^\infty(\G)\to \alpha(A):z \mapsto \varphi(U_\alpha^*z U_\alpha).$$
Then it can be checked that $\widetilde{\varphi}: (B(L^2(A))\bar{\otimes}L^\infty(\G), \id\otimes \Delta)\to (\alpha(A), \id \otimes \Delta)$ is a $\G$-equivariant ucp map satisfying $\widetilde{\varphi}((\pi_A\otimes \id)\alpha(a)) = \alpha(a)$ for $a\in A$. Then the map
$$\widetilde{\varphi}\rtimes \G:  B(L^2(A))\bar{\otimes} B(L^2(\G)) \cong (B(L^2(A))\bar{\otimes}L^\infty(\G))\rtimes_{} \G\to \alpha(A)\rtimes_{} \G \cong A \rtimes \G$$
is a $\check{\G}$-equivariant ucp map such that $(\widetilde{\varphi}\rtimes \G)((\pi_A\otimes \id)(z)) = z$ for $z\in A\rtimes \G$. This means that $A\rtimes \G$ is $W^*$-amenable. Conversely, if $A\rtimes \G$ is $W^*$-amenable, then there exists a $\check{\G}$-equivariant ucp map $$\psi: (B(L^2(A))\bar{\otimes} B(L^2(\G)), \id \otimes \check{\Delta}_r)\to (A\rtimes \G, \alpha^\rtimes)$$ such that $\psi((\pi_A\otimes \id)(z)) = z$ for $z\in A\rtimes \G$. Taking $\check{\G}$-fixed points, we obtain a $\G$-equivariant ucp map $$\widetilde{\psi}: B(L^2(A))\bar{\otimes}L^\infty(\G)\to \alpha(A)$$ and then the map $$\varphi: (B(L^2(A))\bar{\otimes}L^\infty(\G), \beta)\to (\alpha(A), \id \otimes \Delta): z \mapsto \widetilde{\psi}(U_\alpha z U_\alpha^*)$$ is $\G$-equivariant and satisfies $\varphi(\pi_A(a)\otimes 1) = \alpha(a)$ for all $a\in A$, which is to say that $A$ is $\G$-$W^*$-amenable.
    \end{enumerate}

\section{Outlook}
 In this paper, the general theory of equivariant $W^*$-correspondences and their weak containment was developed. Many interesting approximation properties of $W^*$-dynamical systems $(A,\alpha)$ can then be formulated by considering which precise equivariant $A$-$A$-correspondences can weakly contain the trivial equivariant $A$-$A$-correspondence. In follow-up work of the authors, concrete applications of the theory developed in this paper will be investigated in the framework of compact and discrete quantum groups. We announce two such applications:
\begin{itemize}
    \item If $\G$ is a compact or a discrete quantum group and $A\subseteq B$ is a $\G$-equivariant inclusion (of $\sigma$-finite von Neumann algebras), then strong equivariant amenability and equivariant amenability of the inclusion $A\subseteq B$ are equivalent. 
    \item There exists an amenable discrete quantum group that acts non-amenably on a von Neumann algebra.
\end{itemize}
A further venue for research is the considerations of the topics of this paper in the setting of C$^*$-dynamical systems. For the moment, this remains completely unexplored in the setting of locally compact quantum groups.

\emph{Acknowledgements}: The work of K. De Commer was supported by Fonds voor Wetenschappelijk Onderzoek (FWO), grant G032919N. The work of J. De Ro is supported by the FWO Aspirant-fellowship, grant 1162522N. We thank J. Crann for an email exchange concerning the topics of this paper. We also thank the referees for valuable comments improving the exposition of the paper.

\end{document}